\documentclass[12pt]{amsart}
\usepackage[utf8]{inputenc}
\usepackage{enumerate}
\usepackage{graphicx,graphics}
\usepackage[usenames]{color}
\usepackage{amsfonts}
\usepackage{amssymb}
\usepackage{amsthm}
\usepackage{amsmath}
\usepackage{dirtytalk}
\usepackage{mathrsfs}  
\usepackage{lineno}
\usepackage[margin=1 in]{geometry}
\newtheorem{theorem}{Theorem}
\newtheorem{definition}{Definition}
\newtheorem{lemma}{Lemma}
\newtheorem{proposition}{Proposition}
\newtheorem{corollary}{Corollary}
\newtheorem{Problem}{Problem}

\newtheorem{problem}[theorem]{Problem}
\newtheorem{remark}{Remark}
 \usepackage{multicol}  
 \usepackage{bbm}
 \usepackage{tabularx}
\usepackage{fancyhdr}
\usepackage{enumitem}

\usepackage{lipsum}

\newcommand\blfootnote[1]{%
  \begingroup
  \renewcommand\thefootnote{}\footnote{#1}%
  \addtocounter{footnote}{-1}%
  \endgroup
}
\title{Construction schemes: Transferring structures from $\omega$ to $\omega_1$}
\author{Jorge Antonio Cruz Chapital}
\author{Osvaldo Guzm\'an Gonz\'alez}
\author{Stevo Todor\v{c}evi\'c}

\begin{document}
\maketitle
\begin{abstract} A structural analysis of construction schemes is developed. That analysis is used to give simple and new constructions of  combinatorial objects which have been of interest to set theorists and topologists.  We then continue the study of capturing axioms associated to construction schemes. From them, we deduce the existence of several uncountable structures which are known to be independent from  the usual axioms of Set Theory. Lastly, we prove that the capturing axiom $FCA(part)$ is implied by Jensen's $\Diamond$ principle.
\end{abstract}

 \blfootnote{The first and second authors were partially supported by a PAPIIT grant IN101323 and CONACyT grant A1-S-16164. The second author was also supported by PAPIIT grant IA102222. The third author is partially supported by grants from NSERC (455916), CNRS (IMJ-PRG-UMR7586) and SFRS (7750027-SMART).}The purpose of this paper is to continue developing a technique introduced by the third author in \cite{schemeseparablestructures} for constructing mathematical objects of cardinality  $\omega_1$, the first uncountable cardinal. Towards explaining the nature of this construction scheme, suppose we want to build a particular structure  of cardinality $\omega_1$. Probably the most straight forward approach is to build such object by countable approximations in $\omega_1$ many steps. Typical examples are the construction of an Aronszajn tree from \cite{Kunen} or \cite{Jech} and the construction of a Hausdorff gap from \cite{ScheepersGaps}, \cite{discovering2} or \cite{integersdouwen}\footnote{Of course, there are also other constructions of this objects that do not rely in a recursions of size $\omega_1$, see \cite{Walks}.}. The topic of the present work 
is on the method of \textit{construction and capturing schemes}. In our approach, the desired structure is not built from its countable substructures, but rather from its finite ones. Also, the construction takes only countable many steps. The approach is based on amalgamations of (many) isomorphic finite structures a method that has quite different flavour than the method used in typical constructions by countable approximations. To see this, the reader is invited to compare our constructions of, say,  an Aronszajn tree or a Hausdorff gap (see Theorem \ref{hausdorffgapconstruction}) with the  previous constructions mentioned above. For example, the classical construction of a Hausdorff gap  is usually regarded as difficult (for example, see \cite{gapsreconstructed}) while our construction is rather simple and takes only a few lines.\\

The idea of building an uncountable structure by finite approximations is not new. For example, this approach was previously performed with R. Jensen's gap $1$-morasses (see \cite{aspectsofconstructibility}) and its simplification by D. Velleman (see \cite{simplifiedmorasses}). Their motivation was to construct model-theoretic structures by amalgamating {\it two} finite structures at a given stage. This is greatly improved using construction schemes where many finite isomorphic structures are amalgamated at a given stage. Being able to amalgamate more than just two structures greatly improves the versatility of the method.\\

Roughly speaking, a construction scheme is a collection of finite sets with  strong coherent properties. These properties are allowing us recursively follow the construction scheme and  perform amalgamations of several (not just two) isomorphic finite structures. As in the case of Jensen's and Velleman's morasses, the existence of construction schemes on $\omega_1$ can be obtained just from $ZFC$, without the need to appeal to any extra axiom. However, sometimes in applications of our approach,  we need that the construction scheme  have further \say{capturing properties} which informally means that it can \say{capture} parts of an uncountable set  in \say{isomorphic positions}. This will be formally introduced and explained in Section \ref{beyondzfcsection}. For the existence of a \textit{capturing construction scheme} we need to go beyond $ZFC$. It is worth pointing out that the study of construction schemes was originally motivated by problems concerning metric structures such as normed spaces and Boolean algebras. The paper \cite{schemeseparablestructures} contains several applications of capturing construction schemes to functional analysis and topology. The general theory was further developed in \cite{forcingandconstructionschemes}, \cite{lopezschemethesis} and \cite{treesandgapsschemes}.\\

In the present paper, we construct a large amount of combinatorial objects that have been of interest to set theorists and topologists. The constructions developed here are rather different from the original constructions. It is often the case that using construction or capturing schemes greatly simplifies the previously  known constructions (for example, this is the case of Theorems \ref{luzinjonestheorem} and \ref{suslinlatticescheme}). In other cases, there was no direct construction previously known (see Theorem \ref{hausdorffcoherenttheorem}).\\

To follow this paper no previous knowledge of construction and capturing schemes is needed. In fact, we strive to make this a good introduction to this fascinating topic. It is our hope to spread the interest  towards readers  versed in amalgamation techniques in their own fields  giving them a powerful tool for building structures of their own interest.\\

In Section \ref{constructionschemesection} we introduce the notion of construction scheme and derive its basic properties. In Section \ref{ordinalmetricssection} we study the relation of construction schemes and the ordinal metrics studied in \cite{Walks}. In Section \ref{resultinzfcsection} we apply construction schemes to build a large number of structures that are of interest in set theory, topology and infinite combinatorics. The constructions are mostly independent from each other, so that the reader can start from the ones that she or he finds more interest. For the convencience of the reader, we list here the constructions that appear on this sections. References for further study and historic remarks will be provided as we encounter them on the Section \ref{resultinzfcsection}:\begin{enumerate}
    \item \textbf{Special Aronszajn tree.} Aronszajn trees are the most well-known examples of the incompactness of $\omega_1$. An Aronszajn tree is a tree of height $\omega_1$, its levels are countable, yet it has no cofinal branches. A simple way in which we can guarantee that a tree (of height $\omega_1$) has no cofinal branches is to make it \say{special}, which means that it can be covered with countably many antichains (equivalently, they can be embedded in the rational numbers). Aronszajn proved this kind of trees exist.
    
    \item \textbf{Countryman lines}. Let $(X,<)$ be a total linear order. Except for the trivial cases $X^2$ is not a linear order. In this way, it makes sense to ask how many chains we need to cover it. A countryman line is an uncountable linear order whose square can be covered with only countably many chains. These orders seem so paradoxical at first glance that Countryman conjectured they do not exist. However, it was first proved by Shelah that Countryman lines do exist (see \cite{shelahdecomposing}).
    
    \item \textbf{Hausdorff gaps}. An interesting feature of the Boolean algebra $\mathscr{P}(\omega)/fin$ is that it is not complete. The easiest way to see this is as follows: Let $\{A_\alpha\,|\,\alpha\in 2^\omega\}\subseteq [\omega]^{\omega}$ be an $AD$ family. For every $H\subseteq 2^\omega$, consider the set $\mathcal{A}_H=\{[A_\alpha]\,|\,\alpha\in H\}$ (where $[A_\alpha]$ is the class of $A_\alpha$ in $\mathscr{P}(\omega)/fin$). It follows by a simple counting argument that there must be an $H\subseteq 2^\omega$ such that $\mathcal{A}_H$ has no supremum. Although this is a very simple argument, we are often interested in more concrete examples of the incompleteness of $\mathscr{P}(\omega)/fin$. The nicest examples are provided by (Hausdorff) gaps. Moreover, gapes are important since they represent obstructions we may encounter when embedding structures on $\mathscr(\omega)/fin$ (a very illustrative example of this situation is Theorem 8.8 of \cite{PartitionProblems}). The classic construction of a Hausdorff gap requires a very clever argument to take care of $2^\omega$ many tasks in only $\omega_1$ many steps. With construction schemes, we can prove very easily that they exist.
    
    \item \textbf{Luzin-Jones almost disjoint families.} This topic is about almost disjoint families (AD) and their separation properties. A Luzin family is an $AD$ family of size $\omega_1$ in which no two uncountable subfamilies can be separated. On the other hand, a Jones family is an $AD$ family with the property that every countable subfamily of it can be separated from its complement. It is easy to prove that both of this kind of families exist. However, building a family that is both Luzin and Jones at the same time is much more complicated, since there is a tension between this two properties. A very difficult and highly complex construction of a Luzin-Jones family appears in \cite{guzman2019mathbb}. With construction schemes, we will be able to build a Luzin-Jones family very easily. It is worth pointing out that these families can be used to build interesting examples in functional analysis (see \cite{anonstable}).
    \item \textbf{Luzin coherent family of functions.} We now look at a generalization of the Hausdorff gaps discussed previously. A luzin coherent family of functions is a coherent system of functions supported by a pretower, in which we impose a strong non triviality condition. The importance of these families is that they provide many cohomologically different gaps. They were first studied by Talayco in \cite{talayco1995applications}. Later Farah proved that such families exist (see \cite{farah1996coherent}). The proof of Farah is highly non-constructive and indirect, since it appeals to Keisler's completeness Theorem. We will build such families using a construction scheme. No previous direct construction was known.
 \end{enumerate} 

In Section \ref{beyondzfcsection} we study distinct notions of capturing, which are additional strong properties we may demand from a construction scheme (We do this in an axiomatic way). The existence of such schemes can not be deduced from $ZFC$ alone. We then apply such schemes to build several combinatorial structures whose existence is known to be independent from $ZFC.$ Once again, we will provide a brief description of the structure that we will build using capturing schemes:
\begin{enumerate}
    \item\textbf{Suslin trees.} A Suslin tree is an Aronszajn tree in which every antichain is countable. The Suslin Hypothesis ($HS$) is the statement that there are no Suslin trees. We now know that $HS$ is independent from $ZFC$. A related concept, the Suslin lines, were introduced by Suslin while studying the ordering of the real numbers. Kurepa was the one to realize that there is a Suslin line if and only if there is a Suslin tree. Applications and constructions from Suslin trees are abundant in the literature. We will use capturing schemes to build two types of these trees: Coherent Suslin and full Suslin trees. These two families of trees are diametrically opposed. While forcing with a Coherent tree, it completely destroys the $ccc$-ness of it, while with a full Suslin tree, many subtrees of it remain $ccc$.
    \item \textbf{Suslin lower semi-lattices}. If in the definition of a Suslin tree, we relax the condition of being a tree to just a being a lower semi-lattice, we get the notion of a Suslin lower semi-lattice. They were introduced by Dilworth, Odell and Sari (see \cite{dilworth2007lattice}) in the context of Banach spaces. They were then studied by Raghavan and Yorioka (see \cite{raghavan2014suslin}). Among other things they proved that the $\Diamond$-principle implies that $\mathscr{P}(\omega)$ contains a Suslin lower semi-lattice. We were able to obtain the same result from a capturing scheme.
    \item \textbf{Entangled sets}. A well-known theorem of Cantor is that any two countable dense linear orders with no end points are isomorphic. The straight forward generalization to linear orders of size $\omega_1$ is false. For this reason we want to restrict to suborders of the real numbers. We say that  an uncountable $B\subseteq \mathbb{R}$ is $\omega_1$-dense if $|U\cap \mathbb{B}|=\omega_1$ whenever $U$ is an open interval whose intersection with $U$ is nonempty. A remarkable theorem of Baumgartner is that $PFA$ implies that any two $\omega_1$-dense sets of reals are isomorphic (see \cite{ApplicationsofPFA}). This statement is now known as Baumgartner axiom ($BA(\omega_1)$). An entangled set is a subset of $\mathbb{R}$ with very strong combinatorial properties. The existence of an entangled set implies the failure of the Baumgartner axiom. Entangled sets were introduced by Abraham and Shelah in order to show that $BA(\omega_1)$ does not follow from Martin axiom (see \cite{MAdoesnotImplyBA}). We will use a capturing scheme to build an entangled set. In this way, the existence of certain capturing schemes contradict the Baumgartner axiom. 
\item \textbf{Independent coherent family of functions}.
 We now return to the study of gaps that are obtained from a coherent family of functions, as we did in Theorem \ref{hausdorffcoherenttheorem}. However, this time we want our family of gaps to be \say{independent}. This means that we can either fill or freeze any subfamily without filling or freezing any of the remaining gaps in the family. A similar result was obtained by Yorioka assuming the $\Diamond$-principle in \cite{yoriokadestructiblegaps} (the analogue for Suslin trees was proved by Abraham and Shelah in \cite{Shelahabrahamincompactness}).
 \item \textbf{$ccc$ destructible $2$-bounded coloring without injective sets}. A coloring $c:[\omega_1]^2\longrightarrow\omega_1$ is called $2$-bounded if every color appears at most $2$ times. A set $A\subseteq \omega_1$ is $c$-injective if no color appears twice in $[A]^2$. Galvin was the first to wonder if there is a $2$-bounded coloring without an uncountable injective set. He proved that such coloring exists assuming the Countinuum Hypothesis. On the other hand, the third author proved that no such coloring exist under $PFA$.  years later, Abraham, Cummings and Smyth proved that $MA$ is consistent with the existence of a $2$-bounded coloring without uncountable injective sets (see \cite{abrahampolychromatic}). After hearing this result, Friedman asked for a concrete example of a $2$-bounded coloring without an uncountable injective set, but that such set can be added with a $ccc$ partial order. In \cite{abrahampolychromatic} such example is constructed assuming $CH$ and the failure of the Suslin hypothesis. We will also find an example with a capturing scheme. 
 \end{enumerate}

In the book \cite{PartitionProblems}, the third author developed an oscillation theory which is based on an unbounded family of functions. A plethora of applications of this theory has been found through the years. For example the oscillation theory is key in the proof that $PFA$ implies that the continuum is $\omega_2$ (see \cite{TopicsinSetTheory}). we will develop a similar theory using a capturing scheme. An important difference between the classic oscillation theory and the one from capturing schemes, is that this new one is based on a bounded family of functions. Using this new oscillation theory, we can prove the existence of the following objects:
\begin{enumerate}
    \item \textbf{Sixth Tukey type}. The Tukey ordering is a useful tool to compare directed partial orders. Its purpose is to study how a directed partial order behave cofinally. It was introduced by Tukey in \cite{Tukey} in order to study convergence in topology. The Tukey classification of countable directed partial orders is very simple: Every countable directed partial order is Tukey equivalent to $1$ or to $\omega$. The Tukey classification of directed sets of size $\omega_1$ becomes much more interesting. We now have at leas five Tukey types: $1$, $\omega$, $\omega_1$, $\omega\times\omega_1$ and $[\omega_1]^{<\omega}$. We may wonder if there is a directed partial order of size $\omega_1$ that is not Tukey equivalent to one if this five. In \cite{CategorycofinaltypesII}, Isbell proved that $CH$ entails the existence of a sixth Tukey type. This was greatly improved by the third author in \cite{directedsetscofinaltypes}, where he proved that $CH$ implies that there are $2^{\omega_1}$ distinct Tukey types. On the other hand, in the same paper he showed that $PFA$ implies that there are now sixth Tukey types. Here, we found a sixth Tukey type from a capturing scheme.
    \item $\omega_1\not\rightarrow (\omega_1,\omega+2)^2_2$. Given $\alpha<\omega_1$, the partition $\omega_1\rightarrow (\omega_1,\alpha)^2_2$ means that for every $c:[\omega_1]^2\longrightarrow 2$, either there is an uncountable $0$-monochromatic set, or there is a $1$-monochromatic set of order type $\alpha$. A celebrated result of Erd\"os and Rado (extending a theorem by Dushnik, Miller and Erd\"os) is that $\omega_1\rightarrow (\omega_1,\omega+1)^2_2$. We may wonder if this theorem can be improved. This turns out to be independent from $ZFC$. The Proper Forcing Axiom implies that $\omega_1\rightarrow (\omega_1,\alpha)^2_2$ for every $\alpha<\omega_1$ (see \cite{NotesonForcingAxioms}), while $\mathfrak{b}=\omega_1$ implies $\omega_1\not\rightarrow (\omega_1,\omega+2)^2_2$. we will prove a similar result from our oscillation theory. 
    \item \textbf{Non productivity of $ccc$ partial orders}. When is the product of two $ccc$ partial orders again $ccc$? This question has been of interest to set theorists for a long time. On one hand, Martin axiom implies that the product of $ccc$ partial orders is $ccc$. On the other hand, the failure of the Suslin Hypothesis implies the opposite. Consistent examples of two $ccc$ partial orders whose product is not $ccc$ have been constructed by Galvin under $CH$ (see \cite{Kunen}) and by the third author under $\mathfrak{b}=\omega_1$ in  \cite{PartitionProblems}. As an applications of the oscillation theory developed in the section, we encounter new situations in which there are $ccc$ partial orders whose product is not $ccc$. 
    \item \textbf{Suslin pretowers}. In \cite{oscillationsintegers}, the third author developed an analogue of his oscillation theory, now this time based on non-meager towers. While studying this oscillation, Borodulin-Nadzieja and Chodounsk\'y introduced the notion of a Suslin pretower. A pretower $\mathcal{T}=\{T_\alpha\,|\,\alpha\in\omega_1\}$ is Suslin if every uncountable subset of $\mathcal{T}$ contains two pairwise $\subseteq$-incomparable elements. The existence of a Suslin pretower is independent from $ZFC$. It can be proved that the Open Graph Axiom forbidds the existence of such families, while a Suslin pretower can be constructed assuming $\mathfrak{b}=\omega_1$ (see \cite{GapsandTowers}). We will build a Suslin pretower as an application of the oscillation theory obtained from a capturing scheme.
    \item \textbf{S-spaces}. Hereditarily separable and hereditarily Lindel\"of are two properties that in some sense are dual to each other. We may wonder if they always coincide. The question is only of interest in the realm of regular spaces. An $S$-space is a regular hereditarily separable, non Lindel\"of space. The study of $S$-spaces (and the dual notion, $L$-spaces) used to be one of the most active areas of set-theoretic topology. $S$-spaces can be constructed under several set-theoretic hypothesis (like $CH$ or the negation of the Suslin Hypothesis). It was the third author who for the first time succeeded in proving that it is consistent that there may be no $S$-spaces (see \cite{PartitionProblems}). Here, we will apply the oscillation theory obtained from a capturing scheme to construct several $S$-spaces.
    \end{enumerate}
The following table summarizes many of the constructions related to set theory and topology which can be carried out with a capturing scheme.\\\\

\begin{tabularx}{0.8 \textwidth} { 
  | >{\raggedright\arraybackslash}X 
  | >{\centering\arraybackslash}X 
  | >{\raggedleft\arraybackslash}X | }
 \hline
  \textbf{Result} & \textbf{Capturing Axiom }  \\
 \hline
 Sixth Tukey type  & $CA_2$   \\
\hline
$\omega_1\not\rightarrow(\omega_1,\omega+2)^2_2$ & $CA_2$ \\
\hline
Non productivity of $ccc$ & $CA_2$ \\
\hline
Suslin pretowers & $CA_2$ \\
\hline 
$S$-spaces & $CA_2$ \\
\hline
Suslin lower semi-lattices & $CA_2$ \\
\hline
Failure of $BA(\omega_1)$ & $CA_2$ \\
\hline
Destructible gap (see \cite{treesandgapsschemes}) & $CA_2(part)$ \\
\hline
$ccc$ destructible $2$-bounded coloring & $CA_3$ \\
\hline
Full Suslin Tree & $FCA$\\
\hline
Coherent Suslin Tree & $FCA(part)$\\
\hline
Entangled sets & $FCA$\\
\hline 
Independent coherent family of functions & $FCA$\\
\hline
\end{tabularx}\\\\

Finally, in Section \ref{diamondsection} we prove that  capturing construction schemes  exist assuming Jensen's $\Diamond$-principle. Originally, this appeared as claim in  \cite{schemeseparablestructures}, but the proof sketched there is not complete. Here we present a complete proof.
There is however still much to learn and understand on this topic.  Some of the open problems or directions for further research will be mentioned throughout the paper.
\section{Notation and preliminaries}
Our notation is fairly standard. Most of the concepts used will be defined in the sections where they first appear. The undefined notions related to set theory can be found in  \cite{Jech} or \cite{Kunen}. We will try our best to  point out more specific bibliography when it is necessary. Nevertheless, we will review some of the main notions used throughout the paper for the convenience of the reader. \\\\
The quantifiers $\forall^\infty$ and $\exists^\infty$ should be read as $\textit{for all but finitely many }$ and $\textit{ for infinitely many }$ respectively. The quantifier $\exists!$ stands for $\textit{there is a unique}$.\\\\
Given a set $X$, $\mathscr{P}(X)$ denotes its power set. Given a cardinal $\kappa$, $[X]^{\kappa}$ stand for the set of all subsets of $X$ of cardinality $\kappa$ and $[X]^{<\kappa}$ is the collection of all subsets of $X$ of cardinality strictly less than $\kappa.$ By $LIM$ we denote the set of all limit ordinals below $\omega_1$. For $A,B\in \mathscr{P}(\omega)$, we  write $A\subseteq^* B$ if $A\backslash B$ is finite and $A=^*B$ if $A\subseteq^* B$ and $B\subseteq ^*A.$ Following this notation, 
we write $A\subsetneq ^* B$ whenever $A\subseteq^*B$ and $A\not=^*B$. Analogously, if $f,g$ are two functions with domain $X$ and an ordered set as a codomain,  we write $f<^* g$ if $\{x\in X\,|\,g(x)\leq f(x)\}$ is finite and $f<g$ if the previous set is moreover empty.\\\\
A family $\mathcal{D}$ is a $\Delta$-system with root $R$ if $|\mathcal{D}|\geq 2$ and $A\cap B=R$ whenever $A,B\in \mathcal{D}$ are different. If $X$ is a set of ordinals, $A\in [X]^{<\omega}$ is nonempty and $n<|A|$, $A(n)$ stands for the nth element of $A$ with the usual ordering. Given $S\subseteq |A|$, $A[S]$ is defined as $\{A(n)\,|\,n\in S\}$. If $B$ is also a nonempty element of $[X]^{<\omega}$ and $\alpha<\beta$ for each $\alpha\in A$ and $\beta \in B$ we will write $A<B.$ Lastly, if $C,D\in[X]^{<\omega}$ we will write $C\sqsubseteq D$ if $C\subseteq D$ and whenever $c\in C$ and $d\in D$ is such that $d\leq c$, then $d\in C$.\\\\
Let $(X,<)$ be a partially ordered set. Given two elements of $X$, say $x$ and $y$, $x$ and $y$ are said to be comparable if either $x\leq y$ or $y<x$. If there is $z\in X$ such that $z\leq x$ and $z\leq y$, we say that $x$ and $y$ are compatible. Whenever $x$ and $y$ are not comparable we say that they are incomparable, and whenever they are not compatible then they are 
incompatible. By $(x,y)$ and $[x,y]$ we denote the sets $\{z\in X\,|\,x<z<y\}$ and $\{z\in X\,|\,x\leq z\leq y\}$ respectively. $(x,y]$ and $[x,y)$ are defined in an analogous way.  We say that $A\subseteq X$ is an interval if $(x,y)\subseteq A$ whenever $x,y\in A$. We say that $A$ is a chain if every two elements of $A$ are comparable. On the other side, if all the elements of $A$ are pairwise incomparable we call $A$ a pie. A similar notion which is more frequently used is the 
one of antichain. $A$ is said to be an antichain if every two elements of $A$ are incompatible.  $A$ is cofinal if every element of $X$ is below some element of $A$. We define the cofinality of $X$ as $cof(X)=\min\{|A|\,|\,A\subseteq X\textit{ is cofinal}\,\}$. Finally, $X$ is well founded if every nonempty subset of $X$ has a minimal element. In this case, there is a unique function $rank:X\longrightarrow Ord$ such that $rank(x)=\sup\{rank(y)+1\,|\,y<x\}$ for each $x\in X$. For each ordinal $\alpha$, let $X_\alpha=\{x\in X\,|\,rank(x)=\alpha\}$ and $Ht(X)=\min\{\alpha\in Ord\,|\, X_\alpha=\emptyset\}$.\\\\
As a final remark: Whenever $\alpha$ and $\beta$ are ordinals, $\beta\backslash \alpha=\{\xi\in \beta\,|\,\xi\geq \alpha\}=[\alpha,\beta)$. We will alternate between these two notations whenever we think one is more appropriate than the other. In particular, we may write that $k\in \omega\backslash n$ whenever we want to say that $k$ is a natural number greater or equal than $n.$ The other notation will be used whenever we want to remark that the respective set is an interval.
\section{Construction Schemes}\label{constructionschemesection}

\begin{definition}[Type]\label{definitiontype}We call a sequence $\tau=\{(m_k,n_{k+1},r_{k+1})\}_{k\in\omega}\subseteq \omega^3$ a type if:\begin{enumerate}[label=$(\alph*)$]
\item $m_0=1,$ 
\item $\forall k\in \omega\backslash 1 \big(\;n_k\geq 2\;\big),$
\item $\forall r\in \omega\; \exists^\infty k\in \omega\big(\; r_k=r\;\big), $
\item $\forall k\in \omega \big( \; m_k>r_{k+1}\; \big), $
\item $\forall k\in\omega\backslash 1 \big( \; m_{k+1}=r_{k+1}+(m_k-r_{k+1})n_{k+1}\;\big).$
\end{enumerate}
Additionally, we say that a partition of $\omega$, namely $\mathcal{P}$, is compatible with $\tau$ if: 
\begin{enumerate}
    \item[(c')]$\forall P\in \mathcal{P}\,\forall r\in\omega\,\exists^\infty k\in P (\;r_k=r\;).$
\end{enumerate}
\end{definition}
It is easy to see that points $(b)$ and $(e)$ of the previous definition imply that the sequence of  $m_k$'s is increasing.\\

 Observe that if $X$ is an infinite set and $\mathcal{F}\subseteq [X]^{<\omega}$ then $(\mathcal{F},\subsetneq)$ is well founded. Hence, for each ordinal $\alpha$ we can define $\mathcal{F}_\alpha$ as in the preliminaries . In this particular case we have that $\mathcal{F}=\bigcup\limits_{k\in\omega}\mathcal{F}_k$. Furthermore, if $\mathcal{F}$ is cofinal in $[X]^{<\omega}$ then $\mathcal{F}_k\not=\emptyset$ whenever $k\in\omega.$
\begin{definition}[Construction scheme,\,\cite{schemeseparablestructures}] \label{constructionscheme}Let $\{(m_k,n_{k+1},r_{k+1})\}_{k\in\omega}$ be a type,
and $X$ be a set of ordinals.
We say that  $\mathcal{F}\subseteq [X]^{<\omega}$ is a construction scheme  over
$X$ of type $\{(m_k,n_{k+1},r_{k+1})\}_{k\in\omega}$  if:\begin{enumerate}[label=$(\arabic*)$]
\item $\mathcal{F}$ is cofinal in $[X]^{<\omega},$
\item $\forall k\in\omega\;\forall F\in \mathcal{F}_k\big(\;|F|=m_k\;\big)$,
\item $\forall k\in\omega\;\forall F,E\in \mathcal{F}_k\big(\; E\cap F\sqsubseteq E,F\;\big)$,
\item $\forall k\in\omega\;\forall F\in \mathcal{F}_{k+1}\;\exists !F_0,\dots,F_{n_{k+1}-1}\in \mathcal{F}_k$ such that $$F=\bigcup\limits_{i<n_{k+1}}F_i.$$
Moreover, $\{ F_i\}_{i<n_{k+1}}$\footnote{We will call $\{F_i\}_{i<n_{k+1}}$ the canonical decomposition of $F.$} forms a $\Delta$-system with root $R(F)$ such that  $|R(F)|=r_{k+1}$ and $$R(F)<F_0\backslash R(F)<\dots < F_{n_{k+1}-1}\backslash R(F).$$
\end{enumerate}
\end{definition}

The following theorem was proved implicitly in \cite{schemeseparablestructures}. Another proof can be found in \cite{lopezschemethesis}.
\begin{theorem}[\cite{schemeseparablestructures}]\label{theoremexistenceconstructionschemes} For every type $\{(m_k,n_k,r_k)\}_{k\in\omega}$ there is a construction scheme over $\omega_1$ of that type.
\end{theorem}

It turns out that simplified morasses as defined in \cite{simplifiedmorasses} and 
construction schemes are  closely related. In 
fact, just as pointed out in \cite{lopezschemethesis} one might think of construction 
schemes as some sort of generalized simplified $(\omega_1,1)$-morasses. By this, we think it is a good idea to state their definition for the sake of completeness. Readers who want to 
know more about morasses and their applications may look at \cite{devlinconstructibility}, \cite{somebanachspacesaddedbyacohenreal}, \cite{morassesincombinatorialsettheory}, \cite{ablackboxtheoremformorasses}, \cite{widescatteredspacesandmorasses}, \cite{simplifiedmorasses},  \cite{souslintreesconstructedfrommorasses}, and \cite{morassesdiamondandforcing}.
\begin{definition}[simplified $(\omega,1)$-morass]Let $\mathcal{A}\subseteq[\omega_1]^{<\omega}$. We say that:\begin{itemize}
    \item $\mathcal{A}$ is homogeneous if for every $X,Y\in\mathcal{A}$, $rank(X)=rank(A)$ implies that $X$ and $Y$ have the same order type and $\{Z\in\mathcal{A}\,|\,Z\subseteq Y\}=f[\{Z\in\mathcal{A}\,|\,Z\subseteq X\}]$ where $f$ is the only increasing bijection from $X$ to $Y.$
    \item $\mathcal{A}$ is directed if for every $X,Y\in\mathcal{A}$ there is $Z\in\mathcal{A}$ such that $X\cup Y\subseteq Z.$
    \item $\mathcal{A}$ is locally almost directed if for every $X\in\mathcal{A}$ either $\{Z\in\mathcal{A}\,|\,Z\subsetneq X\}$ is directed or there are distinct $Y,Z\in\mathcal{A}$ such that $Y\cup Z\subseteq X$, $Y\cap Z\sqsubseteq Y, Z$ and $Y\backslash Z<Z\backslash Y.$
\end{itemize}
 Finally, we say that $\mathcal{A}$ is a simplified $(\omega,1)$-morass if it is homogeneous, directed, locally almost directed and $\bigcup\mathcal{A}=\omega_1.$   
\end{definition}
We finish this section by proving some basic facts about construction schemes that will be used throughout the paper. Although some of the proofs can be found in   \cite{lopezschemethesis} or \cite{schemeseparablestructures}, we decided to include them here for the convenience of the reader. As a direct consequence of these results, construction schemes over $\omega_1$ of type $\{(m_k,2,r_{k+1})\}_{k\in\omega}$ are, in particular, simplified $(\omega,1)$-morasses.\\\\
For the rest of this section, otherwise stated we fix a type $\{(m_k,n_k,r_{k+1})\}_{k\in\omega_1}$.
\begin{lemma}\label{firstfirstlemma}Let $X$ be a set of ordinals, $\mathcal{F}$ a construction scheme over $X$, $l\in \omega$ and $F\in \mathcal{F}_l$. For each $k\leq l$, it happens that $F=\bigcup\{G\in \mathcal{F}_k\,|\,G\subseteq F\}$.
\begin{proof}
The proof is done by induction over $l$ and is an immediate consequence of part $(4)$ in Definition 2.
\end{proof}
\end{lemma}
Since the sequence $\{m_k\}_{k\in\omega}$ is strictly increasing, the rank of each element of $\mathcal{F}$ is fully determined by its cardinality.
\begin{corollary}\label{corollaryinclusion}Let $X$ be an infinite set of ordinals, $\mathcal{F}$ a construction scheme over $X$ and $k\in\omega$. Then $X=\bigcup\mathcal{F}_k$.
\begin{proof}It is enough to prove that $X\subseteq \bigcup\mathcal{F}_k$. For this, let $\alpha\in X$. Since $X$ is infinite, we can take  a finite $A\subseteq X$ for which $\alpha\in A$ and $|A|\geq m_k$. By part (1) of Definition \ref{constructionscheme}, there is $F\in \mathcal{F}$ so that $A\subseteq F$. Observe that $|F|\geq m_k$. Thus, if $l\in\omega$ is such that $F\in \mathcal{F}_l$, then $l\geq k$. Using Lemma $\ref{firstfirstlemma}$ and the fact that $\alpha\in F$, we conclude  there is $G\in \mathcal{F}_k$ such that $\alpha\in G.$ 
 \end{proof}
\end{corollary}
An important feature of construction schemes is that there is a version of point (3) of Definition \ref{constructionscheme}, but where the elements of $\mathcal{F}$ have different rank. 
\begin{lemma}\label{interseccioninicial}Let $X$ be a set of ordinals, $\mathcal{F}$  a construction scheme over $X$, $k\leq l\in\omega$, $E\in \mathcal{F}_k$ and $F\in \mathcal{F}_l$:
\begin{enumerate}[label=$(\alph*)$]
\item  $E\cap F\sqsubseteq E.$ 
\item If $k<l$ and $E\subseteq F$, there is $i<n_l$ such that $E\subseteq F_i$.
\end{enumerate}
\begin{proof}First we prove part $($a$)$. For this, observe that as a consequence of Lemma \ref{firstfirstlemma} we have $E\cap F=\bigcup\{G\cap E\,|\,G\in\mathcal{F}_k\textit{ and }G\subseteq F\}$. In virtue of part (3) in Definition \ref{constructionscheme}, this equality tells us that $E\cap F$ is a union of initial segments of $E$. Thus, we are 
done.\\
Now, we prove part ($b$). For this, let $\alpha=\max E$.  We know  $E\subseteq F$, so there is $i<n_k$ for which $\alpha\in F_i$. By part ($a$), $\alpha\in E\cap F_i\sqsubseteq E$. This implies $E\cap F_i= E$, so $E\subseteq F_i \subseteq F$.
\end{proof}
\end{lemma}
Note however, that if $k<l$, $F\in \mathcal{F}_l$ and $E\in \mathcal{F}_k$, in general is not true that $E\cap F$ is an initial segment of $F$. This can be easily seen as follows: Take any $F\in \mathcal{F}_{l+1}$ (For any $l\in\omega)$ and let $\{F_i\}_{i<n_{l+1}}$ be its canonical decomposition. Since $n_{l+1}\geq 2$, then $F_1\cap F$ is not an initial segment of $F$.
\begin{lemma}\label{equalconstructions}Let $X,Y$ be infinite sets of ordinals and $\mathcal{F}$, $\mathcal{G}$ be constructions schemes over $X$ and $Y$ respectively. If $X\subseteq Y$ and $\mathcal{F}\subseteq \mathcal{G}$ then $\mathcal{F}=\{G\in\mathcal{G}\,|\,G\subseteq X\}.$
\begin{proof}As an opening observation note that if $\mathcal{F}\subseteq \mathcal{G}$, then they have the same type. Under this hypothesis, it is trivial that $\mathcal{F}\subseteq\{G\in\mathcal{G}\,|\,G\subseteq X\}$. For the other inclusion, let $G\in\mathcal{G}$ be such that $G\subseteq X$. Also let $\alpha=\max G$ and $k\in\omega$ for which $G\in\mathcal{G}_k$. By Corollary \ref{corollaryinclusion}, there is $F\in\mathcal{F}_k$  having $\alpha$ as an element. Since $F\in\mathcal{F}\subseteq\mathcal{G}$, we know  $\alpha\in F\cap G\sqsubseteq G$. By maximality of $\alpha$ in $G$ this means $F\cap G= G$. In other words, $G\subseteq F.$ Since $F$ and $G$ are finite sets of cardinality $m_k$, it follows that $F=G$. In this way, we conclude $G\in \mathcal{F}$.

\end{proof}
\end{lemma}
\begin{proposition}\label{constructionschemefinitesets}Let $X$ be a finite set of ordinals. There is a construction scheme over $X$ if and only if $|X|=m_k$ for some $k\in\omega$. Moreover, if there is a construction scheme over $X$, it is unique.
\begin{proof}For the \textit{only if} part, suppose  there is a construction scheme $\mathcal{F}$ over $X.$ By part (1) of Definition \ref{constructionscheme}, we know  there is $F\in\mathcal{F}$ for which $X\subseteq F.$ This implies  $X=F$, so $|X|=|F|=m_k$ for some $k\in\omega.$\\
We prove the \textit{if} part and uniqueness by induction over $k$.\\\\
{(Base step)} If $k=0$, then $X=\{\alpha\}$. It should be clear that $\mathcal{F}=\{X\}$ is the only construction scheme over $X.$\\\\
{(Induction step)} Consider $k\in\omega\backslash 1$, and suppose we have proved the proposition 
for each $i<k$ and for every finite set of ordinals $Y$ satisfying $|Y|=m_i$. Let $X$ be a 
set of ordinals of cardinality $m_k.$ By part $(e)$ in Definition \ref{definitiontype}, the 
set $$X_i=X[r_k]\cup X\big[[(m_{k-1}-r_k)\cdot(i),(m_{k-1}-r_k)\cdot(i+1))\big]$$ is well defined for 
each $i<n_k$ and has cardinality $m_{k-1}$. The induction hypotheses assure that for each 
$i<n_k$ there is a unique construction scheme $\mathcal{F}^i $ over $X_i$. We claim  $$\mathcal{F}=\bigcup\limits_{i<n_k} \mathcal{F}^i\cup\{X\}$$ is a construction scheme over $X$.\\ Trivially $\mathcal{F}$ satisfies conditions (1), (2)  and (4) of Definition \ref{constructionscheme}. In order to prove condition (3), take $l\in\omega$ and $F,G\in \mathcal{F}_l$. Since $|X|=m_k$ and $F\subseteq X$,  $l\leq k$. If $l=k$ then $F$ and $G$ must be equal to $X$. Thus, we may suppose without loss of generality that $l<k$. In this way, there are $i,j<n_k$ for which $F\in \mathcal{F}^i$ and $G\in \mathcal{F}^j$. If $i=j$ we conclude by recalling that $\mathcal{F}^i$ is a construction scheme. In the case where $i\not=j$, consider $\phi:X_i\longrightarrow X_j$ the only increasing bijection. Observe that $\{\phi[H]\,|\,H\in \mathcal{F}^i\}$ is a construction scheme over $X_j$, so by induction hypotheses this set is equal to $\mathcal{F}^j$. Furthermore, $\phi[F]\in \mathcal{F}^j_l$ and $\phi$ is the identity map when restricted to $X[r_k]$. In this way, $F\cap X[r_k]=\phi(F)\cap X[r_k]$. This means $F\cap G=X[r_k]\cap F\cap G=X[r_k]\cap\phi[F]\cap G$. But $\phi[F]\cap G$ is an initial segment of $G$ and $X[r_k]$ is an initial segment of $X_j$. Thus, their intersection (which equals $F\cap G$) is an initial segment of $G$.\\
In order to prove the uniqueness of $\mathcal{F}$, take $\mathcal{F}'$ an arbitrary construction scheme over $X$. Using part $(4)$ of Definition \ref{constructionscheme} and part $(b)$ of Lemma \ref{interseccioninicial}, we know $$\mathcal{F'}=\bigcup\limits_{i<n_k}\{G\in\mathcal{F}\,|\,G\subseteq X_i\}\cup\{X\}.$$ Observe that $\{G\in\mathcal{F}\,|\,G\subseteq X_i\}$ is a construction scheme over $X_i$ for each $i<n_k$. Since there is a unique construction scheme over $X_i$, this set must be equal to  $\mathcal{F}^i$. This observation finishes the proof.

\end{proof}
\end{proposition}
\begin{definition}Let $X$ be a set of ordinals of cardinality $m_k$ for some $k\in\omega$. We will call $\mathcal{F}(X)$ the only construction scheme over $X.$ \end{definition}
In particular, if $X$ is an arbitrary set of ordinals, $\mathcal{F}$ is a construction scheme over $X$ and $F\in\mathcal{F}$, then $\mathcal{F}(F)=\{G\in\mathcal{F}\,|\,G\subseteq F\}$.
\begin{corollary}\label{corollaryisomorphismfiniteschemes}Let $X,Y$ be finite sets of ordinals of the same cardinality. If
$\mathcal{F}$ and $\mathcal{G}$ are construction schemes over $X$ and $Y$ respectively and
$h:X\longrightarrow Y$ is the only increasing bijection, then $$\mathcal{G}=\{h[F]\,|\,F\in\mathcal{F}\}.$$
\begin{proof}
For this, just note that $\{h[F]\,|\,F\in\mathcal{F}\}$ is a construction scheme over $Y.$ By the uniqueness part of Lemma \ref{constructionschemefinitesets}, we are done.
\end{proof}
\end{corollary}
As a particular instance of the Corollary \ref{corollaryisomorphismfiniteschemes} and the observation previous to it, we have that $\{G\in \mathcal{F}\,|\, G\subseteq F\}=\{F[G]\,|\,G\in \mathcal{F}(m_k)\}$ whenever $F\in\mathcal{F}_k$. This fact will be implicitly used throughout the paper.
We remark that, in general, infinite sets of ordinals may have more than one construction scheme, or none at all.

\section{Construction schemes, ordinal metrics and related functions}\label{ordinalmetricssection}
In this section we study the connection between ordinal metrics and construction schemes. Ordinal metrics were introduced by the third author, and together with the walks on ordinals, they have proven to be invaluable tools in the study of $\omega_1$. The book \cite{Walks} contains a considerable amount of information about them, as well as their relation to walks on ordinals.

\begin{definition}Let $X$ be a set of ordinals and $\rho:X^2\longrightarrow \omega$ be an arbitrary function. For every $F\in[X]^{<\omega}\backslash\{\emptyset\}$ and $k\in\omega$ we define the $k$-closure of $F$ as:\\$$(F)_k=\{\alpha\in X\;|\;\exists \beta\in F\backslash\alpha\,(\,\rho(\alpha,\beta)\leq k\,)\,\}$$
$$(F)^{-}_k=(F)_k\cap \max(F)=(F)_k\backslash\{\max(F)\}.$$\\
Additionally, we say define the diameter of $F$ as $$\rho^F=\max\{\rho(\alpha,\beta)\,|\,\alpha,\beta\in F\}.$$ $F$ is said to be $k$-closed whenever $F=(F)_k.$ Moreover, if $F$ is $k$-closed and $k=\rho^F$ we simply say that $F$ is closed. Whenever $\alpha\in \omega_1$, we will write $(\{\alpha\})_k$ and $(\{\alpha\})^-_k$ simply as $(\alpha)_k$ and $(\alpha)^-_k$ respectively.
\end{definition}
\begin{remark}Given $\rho:X^2\longrightarrow \omega$ and $\beta\in X$, we have that $(\beta)_k=\{\alpha\leq \beta\;|\;\rho(\alpha,\beta)\leq k\}.$
\end{remark} 

\begin{definition}[Ordinal metric]\label{ordinalmetricdefinition}Let $X$ be a set of ordinals. We say that $\rho:X^2\longrightarrow \omega$ is an ordinal metric if \begin{enumerate}[label=(\alph*)]
\item $\forall\alpha,\beta\in X\,(\,\rho(\alpha,\beta)=0\textit{ if and only if }\alpha=\beta\,),$
\item $\forall \alpha,\beta\in X\,(\,\rho(\alpha,\beta)=\rho(\beta,\alpha)\,),$
  \item $\forall\alpha,\beta,\gamma\in X\,(\;\alpha<\beta,\gamma \rightarrow \rho(\alpha,\beta)\leq \max(\rho(\alpha,\gamma),\rho(\beta,\gamma))\;).$
\item $\forall \beta\in X\,\forall k\in\omega\,(\;|(\beta)_k|<\omega\;),$
 \end{enumerate}
\end{definition}
 As the reader may note, the previous definition reassembles that of an ultrametric\footnote{A metric $d$ over a space $X$ is an ultrametric if $d(x,y)\leq \max(d(x,z),d(z,y))$ for all $x,y,z\in X$.}. However, note that one of the instances of the triangle inequality is missing. In this sense, we can interpret $(\beta)_k$ as the ball with radius $k$ centered at $\beta$ intersected with $\beta$. In the theory of metric spaces, we usually want to learn about the nature of points in a given space. To do that, we study the behaviour of their basic neighbourhoods as they become smaller. Ordinal metrics tend to work in the opposite way. Here, usually we want to construct a structure whose elements are parametrized with the points of our space. In order to do that, we make suitable approximations of our structure by analyzing the interaction between $\beta$ and other points in $(\beta)_k$ as $k$ grows larger.\\
 The following Lemma is a direct consequence of the definitions and the triangle inequalities. 
 \begin{lemma}\label{basicpropertiesmetric}Let $\rho:X^2\longrightarrow \omega$ be an ordinal metric over a set of ordinals $X$. Given $A\in[X]^{<\omega}\backslash\{\emptyset\}$, consider $k=\rho^A$ and $\gamma=\max A$. \begin{enumerate}
\item $\rho(\alpha,\gamma)\leq k$ for each $\alpha\in (A)_k$,
\item $(A)_k=(\gamma)_k$,
\item $(A)_k$ is closed and $\rho^{(A)_k}=k$
\end{enumerate}
    
\end{lemma} 
Now, we consider some notions associated to an ordinal metric. In these definitions, whenever we say that $F$ is a maximal closed set, we mean that $F$ is closed and there is no closed set $G$ such that $\rho^F=\rho^G$ and $F\subsetneq G.$

\begin{definition}Let $\rho:X^2\longrightarrow \omega$ be an ordinal metric. We will say that $\rho $ is:\begin{enumerate}
\item locally finite if  $\sup\{\,|F|\,|\,F\in[X]^{<\omega}\backslash\{\emptyset\},\,F\textit{ is closed and }\rho^F\leq k\}<\omega$ for every $k\in\omega.$
\item homogeneous if whenever $F,G\in[X]^{<\omega}\backslash\{\emptyset\}$ are maximal closed sets such that $\rho^F=\rho^G$, then $|F|=|G|$ and if $h:F\longrightarrow G$ is the only increasing bijection from $F$ to $G$ then  $\rho(h(\alpha),h(\beta))=\rho(\alpha,\beta)$ for each $\alpha,\beta\in F.$
    \item regular if for each $0<k\in\omega$ and each maximal closed $F\in [X]^{<\omega}\backslash\{\emptyset\}$ with $\rho^F=k$ there are $j_F\in \omega\backslash 1$ $F_0,\dots,F_{j_F}\in[X]^{<\omega}\backslash\{\emptyset\}$such that\begin{itemize}
    \item $F=\bigcup\limits_{i<j_F}F_i,$
        \item For each $i<{j_F}$, $F_i$ is a maximal closed set with $\rho^{F_i}=k-1,$
        \item $\{F_i\}_{i<{j_F}}$ forms a $\Delta$-system with root $R(F)$ such that 
         $$R(F)<F_0\backslash R(F)<\dots < F_{j_F-1}\backslash R(F).$$ 
    \end{itemize}
\end{enumerate}
\end{definition}
\begin{proposition}\label{ordinalmetricscheme}Let $X$ be an infinite set of ordinals of order type greater than $\omega$ and let $\rho:X^2\longrightarrow \omega$ be a locally finite, homogeneous, 
and regular ordinal metric. The set $$\mathcal{F}^\rho=\bigcup\{F\in[X]^{<\omega}\backslash\{\emptyset\}\,|\,F\textit{ is maximal closed}\}$$
is a construction scheme over $X$.
\begin{proof}
First we prove that $\mathcal{F}=\mathcal{F}^\rho$ satisfies part (1) of Definition  \ref{constructionscheme}. In order to achieve this, take an arbitrary $A\in [X]^{<\omega}$ and let $k=\rho^A$. As stated in Lemma \ref{basicpropertiesmetric}, $(A)_k$ is closed and $\rho^{(A)_k}=k.$ Since $\rho$ is locally finite there must be a closed $F\in[X]^{<\omega}$ with $\rho^F=k$ which contains $(A)_k$, and such that $|F|=\max\{|G|\,|\,G\in[X]^{<\omega}\backslash\{\emptyset\}, G\textit{  is closed },\rho^G\leq k\textit{ and }(A)_k\subseteq G\}$. This last property implies  $F$ is a maximal closed set, so we are done.\\
We claim that $\{F\in \mathcal{F}\,|,\rho^F=k\}\not=\emptyset$ for each $k\in\omega$ and this set coincides with the set of elements of $\mathcal{F}$ with rank $k$ with respect to the inclusion, 
namely $\mathcal{F}_k.$ To prove this, let $\beta\in X$ be such that $\beta\cap X$ is infinite. There is such $\beta$ because the order type of $X$ is greater than $\omega$. By part (d) of Definition \ref{ordinalmetricdefinition}, we know there is $\alpha\in \beta\cap X$ for which $\rho(\alpha,\beta)>k$. Consider $F\in \mathcal{F}$ satisfying $\{\alpha,\beta\}\subseteq F$. Observe that $F^\rho\geq \rho(\alpha,\beta)>k.$ By repeatedly using regularity of $\rho$ we can 
find a maximal closed $G\in[X]^{<\omega}\backslash\{\emptyset\}$ contained in $F$ for which $\rho^G=k$. This proves the first part of the claim. In fact, the previous argument shows that for each $F\in \mathcal{F}$ and every $k<F^\rho$ there is $G\in \mathcal{F}$ contained in $F$ for which $\rho^G=k$. The second part of the claim follows directly from this fact. \\ The homogeneity of $\rho$ implies $|F|=|G|=m_k$ for each $F,G\in \mathcal{F}_k$ and part (1) of Definition \ref{ordinalmetricdefinition} implies  $|F|=1$ for each $F\in \mathcal{F}_0$. In particular, this proves part $(2)$ of Definition \ref{constructionscheme}.\\
For $k\in \omega\backslash 1$ and $F,G\in \mathcal{F}_k$, let $j_F,j_G\in \omega\backslash 1$ and $F_0,\dots,F_{j_F},G_0,\dots,F_{j_G}\in [X]^{<\omega}\backslash\{\emptyset\}$ testify regularity of $\rho$ for $F$ and $G$ respectively. Using the homogeneity of $\rho$, it can be shown that $j_F=j_G$ and if $h:F\longrightarrow G$ is the only increasing bijection, then $G_i=h[F_i]$ for each $i\leq j_F$. So in particular, $h[R(F)]=R(G)$. In this way, if we define $r_k=|R(G)|$ and $n_k=j_G+1$ it is not hard to see that part (4) of Definition \ref{constructionscheme} is satisfied.\\
It only remains to prove part (3) of Definition \ref{constructionscheme}. For this, let $k\in \omega$ and $F,G\in \mathcal{F}_k$. Let $\beta \in F\cap G$ and $\alpha\in G$ be such that $\alpha<\beta$. Since $\alpha,\beta\in G$ and $\rho^G=k$, we know  $\rho(\alpha,\beta)\leq k$. But $F=(F)_k$ and $\beta \in F$, which means that $\alpha \in F$ . Since $\alpha$ was arbitrary, we conclude that $F\cap G\sqsubseteq G$.
 \end{proof}
\end{proposition}
 It turns out that every construction scheme naturally generates a locally finite, homogeneous and regular ordinal metric. 
\begin{definition}Let $\mathcal{F}$ be a construction scheme over a set of ordinals $X$. We define \break $\rho_\mathcal{F}:X^2\longrightarrow \omega$ as:
$$\rho(\alpha,\beta)=\min\{k\in\omega\;|\;\exists F\in\mathcal{F}_k(\;\{\alpha,\beta\}\subseteq F\;)\}.$$
If $\mathcal{F}$ is clear from context, we will write $\rho_\mathcal{F}$ simply as $\rho.$
\end{definition}
Note that $\rho$ is well-defined since $\mathcal{F}$ is cofinal in $[X]^{<\omega}$
\begin{lemma}\label{lemmaclosure}Let $\mathcal{F}$ be a construction scheme over a set of ordinals $X$. For each $\beta\in X$, $k\in \omega$ and $F\in \mathcal{F}_k$ satisfying $\beta\in F$, it happens that $(\beta)_k=F\cap (\beta+1)$ and $(\beta)^-_k=F\cap \beta.$
\begin{proof} By definition of $\rho$,  $F\cap(\beta+1)\subseteq (\beta)_k$ and $(\beta)_k\subseteq \beta+1$. It only remains to prove that $(\beta)_k\subseteq F$. For this, let $\alpha\in (\beta)_k$ and $G\in \mathcal{F}_{\rho(\alpha,\beta)}$ with $\{\alpha,\beta\}\subseteq G$. By definition of $k$-closure, we know $\rho(\alpha,\beta)\leq k$. As a consequence of this, we can use Lemma \ref{interseccioninicial} to conclude $G\cap F\sqsubseteq G$. Since $\beta\in G\cap F$, this means $\alpha\in G\cap F\subseteq F$, so we are done.
\end{proof}
\end{lemma}
The previous lemma will be frequently used throughout the paper without an explicit reference to it. 
\begin{proposition}\label{schemeordinalmetric} Let $\mathcal{F}$ be a construction scheme over a set of ordinals $X$. Then $\rho_\mathcal{F}$ is a locally finite, homogeneous and regular ordinal metric.
\begin{proof}First we prove  $\rho=\rho_\mathcal{F}$ is an ordinal metric. Trivially $\rho$ satisfies conditions (a) and  (b) of Definition \ref{ordinalmetricdefinition}, and (d) is a direct consequence of Lemma \ref{lemmaclosure}.
We will divide the proof of $(c)$ into two cases.\begin{itemize}
    \item [($\beta<\gamma$)] Let $E\in \mathcal{F}_{\rho(\alpha,\gamma)}$ and 
    $F\in \mathcal{F}_{\rho(\beta,\gamma)}$ which contain $\{\alpha,\gamma\}$ 
    and $\{\beta,\gamma\}$ respectively. By Lemma \ref{interseccioninicial}, either $E\cap F\sqsubseteq E$ or $E\cap F\sqsubseteq F$. Since $\gamma\in E\cap F$, this implies that $\alpha\in F $ or $\beta\in E$, so $\rho(\alpha,\beta)\leq \rho(\beta,\gamma)$ 
    or $\rho(\alpha,\beta)\leq \rho(\alpha,\gamma)$.\\
    \item [($\gamma<\beta$)]Let $F\in \mathcal{F}_{\rho(\alpha,\beta)}$, $E\in \mathcal{F}_{\rho(\alpha,\gamma)}$ and $G\in \mathcal{F}_{\rho(\beta,\gamma)}$ which contain $\{\alpha,\beta\}$, 
    $\{\alpha,\gamma\}$ and $\{\beta,\gamma\}$ respectively. As we proved in Lemma \ref{lemmaclosure}, if 
    $\gamma\notin F$ then $\rho(\gamma,\beta)> \rho(\alpha,\beta)$, so we can suppose without loss of generality that 
    $\gamma\in F$. Now, take $\{F_i\}_{i<n_{\rho(\alpha,\beta)}}$ the 
    canonical decomposition of $F$ and  $j<n_{\rho(\alpha,\beta)}$ such 
    that $\beta\in F_j$. Suppose towards a contradiction that 
    $\rho(\alpha,\beta)> \max(\rho(\alpha,\gamma),\rho(\beta,\gamma)$. In this case we can make use of Lemma \ref{interseccioninicial} to conclude 
    $\beta\in F_j\cap G\sqsubseteq G$. As $\gamma<\beta$, we have that $\gamma\in F_j\cap G$ which implies  $\gamma\in F_j\cap E\sqsubseteq E$. Thus $\alpha\in F_j$, but this means $\{\alpha,\beta\}\subseteq F_j$. This is a contradiction to the the minimality of $F$, so we are done.
    \end{itemize}
 Now we prove  $\rho$ is locally finite, homogeneous and regular. First we claim that for each $k\in\omega$ $$\mathcal{F}_k=\{F\in[X]^{<\omega}\backslash\{\emptyset\}\,|\,F\textit{ is maximal closed with }\rho^F=k\}.$$ \\
To prove this, let $G\in[X]^{<\omega}$ be a closed set with $\rho^G=l$ and let $\beta=\max G.$ Take $F\in\mathcal{F}_l$  such that $\beta\in F$. For each $\alpha\in G$, we know  there is $k\leq l$ and  $H\in \mathcal{F}_k$ containing $\{\alpha,\beta\}$.  We know that $\beta\in F\cap H\sqsubseteq H$, which implies  $\alpha\in F\cap H\subseteq F. $ This proves  $G\subseteq F$. Since $G$ was arbitrary, we are done. \\
From this equality it is obvious that $\rho$ is locally finite. Regularity and homogeneity follow directly from part $(4)$ of Definition \ref{constructionscheme} and Corollary \ref{corollaryisomorphismfiniteschemes}.
\end{proof}
\end{proposition}
A careful analysis of the proofs of Propositions \ref{ordinalmetricscheme} and \ref{schemeordinalmetric} yields to the following corollary.
\begin{corollary}Let $X$ be an infinite set of ordinals with order type greater than $\omega$. If $\mathcal{F}$ is a construction scheme over $X$ and $\rho:X^2\longrightarrow \omega$ is a locally finite, homogeneous and regular ordinal metric, then:
\begin{itemize}
    \item $\mathcal{F}^{\rho_\mathcal{F}}=\mathcal{F}$,
    \item $\rho_{\mathcal{F}^\rho}=\rho.$
\end{itemize}
\end{corollary}

We end this subsection by remarking that the missing triangle inequality from Definition \ref{ordinalmetricdefinition} is never consistent for ordinal metrics over $\omega_1$.
\begin{proposition}No ordinal metric on $\omega_1$ is an ultrametric.
\begin{proof}  
 Suppose $\rho$ is as previously mentioned, and let $M$ be a countable elementary submodel of $H(\omega_2)$ having $\rho$ as an element. Let $\gamma\in \omega_1\backslash M $ and consider an arbitrary $\alpha\in M\cap \omega_1$. Let $k$ be $\rho(\alpha,\gamma).$ Using elementarity we can find $\beta\in M\cap \omega_1$ above $(\gamma)_{k}\cap M$ for which $k=\rho(\alpha,\beta)$. By definition of the $k$-closure, we have that $\rho(\beta,\gamma)>k=\max(\rho(\alpha,\beta),\rho(\alpha,\gamma))$. 
 \end{proof}
 \end{proposition}
\subsection{More functions associated to a construction scheme}
Apart from the ordinal metric associated to a construction scheme, there are two important functions that we need to analyze before we enter to the world of applications. The first one being the $\Delta$ function and the second one being the $\Xi$ function. For the rest of this subsection we fix a set of ordinals $X.$
\begin{definition}\label{Deltafunction}Let $\mathcal{F}$ be a construction scheme over $X$. We define $\Delta: \omega_1^2\longrightarrow \omega+1$ as: $$\Delta(\alpha,\beta)=\min\{k\in\omega\,|\,|(\alpha)_k|\not=|(\beta)_k|\,\}\cup\{\omega\}$$
\end{definition}
\begin{lemma}\label{lemmmma}Let $\alpha,\beta\in X$ be distinct ordinals and $k\in\omega$. If $|(\alpha)_k|=|(\beta)_k|$ then $k<\Delta(\alpha,\beta)$.
\begin{proof}Let $F,G\in \mathcal{F}_k$ such that $\alpha\in F$ and $\beta\in G$ and let $h:F\longrightarrow G$ be the increasing bijection. Using the homogeneity of $\rho$, we conclude that $h[(\alpha)_l]=(h(\alpha))_l$ for each $l\leq k$. By hypothesis $h(\alpha)=\beta$, so we are done.
\end{proof}
\end{lemma}
\begin{remark}\label{remarkdelta} By definition of $\rho$, we have that $(\alpha)_{\rho(\alpha,\beta)}\subsetneq  (\beta)_{\rho(\alpha,\beta)}$ whenever $\alpha<\beta$, which implies that $\Delta(\alpha,\beta)\leq \rho(\alpha,\beta)$. In particular, $\Delta(\alpha,\beta)=\omega$ if and only if $\alpha=\beta$.
\end{remark}

\begin{lemma}\label{countrymanlemma3} Let $\alpha,\beta,\delta\in\omega_1$ be distinct ordinals such that $\Delta(\alpha,\beta)< \Delta(\beta,\delta)$. Then $\Delta(\alpha,\delta)=\Delta(\alpha,\beta)$.
\begin{proof} Define $k$ as $\Delta(\alpha,\beta)-1$. By the minimality condition on the definition of $\Delta$, we have that $|(\alpha)_k|=|(\beta)_k|$. Observe that $k$ is strictly smaller than $\Delta(\beta,\delta)$. In this way,  $|(\beta)_k|=|(\delta)_k|$. Consequently,  $\Delta(\alpha,\delta)\geq k+1=\Delta(\alpha,\beta)$. On the other hand, $|(\beta)_{k+1}|=|(\delta)_{k+1}|$. This means that $|(\alpha)_{k+1}|\not=|(\delta)_{k+1}|$, which implies $\Delta(\alpha,\delta)\leq \Delta(\alpha,\beta).$
\end{proof}
\end{lemma}

 \begin{lemma}\label{countrymanlemma1}Let $\mathcal{F}$ be a construction scheme over $X$, $\{\alpha,\beta\}\in [X]^2$, $k=\Delta(\alpha,\beta)-1$ and $\phi$ the increasing bijection from $(\alpha)_{k}$ to $(\beta)_{k}$. Suppose that $\delta,\gamma\in (\alpha)_k$ are such that $\delta<\gamma$ and $\phi(\delta)\not=\delta$. Then the following happens: \begin{enumerate}[label=$(\alph*)$]\item$\phi(\gamma)\not=\gamma,$
\item$\rho(\alpha,\beta)\geq \Delta(\delta,\phi(\delta))$
\item $\Delta(\delta,\phi(\delta))\geq \Delta(\gamma,\phi(\gamma)),$
\item $\Delta(\gamma,\phi(\gamma))\geq \Delta(\alpha,\beta).$
\end{enumerate}
\begin{proof} To prove (a) we first observe $\delta\notin (\alpha)_k\cap(\beta)_k.$ This follows from the fact that $\phi(\delta)\not=\delta$. Since $(\alpha)_k\cap(\beta)_k\sqsubseteq (\alpha)_k$ and $\gamma>\delta$, it is also true that $\gamma\notin(\alpha)_k\cap (\beta)_k$. This last statement is the same as saying $\phi(\gamma)\not=\gamma$.\\
    Now, we prove (b). For this, let $F\in\mathcal{F}_{\rho(\alpha,\beta)}$ for which $\{\alpha,\beta\}\subseteq F.$ By definition of $k$, $(\alpha)_k\cup(\beta)_s\subseteq F$. This means both $\delta$ and $\phi(\delta)$ belong to $F.$ In this way
$(\delta)^-_{\rho(\alpha,\beta)}=\delta\cap F\not=\phi(\delta)\cap 
    F=(\phi(\delta))^-_{\rho(\alpha,\beta)},$ and consequently $\rho(\alpha,\beta)\geq 
    \Delta(\delta,\phi(\delta))$.\\ To prove $(d)$ observe that $\phi[(\gamma)_k]=(\phi(\gamma))_k$. This implies
    $\Delta(\gamma,\phi(\gamma))\geq k+1=\Delta(\alpha,\beta)$.\\
    Finally we prove $(c)$. By hypothesis $\gamma\in (\alpha)_k$, which implies $(\gamma)_k=(\alpha)_k\cap(\gamma+1)$. Furthermore, by point $(d)$ it is also true that $(\gamma)_k\subseteq (\gamma)_{r}$ where $r=\Delta(\gamma,\phi(\gamma))-1$. As a consequence of these facts (and using that $\delta<\gamma$) we  conclude $\delta\in (\gamma)_{r}$. In a similar way, we can deduce that $\phi(\delta)\in (\phi(\gamma))_{r}$. To finish, let $\psi$ be the increasing bijection  from $(\gamma)_{r}$ to $(\phi(\gamma))_{r}$. The homogeneity of $\rho$ implies  $\psi|_{(\gamma)_k}=\phi|_{(\gamma)_k}$. In particular  $\psi(\delta)=\phi(\delta)$, which means  $\psi[(\delta)_r]=(\phi(\delta))_r$. It follows that  $\Delta(\delta,\phi(\delta))\geq r+1=\Delta(\gamma,\phi(\gamma))$, so the proof is over.
\end{proof}
\end{lemma}

\begin{definition}\label{Xifunction} Let $\mathcal{F}$ be a construction scheme over $X$ and $\alpha\in X$. $\Xi_\alpha:\omega\longrightarrow \omega\cup\{-1\}$ is the function defined as:
$$\Xi_\alpha(k)=\begin{cases}0 &\textit{if }k=0\\
-1 &\textit{if }k>0\textit{ and }|(\alpha)_k|\leq r_k\\
\frac{|(\alpha)_k|-|(\alpha)_{k-1}|}{m_{k-1}-r_k} &\textit{otherwise}
\end{cases}$$
\end{definition}
As we shall see in the next proposition, the function $\Xi_\alpha$ tells us in which piece of a canonical decomposition lies $\alpha.$
\begin{proposition}\label{propositionxifunction}Let $\mathcal{F}$ be a construction scheme over $X$, $\alpha\in X$ and $k\in \omega\backslash 1$. If $F\in \mathcal{F}_k$ is such that $\alpha\in F$ and $\{F_i\}_{i<n_k}$ is its canonical decomposition, then:\begin{enumerate}[label=(\alph*)]
    \item $\Xi_\alpha(k)=-1$ if and only if $\alpha\in R(F),$
    \item If $i<n_k$ is such that $\alpha\in F_i\backslash R(F)$ then $\Xi_\alpha(k)=i.$
\end{enumerate}
\begin{proof}(a) is clear by definition of $\Xi$.  To prove (b), suppose that $i$ is as stated. Since $F=R(F)\cup\big(\bigcup\limits_{j<n_k}F_j\big)$ and $R(F)<F_0\backslash R(F)<\dots<F_{n_k-1}\backslash R(F)$, we conclude that
\begin{equation*}\begin{split}(\alpha)_k=F\cap(\alpha+1)&=R(F)\cup\big(\bigcup\limits_{j<i}(F_j\backslash(R(F)\big)\cup (F_i\backslash R(F)\cap(\alpha+1))\\ &=(F_i\cap(\alpha+1))\cup\big(\bigcup\limits_{j<i}F_j\big).
\end{split}
\end{equation*}
The sets appearing just before the last equality are disjoint. This means that $$|(\alpha)_k|-|(\alpha)_{k-1}|=\big(\sum\limits_{j<i}(m_{k-1}-r_k)\big)=i(m_{k-1}-r_k).$$
Dividing by $m_{k-1}-r_k$, we obtain the equality we were looking for.
\end{proof}
\end{proposition}
\begin{lemma}\label{lemmaxi}Let $\mathcal{F}$ be a construction scheme over $X$ , $\alpha<\beta\in X$ and $k\in \omega\backslash 1$. Then:
\begin{enumerate}[label=$(\alph*)$]
\item If $k<\Delta(\alpha,\beta)$, then  $\Xi_\alpha(k)=\Xi_\beta(k).$
\item If $k=\rho(\alpha,\beta)$, then $0\leq \Xi_\alpha(k)<\Xi_\beta(k).$
\item If $k>\rho(\alpha,\beta)$, then either $\Xi_\alpha(k)=-1$ or $\Xi_\alpha(k)=\Xi_\beta(k).$ 
\end{enumerate}
\begin{proof}(a) is just a consequence of the definition of $\Delta(\alpha,\beta)$. To prove (b), fix $F\in \mathcal{F}_{\rho(\alpha,\beta)}$ for which $\{\alpha,\beta\}\subseteq F$. We know there are $i,j<n_k$ such that $\alpha\in F_i$ and $\beta\in F_j$. Since $F_i,F_j\in \mathcal{F}_{\rho(\alpha,\beta)-1}$, minimality of $\rho$ implies that $\alpha\notin F_j$ and $\beta\notin F_i$. In particular, this means $i\not=j $ and $\alpha,\beta\not\in R(F)$. Thus, by Proposition \ref{propositionxifunction} we have that $\Xi_\alpha(k)=i$ and $\Xi_\beta(k)=j$. Moreover, since $\alpha<\beta$ then $i<j$. \\
Now, we prove $(c)$. For this, suppose that $\Xi_\alpha(k)\not=-1$ and let $F\in \mathcal{F}_k$ be such that $\{\alpha,\beta\}\subseteq F$. Since $\alpha<\beta$, we also have $\Xi_\beta(k)\not=-1$. Thus, there is $j<n_k$ satisfying $\beta\in F_j\backslash R(F) $. Notice that $\alpha\in (\beta)_{k-1}=(\beta+1)\cap F_j$ and $\alpha\not\in R(F)$. This means  $\Xi_\alpha(k)=j=\Xi_\beta(k)$, so we are done.
\end{proof}
\end{lemma}
\begin{lemma}\label{lemmaxidelta}Let $\mathcal{F}$ be a construction scheme over $X$,  $\{\alpha,\beta\}\in[X]^2$  and $k=\Delta(\alpha,\beta)$.
Then $\Xi_\alpha(k)$ and $\Xi_\beta(k)$ are both distinct and greater or equal than $0.$
\begin{proof}
Let $F,G\in \mathcal{F}_k$ be such that $\alpha\in F$ and $\beta\in G$. Also, let $\phi:F\longrightarrow G$ be the increasing bijection. Observe that $\Delta(\alpha,\beta)=\Delta(\phi(\alpha),\beta)\leq \rho(\phi(\alpha),\beta)\leq k $. Thus, we can use part $(c)$ of Lemma \ref{lemmaxi} to conclude that $\Xi_{\phi(\alpha)}(k)$ and $\Xi_\beta(k)$ are both distinct and greater or equal than $0$. To finish, just notice that $\Xi_\alpha(k)=\Xi_{\phi(\alpha)}(k).$
\end{proof}
\end{lemma}
In virtue of the previous Lemma, we can improve the part $(d)$ of Lemma \ref{countrymanlemma1}. 
\begin{lemma}\label{countrymanlemma5}Let $\mathcal{F}$ be a construction scheme over $X$, $\{\alpha,\beta\}\in [X]^2$, $\Delta=\Delta(\alpha,\beta)$ and $\phi$ be the increasing bijection from $(\alpha)_{\Delta-1}$ to $(\beta)_{\Delta-1}$. If $\gamma\in (\alpha)_{\Delta-1}$ is such that $\phi(\gamma)\not=\gamma$ then $\Delta(\gamma,\phi(\gamma))>\Delta$ if and only if $\Xi_\gamma(\Delta)=-1$. Furthermore, if $\Xi_\gamma(\Delta)\geq 0$ then $\Xi_\gamma(\Delta)=\Xi_\alpha(\Delta)$ and $\Xi_{\phi(\gamma)}(\Delta)=\Xi_{\beta}(\Delta).$
\begin{proof}
By point $(d)$ of Lemma \ref{countrymanlemma1} we always have that $\Delta(\gamma,\phi(\gamma))\geq \Delta(\alpha,\beta)$. If $\Xi_\gamma(\Delta)=-1$ then the conclusion of part $(c)$ in Lemma \ref{lemmaxidelta} can not hold when applied to $\gamma$ and $\phi(\gamma)$. In virtue of this, $\Delta$ must be distinct from $\Delta(\gamma,\phi(\gamma))$.\\
Now suppose $\Xi_\gamma(\Delta)\geq 0$. By hypothesis we have that $\rho(\alpha,\gamma)\leq \Delta-1$. In this way, we can use part $(c)$ of Lemma \ref{lemmaxi} to conclude that $\Xi_\alpha(\Delta)=\Xi_\gamma(\Delta)$. $\Xi_\beta(\Delta)$ is equal to $\Xi_{\phi(\gamma)}(\Delta)$ for the same reason. By applying once more Lemma \ref{lemmaxidelta} to $\alpha$ and $\beta$ we get that $\Xi_\gamma(\Delta)$ is distinct from $\Xi_{\phi(\gamma)}(\Delta)$. Thus, $\Delta(\gamma,\phi(\gamma))\leq k+1=\Delta(\alpha,\beta)$. 
\end{proof}
\end{lemma}
\begin{corollary}\label{corollaryxi}Let $\mathcal{F}$ be a construction scheme over $\omega_1$, $\xi<\alpha<\beta\in\omega_1$. If $\rho(\xi,\beta)<\rho(\alpha,\beta)$, then $\rho(\xi,\alpha)<\rho(\alpha,\beta)$.
\begin{proof}
Since $\rho$ is an ordinal metric, $\rho(\xi,\alpha)\leq \max\{\rho(\xi,\beta),\rho(\alpha,\beta)\}=\rho(\alpha,\beta)$. Suppose towards a contradiction that $\rho(\alpha,\beta)=\rho(\xi,\alpha)$. Note that $0\leq \Xi_\xi(\rho(\alpha,\beta))<\Xi_\alpha(\rho(\alpha,\beta))<\Xi_\beta(\rho(\alpha,\beta))$ follows from part (b) of Lemma \ref{lemmaxi}. On the other hand, since $\rho(\alpha,\beta)>\rho(\beta,\xi)$, part (c) of Lemma \ref{lemmaxi} implies that either $\Xi_\xi(\rho(\alpha,\beta))=-1$ or $\Xi_\xi(\rho(\alpha,\beta))=\Xi_\beta(\rho(\alpha,\beta))$. Both cases contradict the previous inequality, so we are done.
\end{proof}
\end{corollary}
\section{Theorems in ZFC}\label{resultinzfcsection}
In this section, we will study a variety of uncountable structures in the framework of $ZFC$. Although most of these objects are well studied and several ways to prove their existence are known, we aim to show the value of construction schemes by presenting new constructions of this structures, some of which are considerably shorter than their classical counterparts.
\subsection{A Luzin-Jones family} For a countably infinite set $X$, we say that a family $\mathcal{A}\subseteq [X]^\omega$ is an almost disjoint family if $A\cap B=^*\emptyset$ whenever $A,B\in \mathcal{A}$ are different. Almost disjoint families are one of the central objects of study in modern combinatorial set theory.  Constructing  almost disjoint families with special properties is usually difficult and in most cases had lead to the development of powerful tools (see  \cite{InvariancePropertiesofAlmostDisjointFamilies}, \cite{Madfamilieswithstrongcombinatorialproperties}, \cite{madnessandnormality},
\cite{frechetlike},  \cite{nonpartitionable}, \cite{cechfunction}, \cite{SANEPlayer} and \cite{ThereisavanDouwenMADfamily}). Almost disjoint families have also played a central roll in the solution of many problems of Topology and Analysis. An example of this is the solution of the selection problem posed by J. van Mill and E. Watell in \cite{van1981selections} and solved by M. Hru\v{s}\'ak and I. Mart\'inez-Ru\'iz in \cite{hruvsak2009selections}. The reader interested in learning more about almost disjoint families is refered to \cite{TopologyofMrowkaIsbellSpaces}, \cite{AlmostDisjointFamiliesandTopology}, and \cite{Combinatoricsoffiltersandideals}.
\begin{definition}Let $\mathcal{A}=\{ A_\alpha\}_{\alpha\in\omega_1}$ be an almost disjoint family. We say that:\begin{itemize}
    \item $\mathcal{A}$ is Luzin if $\{\alpha\in\beta\,|\,|A_\alpha\cap A_\beta|\leq n\}$ is finite for each $\beta\in\omega_1$ and $n\in\omega$.
    \item $\mathcal{A}$ is Jones if for each $\beta\in\omega_1$ there is $C\in[\omega]^{\omega}$ such that $A_\alpha\subseteq^* C$ for every $\alpha\leq \beta$ and $A_\gamma\cap C=^*\emptyset$ for each $\gamma>\beta.$
    \item $\mathcal{A}$ is Luzin-Jones if it is both Luzin and Jones.
\end{itemize} 
\end{definition}
    The first construction of a Luzin family was done in \cite{luzin1947subsets} by N. N. Luzin. A Jones family was implicitly constructed by F. B. Jones in \cite{jones1937concerning}. Although, at first glance, Luzin and Jones properties seem to be incompatible, a construction of a Luzin-Jones family was obtained in \cite{guzman2019mathbb} by the second author, M. Hru\v{s}\'ak and P. Koszmider (building from work by Koszmider in \cite{onconstructionswith2cardinals}). We would like to point out that the highly complex construction of the mentioned family is carried out through the use of simplified morasses. Here, we use construction schemes to give an elementary construction of a Luzin-Jones Family. The reader may want to compare both proofs in order to gain insight between the similarity and differences between construction schemes and simplified morasses.
\begin{theorem}\label{luzinjonestheorem}There is a Luzin-Jones family.
\begin{proof}Let $\mathcal{F}$ be a construction scheme of type $\{(m_k,2,r_{k+1})\}_{k\in\omega}$. For every $k\in \omega\backslash 1$, let $$N_k= \{k\}\times m_{k-1}\times \big((m_{k-1}-r_k)k\big).$$
Our goal is to construct a Luzin-Jones family $\{A_\alpha\}_{\alpha\in \omega_1}$ over the union of the $N_k$'s, which we call $N$. For this, take an arbitrary $\alpha\in\omega_1$. Given $k\in \omega$, we define $A_\alpha^k\subseteq N_k$ as follows:
\begin{enumerate}
     \item If $\Xi_\alpha(k)\leq 0$, let $A^k_\alpha=\{k\}\times \{|(\alpha)^-_k|\}\times (m_{k-1}-r_k)k$.
     \item If $\Xi_\alpha(k)=1$, let $A^k_\alpha=\{k\}\times [r_k,m_{k-1})\times [(|(\alpha)^-_k|-m_{k-1})k,\,(|(\alpha)^-_k|-m_{k-1}+1)k).$
\end{enumerate}     
Now, we define $A_\alpha$ as $\bigcup\limits_{k\in\omega\backslash\{0\}}A^k_\alpha$ and $\mathcal{A}$ as $\{A_\alpha\}_{\alpha\in\omega_1}$. \\
First we show that $\mathcal{A}$ is almost disjoint. For this purpose, let $\alpha<\beta\in\omega_1$ and $k\in\omega\backslash 1$ be such that $\alpha\in (\beta)^-_{k}$. By definition of $\rho$,  $\rho(\alpha,\beta)\leq k<k+1$. Thus,  by Lemma \ref{lemmaxi}, either $\Xi_\alpha(k+1)=-1$ or $\Xi_\alpha(k+1)=\Xi_\beta(k+1)$. From this fact, it is straightforward that $A^{k+1}_\alpha\cap A^{k+1}_\beta =\emptyset.$  The previous argument implies $$A_\alpha\cap A_\beta=\bigcup_{k\in\omega\backslash\{\emptyset\}}\big(A^k_\alpha\cap A^k_\beta\big)\subseteq \bigcup_{k\leq \rho(\alpha,\beta)}U_{k+1}.$$

To prove that $\mathcal{A}$ is Luzin, we will first argue that if $\alpha<\beta$ and $k=\rho(\alpha,\beta),$ then $|A_\alpha\cap A_\beta|\geq k$. From this, it will follow that $\{\alpha<\beta\,|\,A_\alpha\cap A_\beta<k\}\subseteq (\beta)_k$. So let $\alpha, \beta$ and $k$ be as previously stated. Just observe that Lemma \ref{lemmaxi} assures  $\Xi_\alpha(k)=0<\Xi_\beta(k)=1$. In this way, $A^k_\alpha\cap A^k_\beta=\{k\}\times \{|(\alpha)^-_k|\}\times [(|(\alpha)^-_k|-m_{k-1})k,\,(|(\alpha)^-_k|-m_{k-1}+1)k).$ Trivially the cardinality of this set is $k$ and it is contained in $A_\alpha\cap A_\beta$. Thus, we are done. \\

We end the proof by showing that $\mathcal{A}$ is a Jones family. For this purpose define  $C_\beta$ as $\bigcup\limits_{k\in\omega\backslash 1}\big(\bigcup\limits_{\alpha\in (\beta)_k } A^k_\alpha\big)$ for each $\beta\in\omega$ and notice that for every $\alpha\in\omega_1$, if $\alpha\leq \beta$ then $A_\alpha\subseteq^* C_\beta$ and if $\alpha>\beta$ then $A_\alpha\cap C_\beta=^*\emptyset.$
\end{proof}
\end{theorem}
\subsection{Countryman lines and Aronszajn trees} Countryman lines are a certain kind of linear orders whose existence was proposed by R. S. Countryman in the 1970's. In \cite{shelahdecomposing}, S. Shelah proved their existence from ZFC for the first time. Years later, an easy construction was found by the third author using walks on ordinals (see \cite{Walks}). Readers interested in learning more about Countryman lines  may search for \cite{hudson2007canonical}, \cite{moore2006five}, or \cite{Walks}. On the other hand, Aronszajn trees are a particular kind of trees constructed by N. Aronszajn in 1930's. Their existence imply that the natural generalization of K\"{o}nig's Lemma to $\omega_1$ is false. For further information about Aronszajn trees, we recommend to search for \cite{someremarksoncoherent},  \cite{Kunen}, \cite{minimalscatteredorders}, \cite{gapstructure}, \cite{structuralanalysisofaronzajn}, \cite{universalaronszajn}, \cite{coherenttreesthatare}, \cite{Walks} or \cite{todorvcevic1984trees}.
\begin{definition}[Countryman line]We say that a totally ordered set $(X,<)$ is a countryman line if $X$ is uncountable and $X^2$ can be covered by countably many chains\footnote{Here, we consider the order over $X^2$ given by $(x,y)\leq (w,z)$ if and only $x\leq w$ and $y\leq z$.}.   
\end{definition}
It is interesting to compare our construction to the one in \cite{Walks}. In our case, the difficulty is mainly in proving that our order is transitive, while decomposing the square into countably many chains is very easy. In the one in \cite{Walks} it is easy to prove that the order is transitive, while the main difficulty is in the decomposition part. It seems both constructions are of a different nature.

\begin{definition}Let $\mathcal{F}$ be a construction scheme. For each $\{\alpha,\beta\}\in[\omega_1]^2$ recursively decide whether $\alpha<_\mathcal{F} \beta$ when one of the following conditions holds:\begin{enumerate}[label=$(\alph*)$]
    \item $|(\alpha)_{\Delta}\cap(\beta)_{\Delta}|\geq r_{\Delta}$ and $|(\alpha)_{\Delta}|< |(\beta)_{\Delta}|$ $($or equivalently $\Xi_{\alpha}(\Delta)<\Xi_{\beta}(\Delta)$$)$\footnote{The equivalence of this two conditions is an easy consequence of Lemma \ref{lemmaxidelta} and the part $(a)$ in Lemma \ref{lemmaxi}.}.
   
    \item  $|(\alpha)_{\Delta}\cap(\beta)_{\Delta}|< r_{\Delta}$ and $\min(\alpha)_{\Delta}\backslash(\beta)_{\Delta}<_\mathcal{F} \min\,(\beta)_{\Delta}\backslash(\alpha)_{\Delta}$.
\end{enumerate}
here $\Delta=\Delta(\alpha,\beta).$
\end{definition}
Paraphrasing, take $E,F\in \mathcal{F}_\Delta$ with $\alpha\in E$
 and $\beta\in F$. We know $E\cap F$ is an initial segment of both. This intersection may or may not contain $R(E)$ (respectively $R(F)$). In case $R(E)\subseteq E\cap F$ (which implies $R(E)=R(F)$), we make $\alpha<_\mathcal{F}\beta$ if $|E\cap \alpha|<|F\cap \beta|$. On the other side, if $R(E)\not\subseteq E\cap F$ then we relate $\alpha$ and $\beta$ in the same way as $\min(E\backslash F)$ and $\min(F\backslash E).$
 \begin{lemma}\label{countrymanlemma2}Let $\{\alpha,\beta\}\in [\omega_1]^2$, $\Delta=\Delta(\alpha,\beta)$ and $\phi$ the increasing biyection from $(\alpha)_{\Delta-1}$ to $(\beta)_{\Delta-1}$. Given $\gamma\in (\alpha)_{\Delta-1}$ for which $\phi(\gamma)\not=\gamma$, the following statements are equivalent:
\begin{enumerate}
    \item $\alpha<_\mathcal{F} \beta$,
    \item $\gamma<_\mathcal{F} \phi(\gamma)$.
\end{enumerate}
\begin{proof}We prove both implications simultaneously. The proof is carried by induction. Let $\alpha,\beta, \Delta,\phi$ and $\gamma$ be as in the hypotheses and suppose that we have proved the lemma for each $\alpha'<\alpha$ and $\beta'<\beta$. Before starting observe that $(\alpha)_\Delta\cap (\beta)_{\Delta}=(\gamma)_{\Delta}\cap (\phi(\gamma))_{\Delta}$. This is because $\phi(\gamma)\not=\gamma.$\\
 In the case $\Xi_{\gamma}(\Delta)\geq 0$, we can summon Lemma \ref{countrymanlemma5} to conclude $\Delta(\gamma,\phi(\gamma))=\Delta$, $\Xi_\gamma(\Delta)=\Xi_{\alpha}(\Delta)$ and $\Xi_{\beta}(\Delta)=\Xi_{\phi(\gamma)}(\Delta)$. The desired equivalence follows by this fact and the opening observation. 
 On the other side, if $\Xi_\gamma(\Delta)=-1$ it is necessarily true that $|(\alpha)_\Delta\cap (\beta)_\Delta|<r_\Delta$. Due to Lemma \ref{countrymanlemma5}, we also have that $k=\Delta(\gamma,\phi(\gamma))>\Delta$. Let $\delta=\min(\alpha)_\Delta\backslash (\beta)_\Delta$. Then $\phi(\delta)=\min(\beta)_{\Delta}\backslash (\alpha)_{\Delta}$ and $\delta\in (\gamma)_{\Delta}\subseteq (\gamma)_{k-1}$.  Let $\psi:(\gamma)_{k-1}\longrightarrow (\phi(\gamma))_{k-1}$ be the increasing bijection. The homogeneity of $\rho$ guarantees that $\phi(\delta)=\psi(\delta)$. By making use of the induction hypotheses we conclude that $\gamma<_\mathcal{F} \phi(\gamma)$ if and only if $\delta<_\mathcal{F} \phi(\delta).$ By definition, and since $|(\alpha)_\delta\cap (\beta)_\delta|<r_\Delta$, we also have that $\alpha<_\mathcal{F} \beta$ if  and only if $\delta<_\mathcal{F} \phi(\delta)$. This completes the proof.
\end{proof} \end{lemma}
As a direct consequence of the homogeneity of $\rho$, we can generalize Lemma \ref{countrymanlemma2} by means of the following corollary.
\begin{corollary}\label{countrymancorollary}Let $\{\alpha,\beta\}\in [\omega_1]^2$, $z<\Delta(\alpha,\beta)$ and $\phi$ the increasing bijection from $(\alpha)_{z}$ to $(\beta)_{z}$. 
Given $\gamma\in (\alpha)_{z}$ for which $\phi(\gamma)\not=\gamma$, the following statements are equivalent:
\begin{enumerate}
\item $\alpha<_\mathcal{F} \beta$,
    \item $\gamma<_\mathcal{F} \phi(\gamma)$.
\end{enumerate}
\begin{proof}Let $\psi:(\alpha)_{\Delta(\alpha,\beta)-1}\longrightarrow(\beta)_{\Delta(\alpha,\beta)-1}$ be the increasing bijection. To prove the corollary, just notice that $\psi|_{(\alpha)_z}=\phi.$

\end{proof}
\end{corollary}

\begin{proposition}$(\omega_1,<_\mathcal{F})$ is a total order.
\begin{proof}The only nontrivial part of this task is to prove transitivity. For this, take 
distinct $\alpha,\beta,\delta\in\omega_1$. We will show by induction 
that neither $\alpha<_\mathcal{F}\beta<_\mathcal{F}\delta<_\mathcal{F}\alpha$ nor $\alpha<_\mathcal{F}\delta<_\mathcal{F}\beta<_\mathcal{F}\alpha$. Suppose we 
have proved that fact for every $\alpha'<\alpha,$  $\beta'<\beta$ and $\delta'<\delta$. 
Let $s_1=\Delta(\alpha,\beta)$, $s_2=\Delta(\alpha,\delta)$ and $s_3=\Delta(\beta,\delta)$. For $x,y\in \{\alpha,\beta,\delta\}$, let $c(x,y)=\min(x)_{\Delta(x,y)}\backslash(y)_{\Delta(x,y)} $. By Lemma \ref{countrymanlemma3}, we can suppose without loss of generality either $s_1=s_2=s_3$ or $s_1=s_2<s_3$. Both cases are handled in a similar way so we will only prove the latter, and to do this we will consider the following subcases:
\begin{enumerate}
    \item If $|(\alpha)_{s_1}\cap(\beta)_{s_1}|\geq r_{s_1}$ and $|(\alpha)_{s_2}\cap(\delta)_{s_2}|\geq r_{s_2}$, we use the fact that since $s_3>s_1$, then $|(\beta)_{s_1}|=|(\delta)_{s_1}|$. From this, either  $\alpha<_\mathcal{F}\beta,\delta$ or $\beta,\delta<_\mathcal{F} \alpha$.\\

    \item If $|(\alpha)_{s_1}\cap(\beta)_{s_1}|< r_{s_1}$ and $|(\alpha)_{s_2}\cap(\delta)_{s_2}|\geq r_{s_2}$ it is not hard to see that $(\alpha)_{s_1}\cap(\beta)_{s_1}=(\beta)_{s_1}\cap(\delta)_{s_1}$. By definition of $<_\mathcal{F}$, $\alpha$ and $\beta$ are related the same way as $c(\alpha,\beta)$ and $c(\beta,\alpha)$. Corollary \ref{countrymancorollary} implies $\delta$ and $\beta$ are related in the same way as $c(\beta,\alpha)$ and $\min(\delta)_{s_1}\backslash\big((\beta)_{s_1}\cap(\delta)_{s_1}\big)$. By induction hypothesis, we are done.\\
    
    \item If $|(\alpha)_{s_1}\cap(\beta)_{s_1}|<|(\alpha)_{s_2}\cap(\delta)_{s_2}|< r_{s_1}$, just as in the previous case it will happen that  $(\alpha)_{s_1}\cap(\beta)_{s_1}=(\beta)_{s_1}\cap(\delta)_{s_1}$. Take $\phi$ the increasing bijection from $(\delta)_{s_1}$ to $(\beta)_{s_1}$. By definition, $\alpha$ and $\beta$ are related in the same manner as $c(\alpha,\beta)$ and $c(\beta,\alpha)$. Furthermore, $c(\alpha,\beta)\in(\delta)_{s_1}$ and $\phi(c(\alpha,\beta))=c(\beta,\alpha)\not=c(\alpha,\beta)$. By Corollary \ref{countrymancorollary}, $\delta$ and $\beta$ are related the same way as $c(\alpha,\beta)$ and $c(\beta,\alpha)$. This means either $\alpha,\delta<_\mathcal{F} \beta$ or $\beta<_\mathcal{F}\alpha,\delta.$\\
 \item If $|(\alpha)_{s_1}\cap(\beta)_{s_1}|=|(\alpha)_{s_2}\cap(\delta)_{s_2}|< r_{s_1}$,  $\alpha$ and $\beta$ are related the same way as $c(\alpha,\beta)$ and $c(\beta,\alpha)$. Furthermore, $\alpha$ and $\delta$ are related the same way as $c(\alpha,\delta)=c(\alpha,\beta)$ and $c(\delta,\alpha)$. To finish, just use the induction hypothesis.
\end{enumerate}
The remaining subcases are equivalent to one of the previous, so we are done.

\end{proof}
\end{proposition}
\begin{theorem}[Countryman line]\label{countrymanconstruction} $(\omega_1,<_\mathcal{F})$ is a Countryman line.
\begin{proof}
The only thing left to do is to prove that $\omega_1^2$ can by partitioned into $\omega$ chains. For this purpose, it is sufficient to prove the same holds for $D=\{(\alpha,\beta)\in\omega_1^2\;|\;\alpha<\beta\}$. For every $x,y,z\in\omega$ define $$P(x,y,z)=\{(\alpha,\beta)\in D\;|\; |(\alpha)_z|=x,\;|(\beta)_z|=y \textit{ and }\rho(\alpha,\beta)=z\}.$$
It is clear that the set of all $P(x,y,z)$ is countable and covers $D$, so we will show each $P(x,y,z)$ is a chain. For this purpose take $(\alpha,\beta),(\delta,\gamma)\in P(x,y,z)$ and suppose  $\alpha<_\mathcal{F}\delta$. Since $|(\beta)_z|=|(\gamma)_z|$ we get that $\Delta(\beta,\gamma)>z$.  Now, let $\phi$ 
be the increasing bijection from $(\beta)_{z}$ to $(\gamma)_{z}$. Then $\phi(\alpha)=\delta$. So by Lemma \ref{countrymanlemma2},  $\beta<_\mathcal{F}\gamma.$
\end{proof}
\end{theorem} Remember that a partial order $(T,<)$ is a tree if $t\downarrow_T=\{x\in T\,|\,x<t\}$ is well-ordered for any $t\in T$. In particular any tree is well founded.  
    \begin{definition}[Aronszajn tree] Let $(T,<)$ be a tree. We say that $T$ is an $\omega_1$-tree if $Ht(T)=\omega_1$,  and $T_\alpha$ is countable for each $\alpha<\omega_1.$ Furthermore, if $T$ is an $\omega_1$-tree without uncountable chains, we call it an Aronszajn tree.
 
    \end{definition}
    There is a natural way to define an Aronszajn tree from a Countryman line (see \cite{todorvcevic1984trees}). In \cite{lopezschemethesis}, the reader can find a construction (using construction schemes) of a family $\{f_\alpha\}_{\alpha\in \omega_1}$ of functions satisfying that:\begin{itemize}
        \item $f_\alpha:\alpha\longrightarrow \omega$ is one to one for each $\alpha\in \omega_1,$
        \item $f_\beta|_\alpha=^*f_\alpha$ for each $\alpha<\beta\in \omega_1.$
    \end{itemize} 
Such a family naturally generates an Aronszajn tree whose elements are all the functions $g:\alpha\longrightarrow \omega$ for which there is $\beta\in \omega_1$ greater or equal to $\alpha$  such that $g=^*f_\beta|_\alpha$. These kind of trees are called coherent Aronszajn trees.\\\\
Another important class of trees are the so called $\textit{special}$. We say that an $\omega_1$-tree, say $(T,<)$, is called special if it can be written as a countable union of antichains. Since any chain intersect each antichain in at most one point, it follows that each special tree is Aronszajn. Coherent Aronszajn trees are often special, but they might not always be (see \cite{Walks} and \cite{someremarksoncoherent}). In the following theorem, we give a simple construction of a special tree from a construction scheme.
\begin{theorem}There is a special Aronszajn tree.

\begin{proof}
    Let $\mathcal{F}$ be a construction scheme of type $\{(m_k,2,r_{k+1})\}_{k\in\omega}$. Given $\beta\in\omega_1$, consider the function $\rho_\beta:\beta+1\longrightarrow \omega$ defined as $\rho_\beta(\alpha)=\rho(\alpha,\beta).$ Let $$T=\{f\in \omega^{<\omega_1}\,|\,\exists \beta\in\omega_1\,(dom(f)=\beta+1\textit{ and }f=^*\rho_\beta)\}.$$
Given $\alpha<\beta\in\omega_1$, it is straightforward that $\rho_\alpha(\xi)=\rho_\beta(\xi)$ for each $\xi\in (\alpha+1)\backslash(\alpha)_{\rho(\alpha,\beta)}$. From this, we conclude that $f|_{\alpha+1}\in T$ for each $f\in T$ and $\alpha\in dom(f).$ Hence, $(T,\subseteq)$ is an $\omega_1$-tree. Now, we will show that it is in fact a special tree. For this, take an arbitrary $f\in \mathcal{F}$ and let $\beta_f$ be such that $dom(f)=\beta_f+1$. Now, let $k_f\in\omega$ be the minimal natural number with the following properties:\begin{enumerate}
    \item $\rho_{\beta_f}(\xi)=f(\xi)$ for each $\xi\notin (\beta_f)_{k_f}$,
    \item $f(\xi)\leq k_f$ for each $\xi\in (\beta_f)_{k_f}$.
\end{enumerate}
For each $k,s\in\omega$, let $T(k,s)$ be the set of all $f\in T$ for which $k=k_f$ and $|(\beta_f)_{k}|=s.$ We claim that this set is an antichain in $T(k,s)$. Indeed, take distinct $f,g\in T(k,s)$. If $\beta_f=\beta_g$ there is nothing to do, so let us assume that $\beta_f<\beta_g$. Since    $|(\beta_f)_{k}|=s=|(\beta_g)_{k}|$, then $\rho(\beta_f,\beta_g)>k$. By definition of $k_f$ and $k_g$, we get that $f(\beta_f)\leq k<\rho(\beta_f,\beta_g)=g(\beta_f)$. This concludes the claim. To finish, just note that $T=\bigcup\limits_{k,s\in\omega}T(k,s)$.
\end{proof}

\end{theorem}

\subsection{Gaps and coherent families}Let $X$ be a countably infinite set. We say that a family $\mathcal{T}\subseteq [X]^{\omega}$ is a (increasing) $\kappa$-pretower (or simply, pretower) if it is well-ordered by $\subsetneq^*$ and its order type is equal to the ordinal $\kappa$. Whenever we 
index a pretower with an ordered set, it is understood that this is done in an order preserving way. Note that each pretower has an asociated rank function defined in terms of $\subsetneq^*$. 
\begin{definition}[Pregap and gap]Let $\mathcal{A}$ be a $\kappa$-pretower and $\mathcal{B}$ be $\lambda$-pretower, both on a countable set $X$. We say that $(\mathcal{A},\mathcal{B})$ is $(\kappa,\lambda)$-pregap (or simply, pregap) if $A\cap B=\emptyset$ for each $A\in \mathcal{A}$ and $B\in \mathcal{B}$ with $rank(A)=rank(B)$. Moreover, we say that $(\mathcal{A},\mathcal{B})$ is a $(\kappa,\lambda)$-gap (or simply, gap) if it is a $(\kappa,\lambda)$-pregap and there is no $C\subseteq X$ such that $A\subseteq^* C$ and  $B\subseteq^* X\backslash C$ for all $A\in \mathcal{A}$ and $B\in \mathcal{B}$. \end{definition}
The existence of gaps can be deduced by a simple application of Zorn's Lemma. However, by doing this, we do not know the order type of its associated pretowers. In 1909, F. Hausdorff gave a clever construction, in ZFC, of an $(\omega_1,\omega_1)$-gap. Such a gap satisfied the following property.
\begin{definition}[Hausdorff condition]Let $(\mathcal{A},\mathcal{B})$ be an $(\omega_1,\omega_1)$-pregap on $\omega$. We say that $(\mathcal{A},\mathcal{B})$ is  Hausdorff if $\{B\in \mathcal{B}\,|\,rank(B)<rank(A)\textit{ and }B\cap A\subseteq k\}$ is finite for each $A\in\mathcal{A}$ and $k\in\omega.$
\end{definition}
A slight variant is the following.
\begin{definition}[Luzin condition]Let $(\mathcal{A},\mathcal{B})$ be an $(\omega_1,\omega_1)$-pregap on an infinite set $X$. We say that $(\mathcal{A},\mathcal{B})$ is Luzin if $\{B\in \mathcal{B}\,|\,rank(B)<rank(A)\textit{ and }|B\cap A|\leq k\}$ is finite for each $A\in\mathcal{A}$ and $k\in\omega.$
\end{definition}
It is not hard to see that each Hausdorff pregap is in fact a gap, and that each Luzin pregap is Hausdorff. In \cite{lopezschemethesis}, there is a construction of a Hausdorff gap using construction schemes. Here, we give an alternative construction. For the reader interested in knowing more about towers and gaps, we recommend him/her to look at \cite{GapsandTowers}, \cite{rothbergergapsinfragmentedideals}, \cite{frankiewicz1994hausdorff}, \cite{Agapcohomologygroup} ,\cite{ScheepersGaps}, \cite{Walks}, \cite{StevoIlias}, \cite{PartitionProblems}, \cite{AnalyticGaps} and \cite{combinatorialprinciplesonomega1}.

\begin{theorem}\label{hausdorffgapconstruction}Let $\mathcal{F}$ be a construction scheme of type $\{(m_k,2,r_{k+1})\}_{k\in \omega}.$ For each $\alpha\in \omega_1$, define $$A_\alpha=\{2k+\Xi_\alpha(k)\,|\,k\in \omega\backslash 1,\textit{ and }\Xi_\alpha(k)\geq 0\},$$
$$B_\alpha=\{2k+(1-\Xi_\alpha(k))\,|\,k\in\omega\backslash 1,\textit{ and }\Xi_\alpha(k)\geq 0\}.$$
Then, $(\{A_\alpha\}_{\alpha\in \omega_1} ,\,\{B_\alpha\}_{\alpha\in \omega_1})$ is a Hausdorff gap.
\begin{proof} Since $r_k=0$ for infinitely many $k's$, it should be clear that each $A_\alpha$ and $B_\alpha$ is infinite. By definition,  $A_\alpha\cap B_\alpha=\emptyset$ for each $\alpha\in \omega.$ In virtue of part (c) in Lemma \ref{lemmaxi}, we have that $A_\alpha\backslash A_\beta\subseteq \{2k+\Xi_\alpha(k)\,|\,k\leq \rho(\alpha,\beta)\textit{ and }\Xi_\alpha(k)\geq 0\,\}$ and $B_\alpha\backslash B_\beta\subseteq \{2k+(1-\Xi_\alpha(k))\,|\,k\leq \rho(\alpha,\beta)\textit{ and }\Xi_\alpha(k)\geq 0\,\}$ whenever $\alpha<\beta.$ This means  $\{A_\alpha\}_{\alpha\in \omega_1}$ and $\{B_\alpha\}_{\alpha\in \omega_1}$ are both $\omega_1$-pretowers. To prove that the Hausdorff condition is satisfied, take $\beta\in \omega_1$ and $k\in\omega$. We claim  $\{\alpha<\beta\,|\,A_\beta\cap B_\alpha\subseteq k\}\subseteq (\beta)_k$. For this, take an arbitrary $\alpha<\beta$ satisfying $\rho(\alpha,\beta)\geq k.$ By part $(b)$ of Lemma \ref{lemmaxi},  $\Xi_\alpha(\rho(\alpha,\beta))=0$ and $\Xi_\beta(\rho(\alpha,\beta))=1$. This means $2\rho(\alpha,\beta)+1\in A_\beta\cap B_\alpha$, so we are done.

\end{proof}
\end{theorem}
Each $(\omega_1,\omega_1)$-pregap over a set $X$, say $(\{A_\alpha\}_{\alpha\in \omega_1},\{B_\alpha\}_{\alpha\in \omega_1})$ , has an associated $\omega_1$-pretower $\{T_\alpha\}_{\alpha\in \omega_1}$ defined as $T_\alpha=A_\alpha\cup B_\alpha$ for each $\alpha\in \omega_1.$ If we let $f_\alpha$ be the characteristic function of $A_\alpha$ with domain $T_\alpha$, then $f_\alpha=^*f_\beta|_{T_\alpha}$ for each $\alpha<\beta\in \omega_1.$ In this way, we can see that the original pregap is in fact a gap if and only if there is no $f:X\longrightarrow 2$ such that $f_\alpha=^*f|_{T_\alpha}$ for all $\alpha\in \omega_1$ (see \cite{TopicsinSetTheory}, pp. 96-98). Based on these observations, various generalizations of pregaps can be formulated. The following definition was formulated by I. Farah in \cite{farah1996coherent}. Previous work in this subject was done by D. E. Talayco in \cite{talayco1995applications} and C. Morgan in \cite{Agapcohomologygroup}.
\begin{definition}A coherent family of functions supported by an $\omega_1$-pretower $\{ T_\alpha\}_{\alpha\in\omega_1\backslash \omega}$ is a family of functions $\{ f_\alpha \}_{\alpha\in\omega_1\backslash \omega}$ such that:
\begin{enumerate}
    \item $\forall \alpha\in \omega_1\backslash \omega\,(\,f_\alpha:T_\alpha\longrightarrow \alpha\,),$
\item $\forall \alpha<\beta\in \omega_1\backslash \omega\,(\,f_\alpha|_{T_\alpha\cap T_\beta}=^*f_\beta|_{T_\alpha\cap T_\beta}\,).$
\end{enumerate}
Additionally, if for all $\xi<\mu<\beta \in\omega_1\backslash \omega$ and for all $k\in\omega$, the set $$\{\alpha\in \beta\backslash \mu\,|\, |f^{-1}_\alpha[\{\xi\}]\cap f^{-1}_\beta[\{\mu\}]|\leq k\}$$
is finite, we say that $\{ f_\alpha\}_{\alpha\in\omega_1\backslash \omega}$ is Hausdorff.
\end{definition}
Suppose that $\{f_\alpha\}_{\alpha\in \omega_1\backslash \omega}$ is a coherent family of functions supported by some pretower. Given $\xi,\alpha \in \omega_1\backslash \omega$ with $\xi<\alpha$, we let $A^\xi_\alpha=f^{-1}_\alpha[\{\xi\}]$. Observe that in this case $\mathcal{A}^\xi=\{A^\xi_\alpha\}_{\alpha\in \omega_1\backslash \xi}$ is a pretower. Furthermore, given $\xi<\mu$, $(\mathcal{A}^\xi,\mathcal{A}^\mu)$ is a pregap\footnote{Here, we make an slight abuse of notation by thinking of $\mathcal{A}^\xi$ as $\{A^\xi_\alpha\}_{\alpha\in \omega_1\backslash \xi.}$ }.
\begin{definition}Let $\mathfrak{F}=\{f_\alpha\}_{\alpha\in \omega_1\backslash \omega}$ is a coherent family of functions supported by an $\omega_1$-pretower. We say that $\mathfrak{F}$ is Luzin (respectively Hausdorff) if $(\mathcal{A}^\xi,\mathcal{A}^\mu)$ is Luzin (respectively Hausdorff) for each $\xi<\mu\in \omega_1\backslash \omega$.
\end{definition}
In \cite{farah1996coherent}, a Luzin coherent family of functions supported by an $\omega_1$-pretower was constructed in $ZFC$. This was achieved using forcing and appealing to Keisler's Completeness Theorem for $L^\omega(Q)$ (see \cite{logicwiththequantifier}). Now, we present a direct construction of this object. No previous direct construction was known.

\begin{theorem}\label{hausdorffcoherenttheorem}There is Luzin coherent family of functions supported by an $\omega_1$-tower.
\begin{proof}Let $\mathcal{F}$ be a construction scheme of type $\{(m_k, 2,r_{k+1})\}_{k\in\omega}$. Let $N=\bigcup\limits_{k\in\omega} N_k$ where each $N_k$ is equal to $\{k\}\times k\times r_k\times r_k\times k.$ Given $\alpha\in \omega_1\backslash \omega$ let $$T_\alpha=\bigcup\{N_k\,|\,\Xi_\alpha(k)\geq 0\}=\{(k,i,j,s)\in \omega^3\,|\,\Xi_\alpha(k)\geq 0\textit{ and }(i,j,s)\in r_k\times r_k\times k\} .$$ Given $\alpha<\beta\in\omega_1$, we have that $T_\alpha\cap N_k\subseteq T_\beta \cap N_k$ for each $k>\rho(\alpha,\beta)$. This holds by the part $(c)$ of Lemma \ref{lemmaxi}. So in particular, $T_\alpha\subseteq^*T_\beta$. This means that $\mathcal{T}=\{T_\alpha\}_{\alpha\in\omega_1}$ is a pretower.\\\\
For each $\alpha\in\omega_1\backslash \omega$ define $f_\alpha:T_\alpha\longrightarrow \alpha$ as follows:
$$f_\alpha(k,i,j,s)=\begin{cases}(\alpha)_{k}(i)&\textit{ if }\Xi_\alpha(k)=0\\
(\alpha)_{k}(j)&\textit{ if }\Xi_\alpha(k)=1
\end{cases}$$
We claim that $\{f_\alpha\}_{\alpha\in\omega_1\backslash \omega}$ is a coherent family of functions supported by $\mathcal{T}$. For this let $\alpha<\beta\in\omega_1$ and take an arbitrary $(k,i,j,s)\in T_\alpha\cap T_\beta$ for which $k>\rho(\alpha,\beta)$ (observe that all but finitely many $k$'s satisfy this property). By definition of $T_\alpha$, it happens that $\Xi_\alpha(k)\geq 0$. Since $k>\rho(\alpha,\beta)$, we can summon Lemma \ref{lemmaxi} to conclude $\Xi_\alpha(k)=\Xi_\beta(k)$.  From this, it is straightforward that $f_\alpha(k,i,j,s)=f_\beta(k,i,j,s)$. \\\\
Finally we prove that the family is Luzin. For this, let $\xi<\mu<\beta\in\omega_1$ and $k\in\omega$. We claim that if $l=\max(k,\rho(\{\xi,\mu,\beta\}))$, then  $\{\alpha\in \beta\backslash \mu\,|\, |f^{-1}_\alpha[\{\xi\}]\cap f^{-1}_\beta[\{\mu\}]|\leq k\}\subseteq (\beta)_{l}$. Indeed, take an arbitrary  $\alpha\in\beta\backslash (\mu\cup (\beta)_l)$. By simplicity, name $\rho(\alpha,\beta)$ as $t$. Due to the part $(b)$ of Lemma \ref{lemmaxi}, $\Xi_\alpha(t)=0 $ and $\Xi_\beta(t)=1$. Furthermore, $\xi,\mu\in (\beta)_t\cap \alpha=(\alpha)^-_t$. In this way, $\Xi_\xi(t)\leq\Xi_\mu(t)\leq \Xi_\alpha(t)<\Xi_\beta(t)$. As $t>\max(\rho(\xi,\beta),\rho(\mu,\beta))$, we conclude that $\Xi_\xi(t)=\Xi_\mu(t)=-1$ using the part $(b)$ of Lemma \ref{lemmaxi}. Thus, both $i=|(\xi)^-_t|$ and $j=|(\xi)^-_t|$ are both strictly smaller than $r_t$. This means $(t,i,j,s)\in T_\alpha\cap T_\beta$ for each $s<t$.  Furthermore, for all such $s$ we have that $f_\alpha(t,i,j,s)=(\alpha)_t(i)=\xi$ and $f_\beta(t,i,j,s)=\mu$. Consequently, $\{t\}\times\{i\}\times\{j\}\times t\subseteq f^{-1}_\xi[\{\alpha\}]\cap f^{-1}_{\mu}[\{\beta\}]$. As this set is of cardinality $t$ which is bigger than $k$, we are done. 
\end{proof}    
\end{theorem}

\section{Beyond ZFC}\label{beyondzfcsection} Theorem \ref{luzinjonestheorem} is a prototypical example of how construction schemes can be applied. We want to build an uncountable object $\{A_\alpha\}_{\alpha\in \omega_1}$; in order to do so, for each $\alpha\in \omega_1$ we define a sequence of approximations of $A_\alpha$, say $\{A^k_\alpha\}_{k\in \omega}$. Typically, the definition of $A^k_\alpha$ solely depends on $|(\alpha)_k|$. This means that whenever $\alpha,\beta\in \omega_1$ and $k<\Delta(\alpha,\beta)$, the respective approximations $A^k_\alpha$ and $A^k_\beta$ are the same. Whenever $k=\rho(\alpha,\beta)$, we define the approximations $A^k_\alpha$ and $A^k_\beta$ in such way that the part (b) of Lemma \ref{lemmaxi} can be used to force a certain interaction between $A_\alpha$ and $A_\beta$. In the case of the Luzin-Jones family, at this step we forced $A_\alpha$ and $A_\beta$ to have at least $k$ elements in their intersection. Finally, for each $k>\rho(\alpha,\beta)$, we make sure to define the approximations in such way that the part (c) of Lemma \ref{lemmaxi} can be used in our advantage. In the case of the Luzin-Jones family, we forced $A^k_\alpha$ and $A^k_\beta$ to be disjoint from the moment $k$ was greater than $\rho(\alpha,\beta).$ All of this means that the only part of the construction where we do not have any control about the interaction between $A^k_\alpha$ and $A^k_\beta$ is when $\Delta(\alpha,\beta)\leq k<\rho(\alpha,\beta).$ In this way, the best scenario happens when $\Delta(\alpha,\beta)=\rho(\alpha,\beta)$. In this section we will be interested in construction schemes for which many pairs of ordinals satisfy this equality.\\\\
Recall that $\rho^C=\max\{\rho(\alpha,\beta)\,|\,\alpha,\beta\in F\}$ for each $C\in[\omega_1]^{<\omega}$. Furthermore, if $F\in \mathcal{F}$ then  $\rho^F=l$ if and only if $F\in \mathcal{F}_l.$
\begin{definition}Let $\mathcal{F}$ be a construction scheme and $\mathcal{C}=\{c_i\}_{i<k}\subseteq [\omega_1]^{<\omega}$. For $F\in \mathcal{F}$, we say that $F$ captures $\mathcal{C}$ (or $\mathcal{C}$  is captured by $F$) if:\begin{enumerate}
    \item $n_{\rho^F}\geq |\mathcal{C}|$,
    \item $\forall i<k\,(\,c_i\subseteq F_i\textit{ and }c_i\backslash R(F)\not=\emptyset)$
    \item $\forall i<k\,(\, \phi_i[c_0]=c_i\textit{ where }\phi_i \textit{ is the increasing bijection from }F_0\textit{ to }F_i\,).$
\end{enumerate} 
Moreover,  if $n_{\rho^F}=k$, we say that $F$ fully captures $\mathcal{C}$. 
\end{definition}
\begin{remark}Given $\mathcal{C}\subseteq [\omega_1]^{<\omega}$ and $F\in \mathcal{F}$, we have that $F$ captures $\mathcal{C}$ if and only if $n_{\rho^F}\geq |\mathcal{C}|$ and there is $c\subseteq m_{\rho^F-1}$ for which $\max(c)\geq r_{\rho^F}$ and $\mathcal{C}=\{F_i[c]\,|\,i<|\mathcal{C}|\}$.
\end{remark}
It is not hard to see that if $\mathcal{C}$ is captured by some element of $\mathcal{F}$, then that element must belong to $\mathcal{F}_{\rho^\mathcal{C}}.$ In this case, every element of $\mathcal{F}_{\rho^\mathcal{C}}$ which contains $\bigcup\mathcal{C}$ will also capture $\mathcal{C}$.\\
Whenever $C\in [\omega_1]^{<\omega}$, we make a small abuse of notation and say that $F$ captures (respectively, fully captures) $C$ instead of saying that the same is true for $\{\{\alpha\}\,|\,\alpha\in C\}$. Regarding this case, we have the following equivalence. Its proof is left to the reader.

\begin{proposition}Let $\mathcal{F}$ be a construction scheme and $C\in [\omega_1]^{<\omega}$ be nonempty.  $C$ is captured by some element of $\mathcal{F}_l$ if and only if $\Xi_{C(i)}(\rho^C)=i$ for each $i<|C|$ and $\Delta(\alpha,\beta)=\rho(\alpha,\beta)=l$ for all $\alpha,\beta\in C$.
\end{proposition}

\begin{definition}[\cite{treesandgapsschemes}]Let $\mathcal{F}$ be a construction scheme, $\mathcal{P}$ a partition of $\omega$ and $n\in\omega\backslash 1$. We say that $\mathcal{F}$ is $n$-$\mathcal{P}$-capturing if for each uncountable $S\subseteq[\omega_1]^{<\omega}$, $P\in\mathcal{P}$ and $k\in\omega$ there are $\mathcal{C}\in[S]^{n}$, $k<l\in P$ and $F\in \mathcal{F}_l$ which captures $\mathcal{C}$. Whenever $\mathcal{P}=\{\omega\}$, we simply say that $\mathcal{F}$ is $n$-capturing.
\end{definition}
\begin{definition}[\cite{schemeseparablestructures}]Let $\mathcal{F}$ be a construction scheme and $\mathcal{P}$ a partition of $\omega$. We say that $\mathcal{F}$ is $\mathcal{P}$-fully capturing if for each uncountable $S\subseteq[\omega_1]^{<\omega}$, $P\in\mathcal{P}$ and $k\in\omega$ there are $\mathcal{C}\in[S]^{<\omega}$, $k<l\in P$ and $F\in \mathcal{F}_l$ which fully captures $\mathcal{C}$. Whenever $\mathcal{P}=\{\omega\}$, we  simply say that $\mathcal{F}$ is fully capturing. 
\end{definition}
Now, we introduce the capturing axioms for construction schemes. This axioms have already been studied in \cite{forcingandconstructionschemes}, \cite{lopezschemethesis}, \cite{treesandgapsschemes} and \cite{schemeseparablestructures}. 
\begin{itemize}
    \item {\bf Fully Capturing Axiom} {[\bf FCA]}: There is a fully capturing construction scheme of every possible type.\\
\item { \bf Fully Capturing Axiom with Partitions [FCA(part)]}: There is a $\mathcal{P}$-fully capturing construction scheme for every type $\tau$  and each partition $\mathcal{P}$ compatible with $\tau$.\\
\item {\bf $n$-Capturing Axiom [CA$_n$]}: There is an $n$-capturing construction scheme of every possible type safisfying that $n\leq n_k$ for each $k\in \omega\backslash1$.\\
\item {\bf $n$-Capturing Axiom with Partitions [CA$_n$(part)]}: There is a $\mathcal{P}$-$n$-capturing construction scheme for every type $\tau$ satisfying that $n\leq n_k$ for each $k\in \omega\backslash 1$  and each partition $\mathcal{P}$ compatible with $\tau$. 
\end{itemize}

 \subsection{Types of Suslin trees}In this subsection we study a particular class of Aronszajn trees called Suslin. T. Jech and S. Tennenbaum independently proved the consistency of the existence of Suslin trees in \cite{Jechsuslin} and \cite{Tennenbaumsuslin} respectively. Later, R. Jensen proved  that the existence of Suslin trees follows from the $\Diamond$-principle. Finally, R. M. Solovay and S. Tennenbaum proved that $ZFC$ is consistent with the non existence of Suslin trees (see \cite{Iteratedcohensuslin}).
\begin{definition}Let $(T,<)$ be an $\omega_1$-tree. we say that $T$ is Suslin if $(T,>)$ does not contain uncountable chains or antichains.
\end{definition}
\begin{definition}Let $(T,<)$ be a Suslin tree. We say that $T$ is coherent if there is a family of functions $\{f_\alpha\}_{\alpha\in\omega_1}$ such that:\begin{enumerate}
    \item $\forall \alpha\in\omega_1\,(\,f_\alpha:\alpha\longrightarrow \omega\,),$
    \item $\forall \alpha<\beta\in\omega_1\,(\,f_\alpha=^*f_\beta|_\alpha\,),$
    \item $T=\{f_\alpha|_\xi\,|\,\xi<\alpha<\omega_1\}$ and $<=\subseteq.$
\end{enumerate}
\end{definition}
Coherent properties have been extensively studied in the past, since they have many interest forcing properties (see \cite{microscopic1}, \cite{Microscopic2},  \cite{variationsouslin}, \cite{katetovproblem}, \cite{chainconditionsmaximalmodels} and \cite{gapstructure}). 
\begin{theorem}[FCA(part)]There is a coherent Suslin tree.
\begin{proof}Let $\mathcal{F}$ be a $\mathcal{P}$-fully capturing construction scheme of type $\{(m_k,n_{k+1},r_{k+1})\}_{k\in\omega}$ with $n_{k+1}\geq 2^{m_k-r_{k+1}}+1$ for each $k\in\omega$ and such that $\mathcal{P}=\{P_c,P_a\}$ for some infinite sets $P_c$ and $P_a$. Given $k\in\omega$, enumerate (possibly with repetitions) the set of all functions from $m_k\backslash r_{k+1}$ into $2$ as $\{g^k_i\}_{0<i<n_{k+1}}$. Now, let $\beta\in\omega_1$ and define $f_\beta:\beta\longrightarrow 2$ as follows:
$$f_\beta(\xi)=\begin{cases}1 &\textit{if } \Xi_\xi(\rho)= 0,\,\Xi_\beta(\rho)=1\textit{ and }\rho\in P_c\\
g^\rho_i(|(\xi)^-_\rho|)&\textit{if }\Xi_\xi(\rho)=0,\,\Xi_\beta(\rho)=i>1\textit{ and }\rho\in P_a\\
0&\textit{ otherwise }
\end{cases}$$
Here, $\rho=\rho(\xi,\beta)$.\\
Take $\alpha<\beta\in \omega_1$. We will prove that $f_\alpha=^*f_\beta|\alpha$. For this, consider an arbitrary $\xi<\alpha$ satisfying $\xi\notin (\alpha)_{\rho(\alpha,\beta)}$. Since $\rho$ is an ordinal metric, we can use the part $(c)$ of Definition \ref{ordinalmetricdefinition} to conclude $\rho(\xi,\alpha)=\rho(\xi,\beta)$. We will make an abuse of notation by calling this number just $\rho$. By the part (c) of Lemma \ref{lemmaxi}, either $\Xi_\alpha(\rho)=-1$ or $\Xi_\alpha(\rho)=\Xi_\beta(\rho)$. But  $\Xi_\xi(\rho)<\Xi_\alpha(\rho)$ due to part $(b)$ of Lemma \ref{lemmaxi}. In this way, $\Xi_\alpha(\rho)$ is equal to $\Xi_\beta(\rho)$. By definition of $f_\alpha$ and $f_\beta$,  $f_\beta(\xi)=f_\alpha(\xi)$.\\
Now, we will prove  $T=\{f_\alpha|_\xi\,|\,\xi<\alpha<\omega_1\}$ is a Suslin tree. First we start by proving  it does not have uncountable chains. For this, take $S\in[T]^{\omega_1}$ and let $$\mathcal{C}=\{\,\{\xi,\alpha\}\in [\omega_1]^2\,|\,f_\alpha|_\xi\in S\}.$$
Without loss of generality we can suppose there is $k\in\omega$ such that $\rho(\xi,\alpha)=k$ for each $\{\xi,\alpha\}\in \mathcal{C}$ and that $\{\xi,\alpha\}\cap \{\delta,\beta\}=\emptyset$ whenever $\{\xi,\alpha\},\{\delta,\beta\}\in \mathcal{C}$. Since $\mathcal{F}$ is, in 
particular, $\mathcal{P}$-$3$-capturing, there are $l\in P_c$ greater than $k$, $F\in \mathcal{F}_l$ and $\{\xi_0,\alpha_0\},\{\xi_1,\alpha_1\}$, $\{\xi_2,\alpha_2\}\in \mathcal{C}$ which are captured by $F$. Observe  $\xi_0<\xi_1<\xi_2$,  $\rho(\xi_0,\alpha_1)=\rho(\xi_0,\alpha_2)=l$,  $\Xi_{\xi_0}(l)=0$, $\Xi_{\alpha_1}(l)=1$ and $\Xi_{\alpha_2}(l)=2$. Hence, $f_{\alpha_1}(\xi_0)=f_{\alpha_1}|_{\xi_1}(\xi_0)=1$ and  $f_{\alpha_2}(\xi_0)=f_{\alpha_2}|_{\xi_2}(\xi_0)=0$. This implies $S$ is not a chain.\\
Finally, we prove  $T$ has no  uncountable antichains. For this, let $A\in [T]^{\omega_1}$. Without loss of generality we can suppose that each element of $A$ is of the form $f_\alpha$ for some $\alpha\in\omega_1$. Now, let $\mathcal{D}=\{\alpha\in\omega_1\,|\,f_\alpha\in A\}$. Since $\mathcal{F}$ is $\mathcal{P}$-fully capturing, there are $l\in P_a$, $F\in \mathcal{F}_l$ and $\{\alpha_0,\dots,\alpha_{n_l-1}\}\subseteq \mathcal{D}$ captured by $F$. Let $g:m_{l-1}\backslash r_l\longrightarrow 2$ be given by:$$g(j)=\begin{cases}f_{\alpha_0}(F(j))&\textit{if } F(j)<\alpha_0\\
0&\textit{otherwise}
\end{cases}$$
 Take $0<i<l$ satisfying $g=g^l_i$. We claim $f|_{\alpha_0}=f_{\alpha_i}|_{\alpha_0}$. For this, take an arbitrary $\xi<\alpha_0$. If $\xi\in F_0\backslash R(F)$, then   $f_{\alpha_i}(\xi)=f_{\alpha_0}(\xi)$ by definition of $g^l_i$. If $\xi \in R(F)$, then $f_{\alpha_i}(\xi)=f_{\alpha_0}(\xi)$ by the part $(a)$ of Lemma \ref{lemmaxi} (since $\Delta(\alpha_0,\alpha_i)=l$). The only case left is when $\xi\notin F\cap \alpha_0=(\alpha)^-_l$. But we proved in a previous paragraph that whenever this situation happens, then $f_{\alpha_0}(\xi)=f_{\alpha_i}(\xi)$. Having proved the claim, we conclude $A$ is not an antichain.
\end{proof}
\end{theorem}
Suppose $T$ is a coherent Suslin tree and $G$ is a generic filter over $T$ (seen as a forcing notion). In the generic extension $T$ has at least $\omega_1$ uncountable branches. Indeed, the function  $\bigcup G:\omega_1\longrightarrow \omega$ determines an uncountable branch through $T.$ Using that $T$ is coherent, it is not hard to see that every finite modification of $\bigcup G$ has the same property. The underlying reason for this is that, in the ground model, $T$ has many nontrivial automorphisms.\\
Our next goal is to construct  Suslin trees which are substantially different from the coherent ones. These trees not only have a trivial automorphism group, but in fact, only one uncountable branch is added when we force with them. The reader interested in knowing more about the automorphisms of Suslin trees is invited to look for \cite{rigiditysouslin}.
\begin{definition}Given $k\in\omega$ and $\{(T_i,<_i)\}_{i< k}$ a family of trees we define their tree product $$\bigotimes\limits_{i< k} T_i=\{t\in \prod\limits_{i<k}T_i\,|\,\forall i<k\,\big( rank(t(0))=rank(t(i))\big)\}.$$
To this set, we associate a canonical order given by: $$s<t\textit{ if and only if }s(i)<_i t(i)\textit{ for each }i<k$$
\end{definition}
It is not hard to see that $\bigotimes\limits_{i<k}T_i$ is always a tree. Furthermore, the tree product of Aronszajn trees is always Aronszajn. Unfortunately, we can not say the same about the tree product of Suslin trees. In fact, the tree product of a Suslin tree with itself is never Suslin.\\
Given a tree $(T,<)$ and $t\in T$, we let $T|
_t$ be the set of all elements of $T$ which are comparable with $t.$ Whenever $\beta\leq Ht(T)$, $T|_\beta$ denotes the set $$\bigcup\limits_{\alpha<\beta}T_\alpha.$$ 
\begin{definition}Let $(T,<)$ be a tree. We say that $T$ is full Suslin if $\bigotimes\limits_{i<k} T|_{t_i}$
is Suslin for every distinct $t_0,\dots,t_{k-1}\in T$, all of the same rank.
\end{definition}
\begin{theorem}[FCA]\label{suslintreescheme}There is a full Suslin tree.
\end{theorem}
Before advancing, we will describe the amalgamation that will be used throughout the proof. For this, let $(T,<_T)$ and $(L,<_L)$ be two finite trees and let $R$ be the intersection of $T$ and $L$. Suppose that both orders agree when restricted to $R$, and that there is  $l<\min(Ht(T),Ht(L))$ such that $R= T|_l=L|_l.$
Let $\mathcal{B}=\{B_t\}_{t\in T_l}$ be a set of branches in $L$ such that $B_t\cap R= t\downarrow_T$ for each $t\in T_l$. For $x,y\in T\cup L$ we let $x<_\mathcal{B} y$ if one of the following happens:
\begin{itemize}
    \item $x,y\in T$ and $x<_T y,$
    \item $x,y \in L$ and $x<_L y,$
    \item there is $t\in T_l$ such that $x\in B_t$ and $t\leq_T y.$
\end{itemize}
We call $(T\cup L, <_\mathcal{B})$ the amalgamation of $T$ and $L$ through $\mathcal{B}$.
\begin{proposition}$(T\cup L,<_\mathcal{B})$ is a tree.
\begin{proof}First we check that $<_\mathcal{B}$ is a partial order. For this, it suffices to prove transitivity. Let $x,y,z\in T\cup L$ be such that $x<_\mathcal{B} y$ and $y<_\mathcal{B} z$. We claim it can not happen that $x,z\in L$ and $y\in T$. Suppose towards a contradiction that it does. Observe  $y<_\mathcal{B} z$ implies  $y\in L$, which means $y\in R$.  But this is impossible because $x<_\mathcal{B} y$ implies  $rank(y)\geq l.$ Thus, the only nontrivial cases occur when $x,y\in L$ and $z\in T$ or when $x\in L$ and $y,z\in T$. In the first case, let $t\in T_l$ for which $y\in B_t$ and $t\leq_T z$. Since $B_t$ is a branch,  $x\in B_t$ and consequently $x<_\mathcal{B} z$. In the second case, take $t\in T_l$ such that $t\leq_T y$ and $x\in B_t$. Just observe $t<_T z$. Thus, $x<_\mathcal{B} z.$\\
To finish the proof, take $y\in T\cup L$. Notice that if $y\in L$, then $y\downarrow_{T\cup L}=y\downarrow_{L}$. On the other hand, if $y\in T\backslash L$ and $t\in T_l$ is such that $t\leq_T y$, then $y\downarrow_{T\cup L}=B_t\cup y\downarrow_{T}$. In both cases, the respective set is a totally ordered. 
\end{proof}
\end{proposition}
Let $(T,<)$ be a finite tree and $l<Ht(T)$. We say that $C\subseteq T$ is $l$-good if $rank(y)\geq l$ for each $y\in C$ and  $(y\downarrow_T)\cap (z\downarrow_T)\subseteq T|_l$ for any other $z\in C$. This last property is equivalent to saying that the only element of $y\downarrow_T\cap T_l$ does not belong to $z\downarrow_T$.
\begin{proof}[Proof of Theorem \ref{suslintreescheme}] Fix a type $\tau=\{(m_k, n_{k+1},r_{k+1})\}_{k\in\omega}$ such that $n_{k+1}\geq m_k(k+1)2^{m_k}$ for each $k\in\omega.$ Let $\mathcal{F}$ be a fully capturing construction scheme of type $\tau.$
For each $F\in \mathcal{F}$, let $$T^F=\bigcup\limits_{l<m_{\rho^F}}\{F(l)\}\times((\rho^F+1)2^l).$$
Given $F,G\in \mathcal{F}$ such that $\rho^F=\rho^G$, let $\phi_{F,G}:T^F\longrightarrow T^G$ be the natural bijection. That is, for each $j<m_{\rho^F}$ it happens that $\phi_{F,G}(F(j),s)=(G(j),s)$.\\ If $\rho^F=k$ then $|T^F|\leq m_k(k+1)2^{m_k}\leq n_{k+1}$. This means we can enumerate  $\mathscr{P}(T^F)$ (possibly with repetitions) as $\{C^F_i\}_{i<n_{k+1}}$. Moreover, we can do this in such way that $C^G_0=\emptyset$ and $C^G_i=\phi_{F,G}[C^F_i]$ whenever $F,G\in \mathcal{F}_k$ and $i<n_{k+1}$.\\
Our goal is to define for each $k\in\omega$ and $F\in \mathcal{F}_k$, partial orders $<_F$ over $T^F$ satisfying the following requirements:

    \begin{enumerate}
        \item $(T^F,<_F)$ is a tree,
        \item $Ht(T^F)=m_k$,
        \item $\forall l< m_k\,(\,T^F_l= \{F(l)\}\times( (k+1)2^l)\,)$
        \item Every non maximal element of $T^F$ has exactly 2 immediate succesors in $T^F.$
    \end{enumerate}
    Also, we will ask that:\begin{enumerate}[label=(\alph*)]
        \item  $\phi_{F,G}$ is an order isomorphism for each $F,G\in \mathcal{F}$ satisfying $\rho^F=\rho^G.$
        \item For each  $k\in\omega$, if $F\in \mathcal{F}_{k+1}$ and $i<n_{k+1}$, then $<_F\cap (F_i\times F_i)=<_{F_i}.$
        \item For each $k\in \omega\backslash 1$, if $F\in \mathcal{F}_{k}$ and $C\subseteq T^{F_0}$ is $r_k$-good, then there is $i<n_k$ such that $x<_{T^F} \phi_{F_0,F_i}(x)$ for each $x\in C.$
    \end{enumerate}
The construction is done by recursion.\\\\
(Base step) For each $F\in \mathcal{F}_0$, let $<_F=\emptyset.$\\\\
(Recursion step) Let $k\in\omega$ and suppose  we have defined $<_G$ for each $G\in \mathcal{F}_k$. Let $F\in \mathcal{F}_{k+1}$ and for each $i<n_{k+1}$ define $T^i$ as $$ \bigcup\limits_{l<i(m_k-r_{k+1})+m_k}\{F(l)\}\times((k+1)2^l).$$
    First, we are going to recursively define the order $<_F$ restricted to each $T^i$, namely $<_i$. If $i=0$, observe that $T^i=T^{F_0}$. Thus, let  $<_i$ be simply $<_{F_0}.$ Now, suppose we have defined $<_i$ for some $i< n_{k+1}-1$ in such way that conditions (1), (2), (3), (4) and (b) are satisfied when relativized in the obvious way. $T^i\cap T^{F_{i+1}}=T^{0}|_{r_{k+1}}$, so we are in conditions to 
    amalgamate these trees as described in previous 
    paragraphs. The only thing missing, is a suitable set of branches. For this, consider $C^{F_0}_{i+1}$. If this set is $r_{k+1}$-good in $T^i$, let $C$ be a maximal $r_{k+1}$-good set which contains it. In
    other case, let $C$ be just an  arbitrary maximal $r_{k+1}$-good set. Since $C$ is $r_{k+1}$-good, for each $x\in T^{i+1}_{r_{k+1}}$, there is at most one element of
    $C$ which is  greater or equal than $\phi^{-1}_{F_0,F_{i+1}}(x).$ By maximality, such element indeed exists. Let us call it $c_x$. For each $x\in T^{i+1}_{r_{k+1}}$, let $B_x$ be a branch in $T^i$ with $c_x\in B_x$ and consider $\mathcal{B}=\{B_x\}_{x\in T^{i+1}_{r_{k+1}}} $. By definition of $c_x$, 
    $\phi^{-1}_{F_0,F_{i+1}}(x)\in B_x.$ Thus, by part $(a)$ of the  conditions, we conclude $B_x\cap T^{0}=x\downarrow_{T^{F_{i+1}}}$ for every $x$. In this way, we get that the relation $<_\mathcal{B}$ turns $T^i\cup T^{F_{i+1}}$ into a tree.  Finally, just extend the order $<_\mathcal{B}$ to $T^{i+1}$ in such way that all the conditions are satisfied. As a final remark, observe that in the case $C^{F_0}_{i+1}\subseteq C$, for each $y\in C^{F_0}_{i+1}$ there is $x$ for which $y\in B_x$, which means  $y=c_x$. Consequently,  $\phi^{-1}_{F_0,F_1}(x)<_{F_0} y$. By part $(a)$ of the conditions, we know $x<_{F_{i+1}}\phi_{F_0,F_1}(y) $. Thus, by definition of $<_\mathcal{B}$, the inequality  $y<_\mathcal{B} \phi_{F_0,F_1}(y)$ holds.\\
    Once we have defined $<_i$ for each $i<n_{k+1}-1$, let $<_F$ be an order over $T$ which extends each $<_i$ and such that the conditions (1),(2),(3), and (4) are satisfied. By construction, we have that conditions (b) and (c) are automatically satisfied. Also it should be clear we can copy the definition of $<_F$ from some $F\in \mathcal{F}_{k+1}$ to any other $G\in \mathcal{F}_{k+1}$, in order to guarantee that condition (a) is fulfilled. In this way, we finish the recursion.\\
    Let $T=\omega_1\times \omega$ and observe $$T=\bigcup\limits_{F\in \mathcal{F}}T_F.$$
    Moreover, $\{T_F\}_{F\in \mathcal{F}}$ is cofinal in $[T]^{<\omega}$. Thus, it is natural to define $<_T$ as $$\bigcup\limits_{F\in \mathcal{F}}<_F.$$
By conditions (1), (a), (b) and cofinality, we conclude  $(T,<_T)$ is a tree. Condition (3) assures  $H(T)=\omega_1$. Furthermore,   $T_\alpha=\{\alpha\}\times \omega$ for each $\alpha\in \omega_1$. Thus, $(T,<_T)$ is an $\omega_1$-tree. The only thing remaining is to prove that all the suitable tree products are Suslin. For this, let $l\in \omega$ and  $t(0),\dots, t({l-1})\in T$ be distinct elements all of the same rank, say $\alpha.$ Given $i<l$, let $q_i\in\omega$ for which $t(i)=(\alpha,q_i)$. We will prove that  $$L=\bigotimes\limits_{i<l}T|_{t_i}$$
is Suslin. By condition (4), every element of $L$ has at least $2$ succesors. It is well-known that it suffices to prove  $L$ has no uncountable antichains. This is because every element of $L$ has at least two elements above  which are not compatible.\\Let $\mathcal{A}\in [L]^{\omega_1}$. For each $x\in \mathcal{A}$ and $i<l$ there is $s^x_i\in \omega$ for which $x(i)=(\beta_x,s_i^x)$, where $\beta_x=rank_T(x(0))$. By a refining argument, we can suppose that for each $i<l$, there is $s_i\in \omega$ with $s_i^x=s_i$ for every $x\in \mathcal{A}.$ Notice that in this case, we 
will have $\beta_x\not=\beta_y$ whenever $x\not= y$. By refining one last time, we can suppose $\beta_x>\alpha$ whenever $x\in \mathcal{A}.$  Let $S=\{\{\alpha,\beta_x\}\,|\,x\in \mathcal{A}\}$. Since $\mathcal{F}$ 
is fully capturing, there are $k>\max(\{s_j\}_{j<l}\cup\{q_j\}_{j<l})+1$, $F\in \mathcal{F}_k$ and $\{x_i\}_{j<n_{k}}\subseteq [\mathcal{A}]^{n_k}$ such that $F$ captures $\{\{\alpha,\beta_{x_j}\}\}_{j<n_k}$. Since $k$ is big enough, we conclude  $t(i),x_j(i)\in T^{F_j}$ for each $j<n_k$ and every $i<l$. By definition of $L$ and because $\alpha\in R(F)$, we have that for each $i\not= j$ $$x_0(i)\downarrow_{T^{F_0}}\cap x_0(j)\downarrow_{T^{F_0}}= t(i)\downarrow_{T^{F_0}}\cap t(j)\downarrow_{T^{F_0}}\subseteq T^{F_0}|_{r_k}.$$
Thus, $C=\{x_0(i)\}_{i<l}$ is $r_k$-good. This means there is $j<n_k$ for which  $$x_0(i)<_F \phi_{F_0,F_j}(x_0(i))=x_j(i)$$
for each $i<l$. Consequently $x_0<x_j$, so we are done.
\end{proof}

Originally, it was proved in \cite{treesandgapsschemes} that the existence of a $3$-capturing construction scheme implies the existence of a Suslin tree, but unfortunately we believe there is a gap in the proof. This left us with the following problem.
\begin{Problem}Does $CA_n(part)$ imply the existence of Suslin trees for some $n\in\omega\backslash\{0\}$?
\end{Problem}
In general, it would be interesting to know which kind of Suslin trees can be constructed with different forms of capturing. 
\subsection{Suslin lower semi-lattices} We say that $(L,<, \wedge)$ is a lower semi-lattice whenever $(L,<)$ is a partial order such that $\sup\{z\in X\,|\,z\leq x\textit{ and }z\leq y\}$ exists and it is equal to $x\wedge y$ for each $x,y\in L$. 

\begin{definition}[\cite{dilworth2007lattice}]Let $(X,<,\wedge)$ be a lower semi-lattice. We say that $(X,<,\wedge)$ is  Suslin  if \begin{enumerate}
    \item $(X,<)$ is well founded,
    
    \item $X$ is uncountable,
    \item $X$ does not contain any uncountable chain or an uncountable set of pairwise incomparable elements.
\end{enumerate}
\end{definition}

For the rest of this subsection, a set of pairwise incomparable elements will be called  $\textit{pie}$\footnote{Notice that  in the forcing sense, antichains are distinct from pies}. Suslin lower semi-lattices were first studied by S. J. Dilworth, E. Odell and B. Sari in \cite{dilworth2007lattice}, in the context of Banach spaces. In \cite{raghavan2014suslin}, D. Raghavan and T. Yorioka proved that, assuming the $\Diamond$-principle, there is a Suslin lower semi-lattice $\mathbb{S}$ which is a substructure of $(\mathscr{P}(\omega),\subseteq,\cap)$ and such that $\mathbb{S}^n$ does not contain any uncountable pie for each $n\in\omega$. A partial order which satisfies this last property is said to be a \textit{powerful pie}.\\
In the following theorem, we show that $CA_2$ is all that is needed in order to construct a Suslin lower semi-lattice with the above mentioned properties.

\begin{theorem}[$CA_2$]\label{suslinlatticescheme} There is $\mathbb{S}\subseteq \mathscr{P}(\omega)$ such that $(\mathbb{S},\subseteq,\cap)$ is a Suslin lower semi-lattice which is powerful pie.
\begin{proof}Fix a type $\tau=\{(m_k,2,r_{k+1})\}_{k\in\omega}$. For each $k\in\omega$, consider $A_k=m_k\times 2^k$ and $U_k=\{k\}\times m_k\times 2^{k-1}$. Since the last expression has no sense when $k=0$, we let $U_0=\{(0,0,0)\}$. Also, let $\phi_k:A_k\longrightarrow A_{k+1}$ be given as:
$$\phi_k(a,b)=\begin{cases}(a,b)&\textit{if }(a,b)\in r_{k+1}\times 2^k\\
(a+(m_k-r_{k+1}),b)&\textit{in other case}
\end{cases}$$
Notice that $A_k\cup \phi_k[A_k]=m_{k+1}\times 2^k.$ As the final part of the preparation, let $$N_k=\bigcup\limits_{i\leq k}U_k.$$
Our first objective is to construct, for each $k\in\omega$, a family $\{S^k_x\}_{x\in A_k}\subseteq \mathscr{P}(N_k)$ in such way that if $\mathbb{S}^k=\{\emptyset\}\cup \{S^k_x\}_{x\in A_k}$,  the following conditions hold:
\begin{enumerate}[label=(\alph*)]
    \item $(\mathbb{S}^k,\subseteq, \cap)$ is a lower semi-lattice.
    \item $\mathbb{S}^k_0=\{\emptyset\}$ and $\mathbb{S}^k_{i+1}=\{S^k_x\}_{x\in\{i\}\times {2^k}}$ for all $i<m_k$ (where $\mathbb{S}^k_i$ is the set of all element of $\mathbb{S}^k$ of rank $i$).
    \item For all $x\in A_k$, $S^{k+1}_x=S^k_x$ and $S^{k+1}_{\phi_{k}(x)}\cap N_k=S^k_x.$ In particular $S^k_x\subseteq S^{k+1}_{\phi_k(x)}.$
    \item The function $\psi_k:\mathbb{S}^k\longrightarrow\mathbb{S}^{k+1}$ given as:$$\psi_k(x)=\begin{cases}\emptyset &\textit{if } x=\emptyset\\
    S^{k+1}_{\phi_k(y)}&\textit{if }x=S^k_y
    \end{cases}$$
    is an embedding for each $k\in \omega.$
\end{enumerate}
The construction is done by recursion over $k$.\\\\
{(Base step)} If $k=0$, then $A_k=\{(0,0)\}$. In this case, we let $S^0_{(0,0)}=\{(0,0,0)\}$. Trivially, all the conditions are satisfied.\\\\
{(Recursion step)} Suppose we have defined $\mathbb{S}^k$ for some $k\in\omega$ in such way that all the conditions are satisfied. In order to define $\mathbb{S}^{k+1}$, we first divide $A_{k+1}$ into three quadrants:
\begin{itemize}
    \item $C_0=A_k$,
    \item $C_1=[m_k,m_{k+1})\times 2^k$,
    \item $C_2=m_{k+1}\times [2^k,2^{k+1}),$
\end{itemize}   
Now, take an arbitrary $x\in A_{k+1}$ and consider the following cases:\begin{enumerate}[label=$(\roman*)$]
    \item If $x\in C_0$, let $S^{k+1}_x=S^k_x$. 
    \item If $x\in C_2$, then $x=(a,b)$ for some $a<m_{k+1}$ and $2^k\leq b<2^{k+1}$. In this case, let $S^{k+1}_x=\{k+1\}\times(a+1)\times\{b-2^k\}$.
    \item If $x\in C_1$, let $z=\phi^{-1}_k(x)$ and consider $D_x=\{b<2^k\,|\,S^k_{(r_{k+1},b)}\subseteq S^k_z\}.$
    Observe that $D_x$ codes the elements of $\mathbb{S}^k_{r_k+1}$ which are below $S^k_z.$ In this case, let $$S^{k+1}_x=S^k_z\cup\big(\bigcup\limits_{b\in D_x}\{k+1\}\times m_k\times \{b\}\big).$$
\end{enumerate}
It follows directly that condition (c) is satisfied. In particular, $\mathbb{S}^k\subseteq \mathbb{S}^{k+1}$. We will prove  $\mathbb{S}^{k+1}$ is a lower semilattice by showing it is closed under intersections. For this, take $x,y\in A_{k+1}$ and consider the following cases:\\
\begin{enumerate}
    \item If $x,y\in C_0$, then $S_x^{k+1},S_y^{k+1}\in \mathbb{S}^{k}$. By the recursion hypotheses, we have  $S^{k+1}_x\cap S^{k+1}_y\in \mathbb{S}^{k}\subseteq \mathbb{S}^{k+1}$.
    \item If $x,y\in C_1$, let $x'=\phi^{-1}_k(x)$ and $y'=\phi^{-1}_k(y)$. We now consider two subcases. If $D_x\cap D_y=\emptyset$, this means $S^{k+1}_x\cap S^{k+1}_y=S^k_{x'}\cap S^k_{y'}\in \mathbb{S}^k$. On the other hand, if the intersection of $D_x$ and $D_y$ is nonempty, this means $S^k_{x'}\cap S^k_{y'}=S^k_w$ for some $w\in[r_{k+1},m_k)\times 2^k$. Since $\mathbb{S}^k$ satisfies the conditions (a) and (b),  $D_x\cap D_y=D_{\phi_k(w)}$. Thus, we conclude  $S^{k+1}_x\cap S^{k+1}_y=S^{k+1}_{\phi_k(w)}.$
    \item If $x,y\in C_2$, we can suppose without loss of generality that the first coordinate of $x$ is smaller or equal than the first coordinate of $y.$ In this case, $S^{k+1}_x\cap S^{k+1}_y=S^{k+1}_x$ if $x$ and $y$ share the second coordinate, or $S^{k+1}_x\cap S^{k+1}_y=\emptyset$ if their second coordinates are equal.

    \item If $x\in C_0$ and $y\in C_1$, by definition, we have that $S^{k+1}_x\cap S^{k+1}_y=S^k_x\cap S^k_{\phi^{-1}(y)}$.
    \item If $x\in C_0$ and $y\in C_2$, then $S^{k+1}_x\cap S^{k+1}_y=\emptyset.$
    \item Finally, if $x\in C_1$ and $y\in C_2$, let $a<m_{k+1}$ and $2^k\leq b<2^{k+1}$ such that $y=(a,b)$. If $b-2^k\in D_x$, we have  $S^{k+1}_x\cap S^{k+1}_y=S^{k+1}_{(c,b)}$ where $c=\min\{a,r_{k+1}\}$. On the other side, if $b-2^k\notin D_x$, the intersection is empty.\\
\end{enumerate}
In this way, we finish the proof of the satisfaction of condition (a). By carefully looking at the equalities in cases (1), (2) and (3), we also conclude that condition (d) is true for $\mathbb{S}^{k+1}$.\\ 
The only thing left to do is to check condition (b). For this,  first observe that for each $(a,b)\in C_0\cup C_2$ it is trivially true that $rank(S^{k+1}_{(a,b)})=a+1$. To prove the same holds for all $x=(a,b)\in C_1$, we use induction over the first coordinate. The base case is when $a=m_k$. Here, we have that $S^{k+1}_{x}=S^{k+1}_{(r_{k+1},b)}\cup S^{k+1}_{(m_k-1,b+2^k)}$. Furthermore, any element 
of $\mathbb{S}^{k+1}$ which is contained in $S^{k+1}_{x}$, is also contained in either $S^{k+1}_{(r_k,b)}$ or $ S^{k+1}_{(m_k-1,b+2^k)}$. This means $$rank(S^{k+1}_{(m_k,b)})=\max\{ r_{k+1}+1, m_k+1\}=m_k+1.$$ 
For the induction step, suppose we have proved that for each $(c,d)\in C_1$ with $c<a$. Let $L$ be the set of all $(a-1,d)$ for which $S^{k+1}_{(a-1,d)}$ is (properly) contained in  $S^{k+1}_x$. By the induction hypotheses, it follows that $rank(S^{k+1}_y)=a$ for each $y \in L$. The key is to notice tha every element of $\mathbb{S}^{k+1}$ which is  properly contained in $S^{k+1}_x$ is also contained in $S^{k+1}_y$ for some $y\in L$. Then $rank(S^{k+1}_x)=a+1$ as a consequence of this fact. \\\\
Now, let $\mathcal{F}$ be a $2$-capturing construction scheme of type $\tau$. For each $k\in\omega$, let $f_k:\omega_1\times \omega \longrightarrow \omega^2$ be given as $f_k(\alpha,b)=(|(\alpha)^-_k|,b)$. Now, for each $x\in \omega_1\times \omega$, let $$S_{x}=\bigcup\limits_{k\in\omega}S^k_{f_k(x)}.$$
Finally, let $\mathbb{S}=\{\emptyset\}\cup\{S_x\}_{x\in \omega_1\times\omega}$. We will show that $\mathbb{S}$ is a Suslin lower semi-lattice. As before, it suffices to show it is closed under intersections. In this way, take $x=(\beta,b),y=(\delta,d)\in \omega_1\times \omega$. By conditions (c) and (d),  $S_x\cap N_k=S^k_{f_k(x)}$ and $S_y\cap N_k=S^k_{f_k(y)}$ for each $k\in\omega$. Observe that $f_k(x),f_k(y)\in A_k$ where $k=\max\{\rho(\beta,\delta),b,d\}$. By condition (a), we know  $S^k_{f_k(x)}\cap S^k_{f_k(y)}\in \mathbb{S}^k$. If this intersection is empty, we can use condition (d) to conclude that $S_x\cap S_y$ is also empty. On the other hand, if $S^k_{f_k(x)}\cap S^k_{f_k(y)}=S^k_{(a,c)}$ for $(a,b)\in A_k$, then $a<|(\beta)^-_k|$. Thus, there is $\alpha\in (\beta)_{k}$ for which $|(\alpha)_k^-|=a$. In this way,  $S^k_{f_k(x)}\cap S^k_{f_k(y)}=S^k_{f_k(\alpha,c)}$. From conditions $(c)$ and $(d)$ we conclude $S_x\cap S_y=S_{(\alpha,c)}$. 
Furthermore, condition (b) implies $\mathbb{S}$ is well founded and its rank function satisfies the following properties:\begin{itemize}
    \item $rank(\emptyset)=0,$
    \item $rank(\alpha,b)=\alpha+1$ if $\alpha\in \omega$,
    \item $rank(\alpha,b)=\alpha$ if $\alpha\in \omega_1\backslash \omega.$
\end{itemize}
We do not have to prove that $\mathbb{S}$ has no uncountable chains since it is a substructure of $(\mathscr{P}(N),\subseteq,\cap)$, where $N=\bigcup\limits_{i\in\omega}N_i$.\\ The only thing left to do is to prove  $\mathbb{S}$ is powerful pie. For this, let $n\in\omega$ and $\mathcal{A}$ be an uncountable subset of $\mathbb{S}^n$. For each $x\in \mathcal{A}$, let $(\alpha^x_i,b^x_i)$ be such that $x(i)=S_{(\alpha^x_i,b^x_i)}$ for all $i<n$. Without loss of generality  $\alpha^x_0\leq\dots\leq \alpha^x_{n-1}$. Furthermore, we can suppose $\alpha^x_0>\alpha^y_{n-1}$ or $\alpha^x_{n-1}<\alpha^y_0$ whenever $x$ and $y$ are distinct elements of $\mathcal{A}.$ Finally, we can suppose  there are $b_0,\dots,b_{n-1}\in \omega$ such that $b^x_i=b_i$ for each $i<n$ and $x\in \mathcal{A}$. Since $\mathcal{F}$ is $2$-capturing, we know there are $k>\max\{b_i\,|\,i<n\}$, $F\in \mathcal{F}_k$ and $x,y\in \mathcal{A}$ such that $F$ captures $\{\{\alpha^x_i\}_{i<n},\{\alpha^y_i\}_{i<n}\}$. In this case, $f_k(x(i)),f_k(y(i))\in A_k$ for each $i<n$. 
    Since the first coordinates of $x$ and $y$ form an non decreasing sequence with respect to the index, we conclude  $|(\alpha_i^x)^-_{k-1}|=|\alpha^x_i\cap F_0|=|\alpha^y_i\cap F_1|=|(\alpha_i^y)^-_{k-1}|$.
Since the second coordinates are equal in each value of $i$, this means $\phi_k(f_{k}(x(i)))=f_{k}(y(i))$. \\
By condition (4), $S^{k}_{f_k(x(i))}=S^k_{f_k(y(i))}\cap N_k\subseteq S^k_{f_k(y(i))}$. By a previous argument, this means  $S_{x(i)}\subseteq S_{y(i)}$ for each $i<n$. Consequently, $x$ and $y$ testify that $\mathcal{A}$ is not a pie. 
\end{proof}
\end{theorem}
For the construction given above, we used a particular type. This was done in order facilitate the calculations. Nevertheless, with some minor modifications, one can check that the proof works for each type. 

 In \cite{forcingandconstructionschemes}, it was proved that the existence of $n$-capturing construction schemes does not imply the existence of $(n+1)$-capturing construction schemes. Since the $\Diamond$-principle implies the existence of fully capturing construction schemes, Theorem \ref{suslinlatticescheme} seems to be a step forward to the solution of the following problem.
 \begin{Problem}[\cite{raghavan2014suslin}]Does $ZFC+CH$ imply the existence of a Suslin lower semi-lattice?
\end{Problem}

\subsection{Entangled sets} Now we pass to the study of entangled sets. This objects were introduced by U. Abraham and S. Shelah in \cite{MAdoesnotImplyBA}, as a way of proving that $MA$ does not imply that every two $\aleph_1$-dense sets\footnote{ $A\subseteq \mathbb{R}$ is said to be $\aleph_1$-dense if $|(a,b)\cap A|=\omega_1$ for each $a<b\in\mathbb{R}$.} of reals are isomorphic, principle that is known as the Baumgartner Axiom (abbreviated as $BA(\omega_1))$. Readers interested in learning more about entangled sets may also look at \cite{WhyYcc}, \cite{guzmanentangled}, \cite{AsperoMotaEntangled}, \cite{RemarksonChainConditionsinProducts}, \cite{PartitionProblems} and \cite{EntangledCohen}.
\begin{definition}Let $k\in\omega$, $t:k\longrightarrow\{>,<\}$ and $a=\{a_0,\dots, a_{k-1}\},b=\{b_0,\dots,b_{k-1}\}\in[\mathbb{R}]^k$ enumerated in an increasing way and disjoint. We say that $(a,b)$ realizes $t$ if  $$a_i\;t(i)\; b_i$$ for each $i<k.$ By $T(a,b)$ we denote the unique  $t$ which is realized by $(a,b)$.
\end{definition}
\begin{definition}Let $\mathcal{E}\in [\mathbb{R}]^{\omega_1}$ and $k\in\omega.$ We say that:\begin{itemize}
    \item $\mathcal{E}$ is $k$-entangled if for each  uncountable family $\mathcal{A}\subseteq [\mathcal{E}]^{k}$ of pairwise disjoint sets and  $t:k\longrightarrow\{>,<\}$ there are $a\not=b\in \mathcal{A}$ such that $T(a,b)=t.$
    \item $\mathcal{E}$ is entangled if it is $k$-entangled for each $k\in\omega$.
\end{itemize}
\end{definition}
The following lemma will be implicitly used in the proof of Theorem \ref{entangledscheme}. In an informal way, it says that it is possible to redefine the notion of entangled when the given set $\mathcal{E}$ is presented to us as an indexed family of reals.
\begin{lemma} Let $k\in\omega$ and $\mathcal{E}\in [\mathbb{R}]^{\omega_1}$ injectively enumerated as $\{r_\alpha\}_{\alpha\in \omega_1}$. Then $\mathcal{E}$ is $k$-entangled if and only if for every uncountable family $\mathcal{C}\subseteq [\omega_1]^{k}$ of pairwise disjoint sets and $t:k\longrightarrow\{>,<\}$ there are distinct $c,d\in \mathcal{C}$ for which $$r_{c(i)}\,t(i)\,r_{d(i)}$$ 
for each $i<k.$
\begin{proof}Suppose $\mathcal{E}$ is $k$-entangled. Let $\mathcal{C}\subseteq[\omega_1]^k$ be an uncountable family of pairwise disjoint sets and $t:k\longrightarrow \{>,<\}$.  Given $c\in \mathcal{C}$, let $h_c:k\longrightarrow k$ be unique function satisfying that $i<j$ if and only if $r_{c(h(i))}<r_{c(h_c(j))}$ for every $i,j<k$. By refining $\mathcal{C}$, we can suppose without loss of generality that there is $h$ for which $h=h_c$ for all $c\in \mathcal{C}$. The key observation is that $h$ codes the increasing enumeration (with respect to the reals) of $\{r_{c(i)}\}_{i<k}$ whenever $c\in \mathcal{C}$. That is, $r_{c(h(0))}<\dots<r_{c(h(k-1))}$. Now, let  $t'=t\circ h$. By the previous observation and since $\mathcal{E}$ is $k$-entangled, there are distinct $c,d\in \mathcal{C}$ for which $$T(\{r_{c(h(i))}\}_{i<k},\{r_{d(h(i))}\}_{i<k})=t'.$$ To finish, take an arbitrary $i<k$ and let $i'=h^{-1}(i)$. Then $r_{c(h(i'))}\,t'(i')\, r_{d(h(i'))}$. But $r_{c(h(i'))}=r_{c(i)}$, $t'(i')=t(i)$ and $r_{d(h(i'))}=r_{d(i)}$, which means that $r_{c(i)}\,t(i)\, r_{d(i)}$. So we are done. The proof of the remaining implication is completely similar. 
\end{proof}
\end{lemma}
As a third application we show that the existence of entangled sets follows from the existence of fully capturing construction schemes. 
\begin{definition}\label{definitionlex} Let $(X,<)$ be a linear order. 
Given functions $f,g\in X^{\omega}$, we say that $f<_{lex} g$ if $f(n)<g(n)$ where $n=\min\{k\in\omega\,|\,f(k)\not=g(k)\}$.
\end{definition}
\begin{remark}\label{remarklex} $(X^\omega,<_{lex})$ is a linear order which can be embedded in $\mathbb{R}$, whenever $X$ is a countable set.
\end{remark}
\begin{theorem}[FCA]\label{entangledscheme}There is an entangled set.
\begin{proof}Let $\mathcal{F}$ be a fully capturing construction scheme of type $\{(m_k,n_{k+1},r_{k+1})\}_{k\in\omega}$ 
such that $n_{k+1}\geq 2^{m_k}+1$ for each $k\in\omega.$ By Remark \ref{remarkdelta}, we can think of $(\mathbb{Z}^\omega,<_{lex})$ as a subset of $\mathbb{R}$. In this way, it is enough to define an entangled set
$\mathcal{E}=\{f_\alpha\}_{\alpha\in\omega_1}\subseteq \mathbb{Z}^\omega$.\\\\
First, enumerate $\mathscr{P}(m_{k}\backslash r_{k+1})$ (possibly with repetitions) as $\{C^k_i\}_{0<i<n_{k+1}}$ for each $k\in\omega$. Now, given $\alpha\in \omega_1$ define $f_\alpha:\omega\longrightarrow \mathbb{Z}$ as follows:
$$f_\alpha(k)=\begin{cases}0 &\textit{if }k=0\textit{ or } \Xi_\alpha(k)=-1\\
\Xi_\alpha(k)&\textit{if }\Xi_\alpha(k)\geq 0 \textit{ and }|(\alpha)^-_{k-1}|\in C^{k-1}_{\Xi_\alpha(k)}\\
-\Xi_\alpha(k)&\textit{if }\Xi_\alpha(k)\geq 0\textit{ and }|(\alpha)^-_{k-1}|\notin C^{k-1}_{\Xi_\alpha(k)}\\
\end{cases}$$

To show that $\mathcal{E}$ is entangled, let $k\in\omega$, $t:k\longrightarrow \{>,<\}$ and $\mathcal{C}\subseteq [\omega_1]^k$ be an uncountable 
family of pairwise disjoint sets. Since $\mathcal{F}$ is fully capturing, there are $l>0$, 
$F\in\mathcal{F}_{l}$ and $\{c_1,\dots,c_{n_{l}-1}\}\in[\mathcal{C}]^{n_l}$ captured by $F$. Since the $c_i's$ are pairwise 
disjoint, then $c_i\subseteq F_i\backslash R(F)$ for all $i<n_{l}$.  Notice that $f_{c_i(j)}|_l=f_{c_r(j)}|_l$ for each $i,r<n_l$ and 
$j<k$. Let $0<z<n_{l}$ for which $C^{l-1}_{z}=\{\,|(c_0(j))^-_{l-1}|\,|
\,j<k\textit{ and }t(j)=\,<\,\}$. By definition, we  have that $f_{c_z(j)}(l)=z$ whenever $t(j)=\,<$ and $f_{c_z(j)}(l)=-z$ in case that $t(j)=\,>$. This implies $T(c_0,c_{z})=t,$ so we are done.
\end{proof}
\end{theorem}
As in the case of Theorem \ref{suslintreescheme}, Theorem \ref{entangledscheme} suggests the following question.
\begin{Problem}Does $CA_n(part)$ imply the existence of entangled sets for some $n\in\omega\backslash 1$?
\end{Problem}

\subsection{A result in polychromatic Ramsey theory}  Ramsey theory was born in \cite{ramseyproblemlogic}. In that paper, F. P. Ramsey wrote a lemma which is now known as Ramsey Theorem. Namely, if $k,l\in \omega$ and $c:[\omega]^k\longrightarrow l$ is an arbitrary coloring, there is $A\in [\omega]^{\omega}$ such that $c|_{[A]^k}$ is monochromatic. Classical Ramsey theory is concerned about results of this flavour. In an informal sense, we could say that its goal is to find order in chaos. On the contrary, polychromatic Ramsey theory is about finding chaos. An example of this can be found in \cite{partitioningpairsofcountableordinals}, where the third author proved that there is a coloring $c:[\omega_1]^2\longrightarrow \omega_1$ such that $c[[A]^2]=\omega_1$ for each $A\in [\omega_1]^{\omega_1}$. For this, he first constructed a function $\phi:[\omega_1]^2\longrightarrow \omega_1$ with the property that $\phi[[A]^2]$ contains a closed an unbounded subset of $\omega_1$ for each uncountable $A$. It is not hard to see that if $\mathcal{F}$ is a construction scheme, the function $\phi$ defined as $\phi(\alpha,\beta)=\min (\beta)_{\Delta(\alpha,\beta)}\backslash \alpha$ satisfies this property.

\begin{definition}Let $c:[\omega_1]^2\longrightarrow \omega_1$ be a coloring, $A\subseteq \omega_1$ and $\kappa$ be a (possible finite) cardinal. We say that:\begin{itemize}
    \item $c$ is $\kappa$-bounded if $|c^{-1}[\{\xi\}]|<\kappa$ for each $\xi\in \omega_1.$
    \item $A$ is injective if $c|_{[A]^2}$ is injective.
\end{itemize}
\end{definition}
The problem whether every $2$-bounded coloring $c:[\omega_1]^2\longrightarrow \omega_1$ has an uncountable injective sets was first asked by F. Galvin in the early 1980's, who proved that $CH$ implies a negative answer to that question. In \cite{positivepartition}, the third author prove that it is consistent, and in particular that $PFA$ implies that every $\omega$-bounded $c:[\omega_1]^2\longrightarrow \omega_1$ has an uncountable injective set.\\
In \cite{abrahampolychromatic}, U. Abraham, J. Cummings and C. Smyth proved that it is consistent that there is a $2$-bounded coloring $c:[\omega_1]^2\longrightarrow \omega_1$ without uncountable injective sets in any $ccc$ forcing extension. After hearing this theorem, S. Friedman asked for a concrete example of a $2$-bounded coloring without an uncountable injective set, but which adquires one in a $ccc$ forcing extension ($ccc$-destructible). Such example was produced in \cite{abrahampolychromatic} assuming $CH$ and the existence of a Suslin tree. Here, we construct one using $3$-capturing construction schemes.
\begin{theorem}[$CA_3$]There is a coloring $c:[\omega_1]^2\longrightarrow \omega_1$ with the following properties:
\begin{enumerate}
    \item $c$ is $2$-bounded,
    \item $c$ has no uncountable injective sets,
    \item $c$ is $ccc$-destructible.
\end{enumerate}
\begin{proof}
Let $\mathcal{F}$  be a $3$-capturing construction scheme of an arbitrary type. Let $\psi:\omega_1\times\omega\times\omega\longrightarrow \omega_1$ be a bijection, and define 
$c:[\omega_1]^2\longrightarrow \omega$ as follows:
$$c(\alpha,\beta)=\begin{cases}\psi(\beta,\rho(\alpha,\beta),|(\alpha)_{\rho(\alpha,\beta)}|)&\textit{if }\alpha<\beta\textit{ and }\Xi_\beta(\rho(\alpha,\beta))\geq 3\\
\psi(\beta,\rho(\alpha,\beta),|(\alpha)_{\rho(\alpha,\beta)-1}|)&\textit{if }\alpha<\beta \textit{ and } \Xi_\beta(\rho(\alpha,\beta))<3
\end{cases}$$
To prove  $c$ is a $2$-bounded coloring, let $\xi\in \omega_1$ and suppose that $\{\alpha_0,\beta_0\},\{\alpha_1,\beta_1\}$ and $\{\alpha_2,\beta_2\}$ are elements of $c^{-1}[\{\xi\}]$. We will show two of those sets are equal. Since $\psi$ is one to one, $\beta_0=\beta_1=\beta_2$ and $\rho(\alpha_0,\beta_0)=\rho(\alpha_1,\beta_1)=\rho(\alpha_2,\beta_2)$. Let us name 
those values as $\beta$ and $k$ respectively. Thus, $\rho(\alpha_i,\alpha_j)\leq\max(\rho(\alpha_i,\beta),\rho(\alpha_j,\beta))=k$ for any $i,j<3$. The first case is if $\Xi_\beta(k)\geq 3$. Here $|(\alpha_0)_k|=|(\alpha_1)_k|=|(\alpha_2)_k|$. In particular,  $\Delta(\alpha_0,\alpha_1)>k$. By Remark \ref{remarkdelta}, the only way in which this can happen is if $\alpha_0=\alpha_1$, so we are done. The remaining case is when $\Xi_\beta(k)<2$. Here,  $0\leq \Xi_{\alpha_i}(k)<2$ for each $i<3$. This means there are $i<j<3$ for which $\Xi_{\alpha_i}(k)=\Xi_{\alpha_j}(k)$. By definition of $c$, we also know  $|(\alpha_i)_{k-1}|=|(\alpha_j)_{k-1}|$. In other words, $\Delta(\alpha_i,\alpha_j)\geq k$. Again, by Remark \ref{remarkdelta}, $\alpha_i=\alpha_j$.\\\\
Now, we prove $c$ has no uncountable injective set. For this, let $S\in[\omega_1]^{\omega_1}$. Since $\mathcal{F}$ is $3$-capturing, there is $l\in \omega$, $F\in \mathcal{F}_l$ and $\{\alpha_0,\alpha_1,\alpha_2\}\in [S]^3$ which is captured by $F$. In particular, $\Xi_{\alpha_2}(l)=2$ and for each $i<2$ the following properties hold:
\begin{itemize}
    \item $\rho(\alpha_i,\alpha_2)=l,$
    \item $ |(\alpha_i)_{l-1}|=|(\alpha_2)_{l-1}|.$
\end{itemize}
In other words, $c(\alpha_0,\alpha_2)=c(\alpha_1,\alpha_2)$.\\\\
Finally, we prove $c$ is $ccc$-destructible. For this, let $\mathbb{P}=\{p\in[\omega_1]^{<\omega}\,|\, p\textit{ is injective}\}$ ordered by reverse inclusion. We claim that $\mathbb{P}$ is $ccc$. Since $\mathcal{F}$ is $3$-capturing ( thus, $2$-capturing ), if $\mathcal{A}\in[\mathbb{P}]^{\omega_1}$ there are $p,q\in \mathcal{A}$, $l\in\omega$ and $F\in\mathcal{F}_l$ capturing $p$ and $q$. By definition of $c$, it is easy to see that $p\cup q$ is injective. Hence, $\mathcal{A}$ is not an antichain, and since $\mathcal{A}$ was arbitrary,  $\mathbb{P}$ is $ccc.$ Finally, since $\mathbb{P}$ is $ccc$ and uncountable, there is $p\in \mathbb{P}$ which forces the generic filter to be uncountable (see \cite{Kunen}). This finishes the proof.
\end{proof}
\end{theorem}
\subsection{Coherent families of functions beyond ZFC} There is a trivial way to destroy an  $(\omega_1,\omega_1)$-gap, say $(\mathcal{A},\mathcal{B})$, via forcing. Just force with a forcing notion which colapses $\omega_1$. In the generic extension, $(\mathcal{A},\mathcal{B})$ will be a countable pregap, which implies that it is no longer a gap. Thus, questions about destructibility of $(\omega_1,\omega_1)$-gaps are only interesting for forcing notions which do not colapse $\omega_1.$
\begin{definition}Let $(\mathcal{A},\mathcal{B})$ be an $(\omega_1,\omega_1)$-pregap. We say that $(\mathcal{A},\mathcal{B})$ is destructible if there is a forcing notion $\mathbb{P}$ which preserves $\omega_1$ in such way that $(\mathcal{A},\mathcal{B})$ is not a gap in some generic extension through $\mathbb{P}$. If this does not happen, the gap is said to be indestructible.
\end{definition}
Sometimes, it is possible to destroy one gap while making another one indestructible. Consider the following notion:
\begin{definition}Let $(\mathcal{A},\mathcal{B})$ be an $(\omega_1,\omega_1)$-pregap with $\mathcal{A}=\{A_\alpha\}_{\alpha\in X}$ and $\mathcal{B}=\{B_\alpha\}_{\alpha\in X}$. We define the following forcing notions:
\begin{itemize}
 \item $\chi_0(\mathcal{A},\mathcal{B})=\{\sigma\in [X]^{<\omega}\,|\, \big(\bigcup\limits_{\alpha\in \sigma}A_\alpha \big)\cap \big(\bigcup\limits_{\alpha\in \sigma}B_\alpha\big)=\emptyset\}$
    \item $\chi_1(\mathcal{A},\mathcal{B})=\{\sigma\in [X]^{<\omega}\,|\,\forall \alpha\not=\beta\in \sigma\,\big( (A_\alpha\cap B_\beta)\cup(A_\beta\cap B_\alpha)\not=\emptyset\big)\}$
\end{itemize}
both ordered by reverse inclusion. Whenever $\{f_\alpha\}_{\alpha\in\omega_1\backslash\omega}$ is a coherent family of functions supported by an $\omega_1$-tower, $i\in 2$ and $s\in [\omega_1\backslash\omega]^2$, we denote $\chi_i(\mathcal{A}^{s(0)},\mathcal{A}^{s(1)}) $ simply as $\chi_i(s).$
\end{definition}
It turns out that the forcing notions $\chi_0$ and $\chi_1$ characterize when a pregap is a gap and when it is destructible respectively.
\begin{theorem} [\cite{independenceforanalysts}, \cite{ScheepersGaps}, \cite{PartitionProblems}, \cite{yoriokadestructiblegaps}]Let $(\mathcal{A},\mathcal{B})$ be an $(\omega_1,\omega_1)$-pregap:\begin{itemize}
    \item $\chi_1(\mathcal{A},\mathcal{B})$ is $ccc$ if and only if $(\mathcal{A},\mathcal{B})$ is a gap. In this case, $\chi_1(\mathcal{A},\mathcal{B})$ forces $(\mathcal{A},\mathcal{B})$ to be indestructible.
    \item $\chi_0(\mathcal{A},\mathcal{B})$ is $ccc$ if and only if $(\mathcal{A},\mathcal{B})$ is destructible. In this case, $(\mathcal{A},\mathcal{B})$ is not a gap in some generic extension through $\chi_0(\mathcal{A},\mathcal{B}).$
\end{itemize}
\end{theorem}
\begin{definition}We say that a family $\{(\mathcal{A}_i,\mathcal{B}_i)\}_{i\in I}$ of $(\omega_1,\omega_1)$-gaps is independent if:$$\prod\limits_{c\in C}\chi_{s(c)}(\mathcal{A}_c,\mathcal{B}_c)$$ 
is $ccc$ for each $\mathcal{C}\in [I]^{<\omega}$ and $s:\mathcal{C}\longrightarrow 2.$ Additionally, a coherent family of functions supported by some pretower, say $\{f_\alpha\}_{\alpha\in\omega_1\backslash \omega}$, is said to be independent if $\{(\mathcal{A}^{c(0)},\mathcal{A}^{c(1)})\}_{c\in [\omega_1\backslash \omega]^2}$ is independent.
\end{definition}
The point of the definition above is that if $\{(\mathcal{A}_i,\mathcal{B}_i)\}_{i\in I}$ is an independent family of $(\omega_1,\omega_1)$-gaps, then for every $F:I\longrightarrow 2$, the finite support product $$\prod\limits_{i\in I}^{FS}\chi_{F(i)}(\mathcal{A}_i,\mathcal{B}_i)$$ is $ccc$. This can be used for coding in a similar way that U. Abraham an S. Shelah did in \cite{Shelahabrahamincompactness}.\\

In \cite{yoriokadestructiblegaps}, T. Yorioka proved, assuming the $\Diamond$-principle, that there is an independent family of $2^{\omega_1}$ gaps. In this subsection, we will construct an independent family of gaps fro $FCA$. For this purpose, we will need the following proposition (see \cite{StevoIlias}).
\begin{proposition}\label{independentccc}Let $\{(\mathcal{A}_c,\mathcal{B}_c)\}_{c\in \mathcal{C}}$ be a finite family of $(\omega_1,\omega_1)$-gaps and $s:\mathcal{C}\longrightarrow 2$. If there is no uncountable antichain $\mathcal{A}$ in $\mathbb{P}=\prod_{c\in \mathcal{C}\chi_{s(c)}(\mathcal{A}_c,\mathcal{B}_c)}$ such that $|g(c)|=1$ for all $g\in \mathcal{A}$ and $c\in \mathcal{C}$, then $\mathbb{P}$ is $ccc.$
\end{proposition}

\begin{theorem}[FCA] There is an independent coherent family of functions supported by an $\omega_1$-pretower.
\begin{proof}The construction here is similar to the one in Theorem \ref{hausdorffcoherenttheorem}. Fix a type $\{(m_k,n_{k+1},r_{k+1})\}_{k\in\omega}$ such that $n_{k+1}>2^{r_{k+1}^2}$ for all $k\in\omega$, and let $\mathcal{F}$ be a fully capturing construction scheme of that type. For each $k\in \omega\backslash 1$, define $N_k$ as $\big(\{k\}\times [r_k]^2)\cup \big(\{k\}\times (m_k-r_k)\big)$ and let $N$ be the (disjoint) union of all $N_k$'s. Also enumerate $\mathscr{P}([r_k]^2)$ (possibly with repetitions) as $\{S^k_i\}_{i<n_k}$ in such way that $S^k_0=S^k_1=[r_k]^2$.\\\\We start by defining an $\omega_1$-pretower over $N$.  
Given $\alpha\in \omega_1\backslash \omega$ and $k\in\omega\backslash 1,$ we define $T^k_\alpha$ as follows:\\
\begin{enumerate}
    \item If $\Xi_\alpha(k)=-1$, let $T^k_\alpha=\emptyset.$
    \item If $i=\Xi_\alpha(k)\geq 0$, let $T^k_\alpha=\big(\{k\}\times S_i^k\big)\cup \big(\{k\}\times (|(\alpha)^-_k|-r_k)\big). $\\
\end{enumerate}
Finally, let $T_\alpha=\bigcup_{k\in\omega\backslash 1} T^k_\alpha$. The proof of $\mathcal{T}=\{T_\alpha\}_{\alpha\in \omega_1\backslash \omega}$ being an $\omega_1$-tower is similar to the one in Theorem \ref{hausdorffcoherenttheorem}. In fact, one can show $T^k_\alpha\subseteq T^k_\beta$ whenever $k>\rho(\alpha,\beta)$. The details of this fact are left to the reader.\\

  We are going to define the coherent family of functions supported by $\mathcal{T}$. For this, let $\beta\in\omega_1\backslash \omega$. In order to define $f_\beta:T_\beta\longrightarrow \beta$ consider $\alpha<\beta$ and an arbitrary $k\in\omega$. We will describe the fiber of $\alpha$ in $f_\beta$ restricted to $T^k_\beta$ (i.e., $f^-1_\beta[\{\alpha\}]\cap T^k_\beta)$. If $k<\rho(\alpha,\beta)$, we let $f^{-1}_\beta[\{\alpha\}]\cap T^k_\beta=\emptyset$.  If $k\geq \rho(\alpha,\beta)$, we consider the following cases:
\begin{enumerate}[label=$(\alph*)$]
\item If $\Xi_\alpha(k)\geq 0$, let $f_\beta^{-1}[\{\alpha\}]\cap T^k_\beta=\{(k,|(\alpha)^-_k|-r_k)\}$.
\item If $\Xi_\alpha(k)=-1$ and $\Xi_\beta(k)=0$, let $M^{(\alpha,\beta)}_k=\{s\in S^k_i\,|\, s(0)=|(\alpha)^-_k|\}$ and $$f_\beta^{-1}[\{\alpha\}]\cap T^k_\beta=\{k\}\times M^{(\alpha,\beta)}_k.$$
\item If $\Xi_\alpha(k)=-1$ and $i=\Xi_\beta(k)>0$, let $M^{(\alpha,\beta)}_k=\{s\in S^k_i\,|\, s(1)=|(\alpha)^-_k|\}$  and $$f_\beta^{-1}[\{\alpha\}]\cap T^k_\beta=\{k\}\times M^{(\alpha,\beta)}_k.$$
\item If $\Xi_\beta(k)=-1$, $T^k_\beta=\emptyset$ so there is nothing to do.
\end{enumerate}
It is easy to check that each function is well defined and  $\{f_\alpha\}_{\alpha\in\omega_1\backslash \omega}$ is a family of coherent functions supported by $\mathcal\{T_\alpha\}_{\alpha\in\omega_1\backslash \omega}$. In fact, the reader can verify that for all $\alpha<\beta\in \omega_1\backslash \omega$ and $k>\rho(\alpha,\beta)$ , $f_\alpha|_{T^k_\alpha\cap T^k_\beta}=f_\beta|_{T^k_\alpha\cap T^k_\beta}$.
 In order to prove that our family of functions is independent, let $\mathcal{C}\in[[\omega_1\backslash \omega]^2]^{<\omega}$ be nonempty and $s:\mathcal{C}\longrightarrow 2.$ We will show $$\mathbb{P}=\prod_{c\in \mathcal{C}}\chi_{s(c)}(c)$$
 is $ccc$. For this, take an uncountable $\mathcal{A}\subseteq \mathbb{P}$. We will prove that $\mathcal{A}$ is not an antichain. Due to Proposition \ref{independentccc}, we may assume $g(c)$ has a unique element for each $g\in \mathcal{A}$ and every $c\in C$. In order to keep the proof as clean as possible, we will identify $g(c)$ with its only element. Without loss of generality we can suppose that for each $g,h\in \mathcal{A}$ and $c,v\in \mathcal{C}$,  $g(c)<g(v)$ if and only if $h(c)<h(v)$. Moreover, we can suppose  $g(c)>c(1)$ whenever $c\in \mathcal{C}$.
 Consider $\{\bigcup \mathcal{C}\cup im(g)\,|\, g\in \mathcal{A}\}$ and $k>\max\{\rho(c(0),c(1))\,|\,c\in \mathcal{C}\}$. Since $\mathcal{F}$ is fully capturing, there is $l>k$, $F\in \mathcal{F}_l$ and $g_0,\dots,g_{n_l-1}\in \mathcal{A}$ such that $F$ captures $$\{ \bigcup \mathcal{C}\cup im(g_0),\dots, \bigcup \mathcal{C}\cup im(g_{n_l-1})\}.$$ Let $S=\{\{|(c(0))^-_l|,|(c(1))^-_l|\}\,|\,c\in \mathcal{C}\textit{ and }s(c)=1\}$ and $0<i<n_l$ for which $S=S^l_i$.  Recall that $g_0(c)<g_0(v)$  if and only if $g_i(c)<g_i(v)$ whenever $c,v\in \mathcal{C}$. As a consequence of this,  $\rho(g_0(c),g_i(c))=\Delta(g_0(c),g_i(c))$ for each $c\in 
 \mathcal{C}$. In particular, this means  $$f^{-1}_{g_0(c)}[\{c(0)\}]\cap f^{-1}_{g_i(c)}[\{c(1)\}]=\big(f^{-1}_{g_0(c)}[\{c(0)\}]\cap f^{-1}_{g_i(c)}[\{c(1)\}]\big)\cap T^l_{g_0(c)}.$$
Since $\Xi_{g_0(c)}(l)=0$ and $\Xi_{c(0)}(l)=-1$, we know $$f^{-1}_{g_0(c)}[\{c(0)\}]\cap T^l_{g_0(c)}=\{l\}\times \{s\in [r_k]^2\,|\,s(0)=|(c(0))^-_l|\}.$$ This is because $S^l_0=[r_l]^2$. On the other hand $$f^{-1}_{g_0(c)}[\{c(0)\}]\cap T^l_{g_i(c)}=\{l\}\times \{s\in S\,|\, s(1)=|(c(1))^-_l|\}$$  
Thus, $f^{-1}_{g_0(c)}[\{c(0)\}]\cap f^{-1}_{g_i(c)}[\{c(1)\}]\subseteq \{|(c(0))^-_l|,|(c(1))^-_l|\}$. Furthermore, this intersection is nonempty if and only if $s(c)=1$. This means  $g_0$ and $g_i$ are compatible in $\mathcal{P}$. Thus, the proof is over.
 \end{proof}
\end{theorem}
\section{Oscillation theory on construction schemes}\label{oscilationsection}
In \cite{PartitionProblems}, the third author developed a very powerful \say{oscillation theory} and deduced very interesting theorems from it. Usually, the oscillation theory is based on an unbounded family of functions (although there are other variations, see for example \cite{oscillationsintegers}). Here we will develop an oscillation theory base on a bounded family of functions defined from a $2$-capturing construction scheme.\\
For the rest of this section we fix a type $\{(m_k,n_{k+1},r_{k+1})\}_{k\in\omega}$ and $\mathcal{F}$ a $2$-capturing construction scheme of that type. Given $\alpha\in\omega_1$, we define $f_\alpha:\omega\longrightarrow \omega$ as:
$$f_\alpha(l)=|(\alpha)_l|.$$
Finally, we define $\mathcal{B}_\mathcal{F}$ as $\{f_\alpha\}_{\alpha\in\omega_1}.$  The following lemma is a direct consequence of the definitions of $\rho$ and $\Delta$.
\begin{lemma}\label{lemmaf}Let $\alpha<\beta\in \omega_1$. Then:\begin{enumerate}[label=$(\arabic*)$]
\item $f_\alpha(i)=f_\beta(i)$ if $i<\Delta(\alpha,\beta)$,
\item$f_\alpha(j)<f_\beta(j)$ whenever  $j\geq \rho(\alpha,\beta),$
\item $f_\alpha<f_\beta$ provided that $\Delta(\alpha,\beta)=\rho(\alpha,\beta)$
\item In particular, $f_\alpha<^*f_\beta$.
\end{enumerate}
Furthermore, $\mathcal{B}_\mathcal{F}$ is bounded by the function in $\omega^\omega$ which sends each $i$ to $m_i.$
\end{lemma}
It is interesting that even though $\mathcal{B}_\mathcal{F}$ is bounded, its oscillation theory mirrors the oscillation theory of \cite{PartitionProblems} for unbounded families. Since $\mathcal{F}$ is 2-capturing, given any $\mathcal{A}\in [\mathcal{B}_\mathcal{F}]^{<\omega}$, there are $\alpha<\beta\in \mathcal{A}$ with $\Delta(\alpha,\beta)=\rho(\alpha,\beta).$ Thus, we have the following corollary.
\begin{corollary}\label{corollaryBuncountableantichains}$(\mathcal{B}_\mathcal{F}, \leq )$ has no uncountable pies\footnote{Recall that by pie, we mean a set of pairwise incomparable elements}.
\end{corollary}
\begin{definition}Let $f,g\in\omega^\omega$ and $k\in\omega$. We define the following objects: $$\overline{osc}_k(f,g)=\{s\in\omega\backslash k\,|\,f(s)\leq g(s)\textit{ and }f(s+1)>g(s+1)\,\},$$
$$osc_k(f,g)=|\overline{osc}_k(f,g)|. $$

\end{definition}
Given $\alpha,\beta\in \omega_1$ and $k\in\omega $ we will write $osc_k(\alpha,\beta)$ and $\overline{osc}_k(\alpha,\beta)$ instead of $osc_k(f_\alpha,f_\beta)$ and $\overline{osc}_k(f_\alpha,f_\beta)$ respectively. These two objects will be written as $osc(\alpha,\beta)$ and $\overline{osc}(\alpha,\beta)$ whenever $k=0.$\\\\
In the following theorem, let $osc_k(a,b)=\{osc_k(\alpha,\beta)\,|\,\alpha\in a\textit{ and }\beta\in b\}$ whenever $a,b\in[\omega_1]^{<\omega}.$
\begin{proposition}\label{fulloscilation}Let $n,k\in\omega$ and $\mathcal{A}\in[\omega_1]^{n}$ be an uncountable family of pairwise disjoint sets such that  $\rho^a=k$ for each $a\in \mathcal{A}$. Given $l\in\omega$, there are $a<b\in \mathcal{A}$ such that $osc_k(a,b)\subseteq [l,2l].$
\begin{proof}The proof is by induction over $l.$\\\\
{(Base step)} If $l=0$, let $a<b\in \mathcal{A}$ and $F\in \mathcal{F}$ be such that $F$ captures $\{a,b\}.$ Observe that $\rho^F> k$. The following properties hold for all $i,j<n$ due to the capturing assumption on $F$ and by Lemma \ref{lemmaf}:
\begin{enumerate}
    \item $f_{a(i)}|_{[k,\rho^F)}=f_{b(i)}|_{[k,\rho^F)},$
    \item $f_{a(i)}|_{\omega\backslash\rho^F}<f_{b(j)}|_{\omega\backslash \rho^F},$
    \item $f_{a(i)}|_{[k,\rho^F)}<f_{a(j)}|_{[k,\rho^F)}$ provided that $i<j$.
\end{enumerate}
From these facts, we conclude that  $osc_k(a(i),b(j))=0$ for all $i,j<n.$\\\\
{(Induction step)} Suppose we have proved the proposition for some $l\in\omega$. Using the induction hypothesis, we can recursively construct an uncountable $\mathcal{C}\subseteq [\mathcal{A}]^{2}$ of pairwise disjoint sets such $osc_k(a,b)\subseteq [l,2l]$  for each $\{a,b\}\in \mathcal{C}$  with $a<b$. Since $\mathcal{F}$ is $2$-capturing, we can find an uncountable family $\mathcal{D}\subseteq [\mathcal{C}]^2$ and $r\in\omega\backslash(k+1)$ with the following properties:
\begin{itemize}
\item Whenever $\{\{a,b\},\{c,d\}\}\in \mathcal{D}$, there is $F\in \mathcal{F}_r$ which captures $\{a\cup b, c\cup d\}.$ In particular, this implies that $r>\rho^{a\cup b}$ and $a\cup b<c\cup d.$
\item For each $x,y\in D$, $\bigcup x\cap \bigcup y=\emptyset.$
\end{itemize}
Using once again that $\mathcal{F}$ is $2$-capturing, we can get $F\in \mathcal{F}$ with $\rho^F>r$ and two elements of $ \mathcal{D}$, say $\{\{a_0,b_0\},\{c_0,d_0\}\},\{\{a_1,b_1\},\{c_1,d_1\}\}$, for which  $F$ captures $$\{a_0\cup b_0\cup c_0\cup d_0, a_1\cup b_1\cup c_1\cup d_1\}.$$
We claim that $osc(c_0,b_1)\subseteq [l+1,2(l+1)]$. For this, take $i,j<n$.  The following properties hold:
\begin{enumerate}
\item $f_{b_1(j)}|_{[k,\rho^F)}=f_{b_0(j)}|_{[k,\rho^F)}$,
\item $f_{c_0(i)}|_{[k,r)}=f_{a_0(i)}|_{[k,r)},$
\item $f_{b_1(j)}|_{\omega\backslash \rho^F}>f_{c_0(i)}|_{\omega\backslash \rho^F}$,
\item $f_{b_0(j)}|_{[r,\rho^F)}<f_{c_0(i)}|_{[r,\rho^F)}$,
\end{enumerate} 
We use these properties to calculate the oscillation. First observe that  we can use the part $(2)$ of Lemma \ref{lemmaf} to conclude $\overline{osc}_k(a_0(i),b_0(j))\subseteq [k,\rho^{a_0\cup b_0})$ and $\overline{osc}_k(c_0(i),b_1(j))\subseteq [k,\rho^F)$. By the properties (1) and (2) written above and since $r>\rho^{a_0\cup b_0}$,  $$\overline{osc}_k(c_0(i),b_1(j))\cap [k,r-1)=\overline{osc}_k(a_0(i),b_0(j)).$$ Due to properties (1), (3) and (4) we also have that $\rho^F-1\in \overline{osc}_k(c_0(i),b_1(j))$. In fact, properties (1) and (4) also imply that $\rho^F-1$ is the only element in the interval $[r,\rho^F)$ which belong to $\overline{osc}_k(c_0(i),b_1(j)).$ By joining all the previous observations, we get: $$\overline{osc}_k(a_0(i),b_0(j))\cup\{\rho^F-1\}\subseteq \overline{osc}_k(c_0(i),b_1(j))\subseteq \overline{osc}_k(a_0(i),b_0(j))\cup\{\rho^F-1\}\cup\{r-1\}.$$
Which means that $l+1\leq osc_k(c_0(i),b_1(j))\leq 2l+2$. This completes the proof.
\end{proof}
\end{proposition}
By a careful analysis of the argument of the preceding theorem, one can show that whenever $\mathcal{A}\in [\omega_1]^{\omega}$ and $l\in\omega$, then there are $\alpha<\beta\in \mathcal{A}$ for which $osc(\alpha,\beta)=l$. Unfortunately, this property does not hold for arbitrary uncountable families of finite sets. Nevertheless, the previous result is enough to redefine a \say{corrected} oscillation.\\\\
The following lemma is easy.
\begin{lemma}\label{lemmapartition}There is a partition $\{P_n\}_{n\in\omega}$ of $\omega$ such that for every $k,n\in\omega$ there is $l\in\omega$ such that $[l,2l+k]\subseteq P_n.$
\end{lemma}
\begin{theorem}[$CA_2$]\label{coloringca2}There is a coloring $o:[\omega_1]^2\longrightarrow \omega $ such that for every uncountable family $\mathcal{A}\subseteq [\omega_1]^{<\omega}$ of pairwise disjoint sets and each $n\in\omega$, there are $a<b\in \mathcal{A}$ for which $\{o(\alpha,\beta)\,|\alpha\in a\textit{ and }\beta\in b\}=\{n\}$.
\begin{proof}
Let $\{P_n\}_{n\in\omega}$ be a partition of $\omega$ as in Lemma \ref{lemmapartition}. Let $o:[\omega_1]^2\longrightarrow \omega$ be defined as:$$o(\alpha,\beta)=n\textit{ if and only if }osc(\alpha,\beta)\in P_n.$$
We claim  that $o$ satisfies the conclusion of the theorem. Indeed, let $\mathcal{A}$ be an uncountable family of pairwise disjoint finite subsets of $\omega_1$ and $n\in\omega$. By refining $\mathcal{A}$ we may suppose that there is $k\in\omega$ such that $\rho^a=k$ for every $a\in \mathcal{A}$. Let $l\in\omega$ be such that $[l,2l+k]\subseteq P_n$. Due to Proposition \ref{fulloscilation}, there are $a<b\in \mathcal{A}$ such that $osc_k(a,b)\subseteq[l,2l].$ Given $\alpha\in a$ and $\beta\in b$, it is easy to see that $\overline{osc}(\alpha,\beta)$ has at most $k$ more elements than $\overline{osc}_k(\alpha,\beta)$. In this way, $osc(\alpha,\beta)\in [l,2l+k]\subseteq P_n$. In other words, $o(\alpha,\beta)=n$. This finishes the proof.
\end{proof}
\end{theorem}
The existence of a coloring with the properties stated above, already implies the existence of a much more powerful coloring. As we shall mention later, such a coloring can be used to build topological spaces with important properties. 
\begin{corollary}[$CA_2$]\label{coloringo*} There is a coloring $o^*:[\omega_1]^2\longrightarrow \omega$ such that for all $n\in\omega$, $h:n\times n\longrightarrow \omega$ and any uncountable family $\mathcal{A}\subseteq[\omega_1]^{n}$ of pairwise disjoint sets, there are $a<b\in \mathcal{A}$ for which $$o^*(a(i),b(j))=h(i,j)$$
for all $i,j<n.$
\begin{proof}Let $\{h_n\}_{n\in\omega}$ be an enumeration of all $h:X\longrightarrow \omega$ for which $X\subseteq \omega^{<\omega}$ is finite and its elements are pairwise incomparable. Let us call $X_n$ the domain of $h_n$. Note that for each $f\in\omega^\omega$ and every $n\in\omega$ there is at most one $\sigma \in X_n$ which is extended by $f.$ Take a coloring $o$ satisfying the conclusion of Theorem \ref{coloringca2}. We define $o^*:[\omega_1]^2\longrightarrow \omega$ as follows: Given distinct $\alpha,\beta\in \omega_1$, if there are $\sigma_\alpha,\sigma_\beta\in X_{o(\alpha,\beta)}$ for which $\sigma_\alpha\subseteq f_\alpha$ and $\sigma_\beta\subseteq f_\beta$, put $$o^*(\alpha,\beta)=h_{o(\alpha,\beta)}(\sigma_\alpha,\sigma_\beta).$$
In any other case, let $o^*(\alpha,\beta)=17$.
In order to prove that $o^*$ satisfies the conclusion of the corollary, let $n\in\omega$, $h:n\times n\longrightarrow \omega$ and  $\mathcal{A}\subseteq[\omega_1]^n$ be an uncountable family of pairwise disjoint sets. By refining $\mathcal{A}$ we may suppose there is $k\in\omega$ with the following properties:\begin{enumerate}[label=$(\arabic*)$]
    \item $\forall a\in \mathcal{A}\,\forall i\not=j<n\, (f_{a(i)}|_k\not=f_{a(j)}|_k)$,
    \item $\forall a,b\in \mathcal{A}\,\forall i<n\,(f_{a(i)}|_k=f_{b(i)}|_k).$

\end{enumerate}
Fix $a_0\in \mathcal{A}$. Let $X=\{f_{a_0(i)}|_k\,|\,i<n\}$ and define $h:X\times X\longrightarrow \omega$ as:
$$h(f_{a_0(i)}|_k,f_{a_0(j)}|_k)=h(i,j).$$
We know there is $m\in\omega$ for which $X=X_m$ and $h=h_m$. For such $m$, there are $a<b\in \mathcal{A}$ such that $o(a(i),b(j))=m$ for all $i,j<n$. For all such $i$ and $j$, we have that $f_{a_0(i)}|_k\subseteq f_{a(i)}$ and $f_{a_0(j)}|_k\subseteq f_{a(j)}$. In this way, $o^*(a(i),b(j))=h(f_{a_0(i)}|_k,f_{a_0(j)}|_k)=h(i,j)$. So we are done.

\end{proof}
\end{corollary}
As an application, we get the following:
\begin{corollary}[$CA_2$]\label{cccnotproductive}ccc is not productive.
\begin{proof}Let $o$ be a coloring of Theorem \ref{coloringca2}. For each $n\in\omega$, let $\mathbb{P}_n=\{p\in [\omega_1]^{<\omega_1}\,|\,\forall\alpha,\beta\in p\,(\textit{if }\alpha\not=\beta\textit{ then }o(\alpha,\beta)=n\,)\}$. In particular, $\mathbb{P}_0$ and $\mathbb{P}_1$ are $ccc$ but $\mathbb{P}_0\times\mathbb{P}_1$ is not. The set $\{(\alpha,\alpha)\,|\,\alpha\in\omega_1\}$ testifies this last fact.
    
\end{proof}
    
\end{corollary}
    
\subsection{Tukey order}We say that a partial order $(D,\leq)$ is (upwards) directed if for every $x,y\in X$ there is $z\in D$ bigger than $x$ and $y$.
\begin{proposition} $(\mathcal{B}_\mathcal{F},\leq)$ is directed.
\begin{proof}Let $\alpha<\beta\in \omega_1$ and let $F\in \mathcal{F}_{\rho(\alpha,\beta)}$ be such that $\{\alpha,\beta\}\subseteq F$. To finish, just observe that if $\delta=\max F$ and $s\leq \rho(\alpha,\beta)$, then $F=\bigcup\{G\in\mathcal{F}_s\,|\,G\subseteq F\}$, so there is $G\in \mathcal{F}_s$ such that $G\subseteq F$ and $\delta\in G$, and it must happen that $\delta=\max G$. This implies that $f_\delta(s)=|(\{\delta\})^-_s|=|(\delta+1)\cap G|=m_s.$ Hence, $f_\delta\geq f_\alpha,f_\beta.$
\end{proof}
\end{proposition}
Recall the following notion for comparing directed partial orders.
\begin{definition}[\cite{Tukey}]Let $(D,\leq_D)$ and $(E,\leq_E)$ be directed partial orders. We say that $E$ es Tukey below $D$, and write it as $E\leq_T D$ if there is $\phi:D\longrightarrow E$ such that $\phi[X]$ is cofinal in $E$ for each cofinal $X\subseteq D.$ Furthermore, we say that $E$ is Tukey equivalent to $D$, and write it  as $E\equiv_T D$, if $E\leq_T D$ and $D\leq_T E$. 
\end{definition}
The study of Tukey order was iniciated by J. Tukey in \cite{Tukey}. Among other things, he proved that the sets $1$, $\omega$, $\omega_1$, $\omega\times \omega_1$ and $[\omega_1]^{<\omega}$ are non Tukey equivalent when equipped with their natural orderings. In \cite{CategorycofinaltypesII}, J. R. Isbell showed that under $CH$, there is at least one directed partial order of cardinality $\omega_1$ which is non Tukey equivalent to any of the previous mentioned. He later improved his result in \cite{sevencofinaltypes}.  In \cite{PartitionProblems}, the third author proved the existence of such a directed partial order under $\mathfrak{b}=\omega_1$. In \cite{directedsetscofinaltypes}, he proved that consistently every directed partial order of cardinality $\omega_1$ is Tukey equivalent to one of the first five we mentioned. From now on, we will call such an order, a sixth Tukey type. The reader interested in learning more about the Tukey ordering and related topics is invited to search for \cite{surveytukeytheory}, \cite{tukeytypesofultrafilters}, \cite{tukeyorderamong}, \cite{Combinatoricsoffiltersandideals}, \cite{analyticidealscofinal}, \cite{CombinatorialDichotomiesandCardinalInvariants}, \cite{idealcompacttukey},  \cite{cofinaltopological} and\cite{avoidingfamilies}.
\begin{proposition}[\cite{Tukey}]Let $(D,\leq_D)$ and $(E,\leq_E)$ be directed partial orders.\begin{itemize}
    \item $E\leq_T D$ if and only if there is $\phi:E\longrightarrow D$ such that $\phi[X]$ is unbounded in $D$ for each unbounded $X\subseteq E.$
    \item $E\equiv_T D$ if and only if there is a partially ordered set $C$ in which both $D$ and $E$ can be embedded as cofinal subsets.
\end{itemize}
\end{proposition}
Let $(D,\leq)$ be a directed partial order. We say $S\subseteq D$ is $\omega$-bounded if every countable subset of $S$ is bounded in $D.$
\begin{proposition}\cite{directedsetscofinaltypes}Let $(D,\leq)$ be a directed set with $|D|=\omega_1.$ Then:
\begin{enumerate}[label=$(\arabic*)$]
    \item $D\leq_T 1$ if and only if $D$ has a greatest element.
    \item $D\leq_T \omega$ if and only if $cof(D)\leq \omega.$
    \item $D\leq \omega_1$ if and only if $D$ is $\omega$-bounded.
    \item $D\leq_T \omega\times\omega_1$ if and only if $D$ can be covered by countably many $\omega$-bounded sets.
    \item $[\omega_1]^{<\omega}\leq_T D$ if and only if there is $A\in[D]^{\omega_1}$ for which every  $X\in[A]^{\omega}$  is unbounded in $D.$
\end{enumerate}
\end{proposition}
As a consequence we get:
\begin{corollary}\label{corollarytukey}Let $(D,\leq)$ be a directed set with $|D|=\omega_1.$ Then either $D\equiv_T 1$, $D\equiv_T \omega$, $D\equiv_T \omega_1$ or $\omega\times \omega_1\leq_T D\leq_T [\omega_1]^{\omega}$.
\end{corollary}
Before we continue, we recall some notation.
\begin{definition}Let $\alpha,\beta$ and $\gamma$ be ordinals. $\gamma\rightarrow (\alpha,\beta)^2_2$ stands for the following statement:\begin{center}
 For all $c:[\gamma]^2\longrightarrow 2$, there is a $0$-monochromatic subset of $\gamma$ of order type $\alpha$ or there is a $1$-monochromatic subset of $\gamma$ of order type $\beta.$ 
\end{center}
Its negation is written as $\gamma\not\rightarrow (\alpha,\beta)^2_2$.
\end{definition}
The following  theorem is due to P. Erd\"os and R. Rado (see \cite{partitioncalculus}). It is a generalization of a well known theorem of B. Dushnik, E. Miller and P. Erd\"os (see \cite{partiallyorderedsets}). The reader can find a proof in \cite{Jech} and \cite{Kunen}.
\begin{theorem}[\cite{partitioncalculus}]\label{erdosdusnik} $\omega_1\rightarrow (\omega_1,\omega+1)^2_2$.
\end{theorem}
In \cite{positivepartition}, the third author proved that it is consistent that it is consistent to have $\omega_1\rightarrow (\omega_1,\alpha)^2_2$ for each $\alpha<\omega_1$. On the other side, he prove in \cite{PartitionProblems} that $\mathfrak{b}=\omega_1$ implies that $\omega_1\not\rightarrow(\omega_1,\omega+2)^2_2$. In the following theorem, we will show that the same is true under $CA_2$. Although a coloring testifying the negation of the previous partition relation can be extracted from some of our other constructions, we decided to include a direct proof due to its simplicity.
\begin{theorem}[$CA_2$] $\omega_1\not\rightarrow (\omega_1,\omega+2)^2_2$. 
\begin{proof}Let $\mathcal{F}$ be a $2$-capturing construction scheme of an arbitrary type. Let $c:[\omega_1]^2\longrightarrow 2$ defined as:
$$c(\alpha,\beta)=\begin{cases}1 & \textit{ if }\Delta(\alpha,\beta)=\rho(\alpha,\beta)\\
0 &\textit{ otherwise}  
\end{cases}$$
 Since $\mathcal{F}$ is $2$-capturing, it is easy to see that there are no uncountable $0$-monochromatic subsets of $\omega_1$. Suppose towards a contradiction that there is a $1$-monochromatic set, say $X$, of order type $\omega+2$. Let $\beta$ and $\gamma$ be the last two elements of $X$. and consider $\alpha\in X\backslash (\gamma)_{\rho(\beta,\gamma)}$. Since $\rho$ is an ordinal metric, $\rho(\alpha,\beta)=\rho(\alpha,\gamma)$. In this way, $\Delta(\beta,\delta)=\rho(\beta,\delta)< \rho(\alpha,\beta)=\Delta(\alpha,\beta)$. By Lemma \ref{countrymanlemma3}, we conclude that $\Delta(\alpha,\delta)=\Delta(\beta,\delta)$. This is a contradiction to $X$ being $1$-monochromatic, so we are done.
\end{proof}
\end{theorem}
Now, we show the existence of sixth Tukey type. 
\begin{theorem}$(\mathcal{B}_\mathcal{F},\leq)$ is a sixth Tukey type.
\begin{proof}By Corollary \ref{corollarytukey}, it is enough to show $\mathcal{B}\not\leq_T\omega\times\omega_1$ and $[\omega_1]^{<\omega}\not\leq_T \mathcal{B}_\mathcal{F}$. First we prove $\mathcal{B}_\mathcal{F}$ does not contain any uncountable $\omega$-bounded set. We argue by contradiction.  Assume there is $A\in[\omega_1]^{\omega_1}$ for which $\{f_\beta\}_{\beta\in A}$ is $\omega$-bounded. Recursively,  we build $\{(\alpha_\xi,\beta_\xi)\}_{\xi<\omega_1}$ satisfying:\begin{enumerate}[label=$(\arabic*)$]
    \item $\{\alpha_\xi\}_{\xi\in\omega_1}$ and $\{\beta_\xi\}_{\xi\in\omega_1}$ are increasing,
    \item $\forall \xi\in \omega_1\, (\,\beta_\xi<\alpha_\xi\,),$ 
    \item $\{\alpha_\xi\}_{\xi\in\omega_1}\subseteq A,$
    \item $\forall \xi\in\omega_1\,(\, f_{\beta_\xi}\textit{ bounds }\{f_{\alpha_\nu}\,|\,\nu<\xi\}\,).$\\
\end{enumerate}
Since $\mathcal{F}$ is $2$-capturing, we can find $l\in\omega$, $F\in\mathcal{F}_{l+1}$ and $\delta<\gamma$ such that $F$ captures $\{\beta_\delta,\alpha_\delta\},\{\beta_\gamma,\alpha_\gamma\}$. It follows that $f_{\beta_\gamma}(l)=f_{\beta_\delta}(l)<f_{\alpha_\delta}(l).$ But this is a contradiction since $f_{\beta_\gamma}$ was supposed  to bound $f_{\alpha_\delta}$.\\
To prove $[\omega_1]^{<\omega}\not\leq_T \mathcal{B}_\mathcal{F}$, let $A\in[\omega_1]^{\omega_1}$. We need to show that $\{f_\alpha\}_{\alpha\in A} $ contains a countable bounded set. Define $d:[A]^2\longrightarrow 2$ as:$$d(\alpha,\beta)=\begin{cases}0 & \textit{ if } f_\alpha\not\leq f_\beta\textit{ and } f_\beta\not\leq f_\alpha\\ 1 & \textit{ otherwise }
\end{cases}$$
By Theorem \ref{erdosdusnik}, there are two possibilities :\begin{enumerate}[label=$(\alph*)$]
    \item $A$ contains a $0$-monochromatic uncountable set.
    \item $A$ contains a $1$-monochromatic set of order type $\omega+1.$
\end{enumerate}
Every $0$-monochromatic set is an antichain in $\mathcal{B}_\mathcal{F}$, so by Corollary \ref{corollaryBuncountableantichains} there can not be uncountable  $0$-monochromatic sets. Hence, there is  $1$-monochromatic subset of $A$, say $X$, of order type $\omega+1$. Observe that $f_\beta$ bounds $\{f_\alpha\}_{\alpha\in X}$ where $\beta=\max X.$
\end{proof}
\end{theorem}
\begin{corollary}[$CA_2$]There is a sixth Tukey type. 
\end{corollary}
\begin{problem}Suppose $\mathcal{F}$ and $\mathcal{G}$ are two $2$-capturing construction. What can we say about the Tukey relation between $\mathcal{B}_\mathcal{F}$ and $\mathcal{B}_\mathcal{G}$? What if they have the same type?
\end{problem}
\begin{problem} Does $CA_2$ imply the existence of $2^{\omega_1}$ non equivalent Tukey types of size $\omega_1$?
\end{problem}
By a result on \cite{CombinatorialDichotomiesandCardinalInvariants}, we can improve a theorem of \cite{lopezschemethesis} by means of the following corollary.
\begin{corollary}$PID+\min(\mathfrak{b},cof(\mathcal{F}_\sigma))>\omega_1$ implies there are no $2$-capturing construction schemes.
    
\end{corollary}
Above, by $PID$ we denote the \say{$P$-ideal Dichotomy}, which is a very strong dichotomy introduced by the third author. For more information see \cite{PartitionPropertiesCompatiblewithCH}, \cite{Lspacespid}, \cite{CombinatorialDichotomiesandCardinalInvariants}, \cite{NotesonForcingAxioms} and \cite{ADichotomyforPIdealsofCountableSets}.
\subsection{Suslin Towers}Let $\mathcal{T}$ be an $\kappa$-pretower. We say $\mathcal{T}$ is Suslin if for every  uncountable $\mathcal{A}\subseteq \mathcal{T}$ there are distinct $A,B\in \mathcal{A}$ with $A\subseteq B$. Suslin towers were studied in \cite{GapsandTowers}. There, P. Borodulin-Nadzieja and D. Chodounsk\'y proved, in particular, that Suslin $\omega_1$-pretowers 
exist under $\mathfrak{b}=\omega_1.$ In fact, whenever $\mathcal{B}$ is an increasing family of functions in $\omega^\omega$ with respect to $<^*$ of order type $\omega_1$, then the family $\{T_f\}_{f\in \mathcal{B}}$, where $$T_f=\{(n,m)\,|\,m\leq f(n)\}$$
is an $\omega_1$-pretower. Also notice that if $f,g\in \mathcal{B}$ are such that $f\leq g$, then $T_f\subseteq T_g$. Consequently, if $(\mathcal{B},\leq)$ has no uncountable pies then $\{T_f\}_{f\in \mathcal{B}}$ is Suslin. Thus, we have the following corollary.
\begin{corollary}[$CA_2$]There is a Suslin $\omega_1$-pretower.   
\end{corollary}
The previous contradicts contradicts the Open Graph Axiom ($OGA$). As with the case of $PID$, $OGA$ (also known as $OCA$, $OCA_{[T]}$ or  Tordorcevic Axiom $TA$) is an axiom introduced by the third author in \cite{PartitionProblems}. The reader may consult \cite{OCAandTowers}, \cite{OCAPmax}, \cite{OpenColoringAxioms}, \cite{TheProperForcingAxiom}, \cite{ProofofNogura}, \cite{CombinatorialDichotomiesinSetTheory} and \cite{StevoIlias} for more information about $OGA$.
\subsection{S-spaces} The results contained in this subsection as well as their proofs are completely based in \cite{PartitionProblems}. For that reason, most of the proofs will be omitted here.
\begin{definition}Let $(X,\tau)$ be a topological space with $|X|=\omega_1$ and $\{x_\alpha\}_{\alpha\in\omega_1}$ be an enumeration of $X.$ We say that:
\begin{itemize}
    \item $X$ is left-open (right-separated) if $\{x_\xi\,|\,\xi\leq \alpha\}$ is open for every $\alpha\in \omega_1$.
    \item  $X$ is an $S$-space if it is $T_3$, hereditarily separable and not Lindel\"of.
    \item $X$ is a strong $S$-space if all of its finite powers are $S$-spaces.
\item $X$ is right-open (left-separated) if $\{x_\xi\,|\,\xi\geq \alpha\}$ is open for every $\alpha\in \omega_1$.
\item  $X$ is an $L$-space if it is $T_3$, hereditarily Lindel\"of and not separable.
\item $X$ is a strong $L$-space if all of its finite powers are $L$-spaces.

\end{itemize}
\end{definition}
The existence of an $S$-space used to be one of the main open problems in set-theoretic topology. Such spaces exist under a large variety of axioms (like $CH$ and some parametrized diamonds of \cite{ParametrizedDiamonds}). This question was finally settled when the third author proved that the Proper Forcing Axiom (PFA) implies that there are no $S$-spaces. Here we will construct $S$-spaces using $2$-capturing construction schemes. To learn more about $S$-spaces (and $L$-spaces) the reader may consult \cite{martinaxiomfirstcountable}, \cite{basicsandl}, \cite{reformulationsandl} and \cite{PartitionProblems}. The following is known:
\begin{lemma}\label{lemmaspacequivalence}Let $(X,\tau)$ be a topological space with $|X|=\omega_1$ and $\{x_\alpha\}_{\alpha\in\omega_1}$ be an enumeration of $X.$ If:
\begin{enumerate}[label=$(\alph*)$]
    \item $X$ is $T_3$,
    \item $X$ is left open,
    \item $X$ does not have uncountable discrete sets.
    \end{enumerate}
Then $(X,\tau)$ is an $S$-space.
\end{lemma}
Let $o^*$ be a coloring satisfying the conclusion of Corollary \ref{coloringo*}. For each $\alpha\in\omega_1$ define $x_\alpha:\omega_1\longrightarrow 2$ as follows:
$$x_\alpha(\beta)=\begin{cases}\min(o^*(\alpha,\beta),1) &\textit{ if }\alpha<\beta\\
0 &\textit{ if }\alpha>\beta\\
1 &\textit{ if }\alpha=\beta
\end{cases}$$
Consider $X=\{x_\alpha\}_{\alpha\in\omega_1}$ endowed with the product topology inherited by $2^{\omega_1}$. Then it is clear that $X$ is $T_3$. Furthermore, for each $\alpha\in \omega_1$, the set $\{x_\beta\,|\,\beta\in \omega_1\textit{ and }x_\beta(\alpha)=1\}$ is an open set contained in $\{x_\xi\,|\xi\leq \alpha\}$ and having $\alpha$ as an element. Hence, $X$ is also left-open. Making use of Lemma \ref{lemmaspacequivalence} and the Ramsey property associated to $o^*$ it is easy to see that $X$ is a strong $S$-space. For more detail the reader can read the Chapter 2 of \cite{PartitionProblems}.\\\\
In the same way, we can define for each $\alpha\in\omega_1$, $y_\alpha$ as follows:
$$y_\alpha(\beta)=\begin{cases} \min(o^*(\alpha,\beta),1) &\textit{ if }\alpha>\beta\\
0 &\textit{ if } \alpha<\beta\\
1 &\textit{ if }\alpha=\beta
\end{cases}$$
By a similar argument, one can show that $\{y_\alpha\}_{\alpha\in \omega_1}$ is a strong $L$-space.
 It is worth pointing out that unlike the case of $S$-spaces, the $L$-spaces exist in $ZFC$ (see \cite{solutionlspace}).\\\\
Now, we will build another $S$-space using our family $\mathcal{B}_\mathcal{F}$. For every $\alpha\in\omega_1$, define $C(\alpha)=\{f_\xi\,|\,f_\xi\leq f_\alpha\}$. Now, let $\tau$ be the topology over $\mathcal{B}_\mathcal{F}$ obtained by refining the canonical Baire topology of $\omega^\omega$ restricted to $\mathcal{B}_\mathcal{F}$ by declaring the sets $C(\alpha)$ open. It is straight forward that each $\alpha\in\omega_1$ has as a local base the following family:
$$\{C(\alpha)\cap [s]\,|\,s\in\omega^{<\omega}\textit{ and }f_\alpha\in[s]\}.$$
Here, $[s]=\{f\in\omega^\omega\,|\,s\subseteq f\}$.\\\\
The following is based on the third author's proof that $\mathfrak{b}=\omega_1$ implies that there is an $S$-space.
\begin{proposition} $(\mathcal{B}_\mathcal{F}, \tau)$ is an $S$-space.
\begin{proof}
For each $\alpha\in\omega_1$ we have that $C(\alpha)$ is closed in the Baire topology, so it is clopen in $\tau.$ In this way, $\mathcal{B}_\mathcal{F}$ is $0$-dimensional, so it is regular. Also, $\mathcal{B}_\mathcal{F}$ is left open by definition. The only thing left to prove is that $\mathcal{B}_\mathcal{F}$ does not contain any uncountable discrete set.\\\\
Let $S\in[\omega_1]^{\omega_1}$ and assume $\{f_\alpha\}_{\alpha\in S}$ is discrete. For each $\alpha\in S$, we can find $s_\alpha\in \omega^{<\omega}$ such that if $U_\alpha=C(\alpha)\cap [s_\alpha]$ then $U_\alpha\cap Y=\{\alpha\}$ for each $\alpha\in\omega_1.$ Let $W\in [S]^{\omega_1}$ and $s\in\omega^{<\omega}$ such that $s_\alpha=s$ for all $\alpha\in W.$ By Theorem \ref{fulloscilation}, there are $\alpha<\beta\in W$ for which $osc(\alpha,\beta)=0$. Since $\alpha<\beta$, this means $f_\alpha<f_\beta,$ so $f_\alpha\in C(\beta)$. Thus $f_\alpha\in U_\beta$, which is a contradiction.
\end{proof}
\end{proposition}
\begin{corollary}[$CA_2$]There is a  first countable $S$-space. 
\end{corollary}
Now, we present the construction of a distinct  $S$-space. For this construction, we adapt the ideas from Chapter 2 of \cite{PartitionProblems}.
\begin{definition}let $l\in\omega\backslash\{0\}$ and define $H_l:\omega_1\longrightarrow [\omega_1]^{<\omega}$ as:$$H_l(\beta)=\{\alpha\in (\beta)^-_{l+1}\,|\,|(\alpha)_l|=|(\beta)_l|\}.$$
Finally, let $H(\beta)=\bigcup\limits_{l\in\omega\backslash\{0\}}H_l(\beta)$. 
\end{definition}
The following lemma is direct consequence of the definition.
\begin{lemma}Let $\beta\in\omega_1$ and $l\in\omega\backslash\{0\}.$ Then $H(\beta)=\{\alpha\in(\beta)^-_{l+1}\,|\,\Delta(\alpha,\beta)=l+1\}$
\end{lemma}
As we will see, this function shares some similarities to the one from Chapter 2 of \cite{PartitionProblems} (although ours is slightly simpler). Since $\mathcal{F}$ is $2$-capturing, we also have the following.
\begin{lemma}\label{Hbetalemma}Let $S\subseteq [\omega_1]^{<\omega}$ be an uncountable family of pairwise disjoint sets. Then
there exist $a,b\in S$ of the same cardinality $n$ such that $a<b$ and such that
$a(i)\in H(b(i))$
for all $i<n.$\end{lemma}
\begin{definition}\label{Cbeta}For each $\beta\in\omega_1$, we recursively define $C(\beta)\subseteq \beta+1$. Assume $C(\gamma)$ has been defined for each $\gamma<\beta.$ $C(\beta)$ will be the set containing $\beta$ and all $\alpha<\beta$ for which there is $\gamma\in H(\beta)$ such that:\begin{enumerate}[label=$(\alph*)$]
\item $\alpha\in C(\gamma)$, 
\item If $\xi\in H(\beta)\cup \{\beta\}$ and $\xi\not=\gamma$, then $\Delta(\alpha,\gamma)>\Delta(\alpha,\xi).$
\end{enumerate}
Finally, let $C_k(\beta)=\{\alpha\in C(\beta)\,|\,\Delta(\alpha,\beta)\geq k\}$ for each $k\in\omega.$
\end{definition}
Note that $\beta\in C_k(\beta)$ for all $k\in\omega$, and if $k<l$ then $C_l(\beta)\subseteq C_k(\beta).$ \begin{lemma}Let $\alpha,\beta\in\omega_1$ and $k\in\omega.$ If $\alpha\in C_s(\beta)$, there is $l\in\omega$ such that $C_s(\alpha)\subseteq C_k(\beta).$

\end{lemma}
For each $\beta\in\omega_1$ and $k\in\omega$ let $\hat{C}_k(\beta)=\{f_\alpha\,|\,\alpha\in C_k(\beta)\}$ and $\hat{C}(\beta)=\{f_\alpha\,|\,\alpha\in C(\beta)\}.$ By the previous Lemma, if $f_\alpha\in \hat{C}_{k_0}(\beta_0)\cap \hat{C}_{k_1}(\beta_1)$ there is $s_0,s_1\in\omega$ such that $\hat{C}_{s_i}(\alpha)\subseteq \hat{C}_{k_i}(\beta_i)$ for each $i\in 2$. Observe that if $s>s_0,s_1$ then $\hat{C}_s(\alpha)\subseteq \hat{C}_{k_0}(\beta_0)\cap \hat{C}_{k_1}(\beta_1)$. This means $\{\hat{C}_k(\beta)\,|\,k\in\omega\textit{ and }\beta\in\omega_1\}$ forms a base for a topology $\tau_C$ over $\mathcal{B}_\mathcal{F}.$ \\
The following lemma follows directly from the fact that $\hat{C}_k(f_\beta)\subseteq [f_\beta|_k]$ for each $\beta\in\omega_1$ and $k\in\omega.$
\begin{lemma}\label{bairerefinement}Let $s\in\omega^{<\omega}$. Then $[s]\cap \mathcal{B}_\mathcal{F}$ is open in $(\mathcal{B}_\mathcal{F},\tau_C).$ In particular, $(\mathcal{B}_\mathcal{F},\tau_C)$ is Hausdorff.
\end{lemma}
\begin{lemma}Let $\beta\in\omega_1$. $\hat{C}(\beta)$ is compact.
 \end{lemma}

 \begin{proposition}$(\mathcal{B}_\mathcal{F},\tau_C)$ is a locally compact strong $S$-space.
 \begin{proof}$\mathcal{B}_\mathcal{F}$ is locally compact by the previous result. In particular, it is regular and clearly it is left-open. Fix $n\in \omega$. It remains to prove $\mathcal{B}_{\mathcal{F}}^{n+1}$ has no uncountable discrete subspaces. For this, let $S\subseteq \omega_1^n$ be uncountable, and assume towards a contradiction that $\{(f_{x(0)},\dots,f_{x(n)})\}_{x\in S}$ is discrete. Without loss of generality we can suppose $x(i)<x(j)$ whenever $i<j\leq n$ and $\{x(i)\}_{i\leq n}\cap \{y(i)\}_{i\leq n}=\emptyset$ for all $x,y\in S$ with $x\not=y$. Furthermore, be a refining argument we can also suppose there is $k\in\omega$ such that for all $x,y\in S$, the following happens: \\\begin{enumerate}
     \item $\big(\prod\limits_{i\leq n}C_k(x(i))\big)\cap S=\{x\}$,
     \item $f_{x(i)}|_k=f_{y(i)}|_k$ for every $i\leq n$.\\
 \end{enumerate}
     Due to Lemma \ref{Hbetalemma}, we know there are distinct $x,y\in S$ such that $x(i)\in H(y(i))$ for every $i\leq n$. For any such $i$, we know  $\Delta(x(i),\xi)=\omega$ if and only if $\xi=x(i)$. Since $x(i)$ clearly belongs to $C(x(i))$, it follows from the definition that $x(i)\in C(y(i))$. But $f_{x(i)}|_k=f_y(i)|_k$, so $x(i)$ in fact is an element of $ C_k(x(i))$. In this way, we conclude that $x\in \big(\prod\limits_{i\leq n}C_k(y(i))\big)\cap S$, which is a contradiction. 
     \end{proof}
  \end{proposition}
  \begin{proposition}The Alexandroff compactification of $(\mathcal{B}_\mathcal{F},\tau_C)$ is a compact strong $S$-space.
   \end{proposition}
   \begin{corollary}[$CA_2$]There is a scattered compact strong $S$-space $K$ whose function space $C(K)$
is hereditarily weakly Lindel\"of and whose space $P(K)$ of all probability
measures is also a strong $S$-space.
\end{corollary}
   \subsection{A second look into Baumgartner Axiom} As we saw on Theorem \ref{entangledscheme}, FCA implies the existence of entangled sets, which means FCA also implies the failure of $BA(\omega_1)$. Although we do not know if $CA_2$ implies the existence of an entangled set, we will prove in this subsection that it does imply the negation of $BA(\omega_1).$ Our proof is based on the third author's proof that $\mathfrak{b}=\omega_1$ implies the failure of $BA(\omega)$ (see \cite{homeomorphismsofmanifolds}).  Remember that at the beginning of this section we fixed $2$-capturing construction scheme, namely $\mathcal{F}$.     \begin{lemma}\label{countableZ}The set $Z=\{\alpha\in\omega_1\,|\,\forall^{\infty}k\in\omega\,\exists F\in\mathcal{F}_{k+1}\,(\,\alpha\in F_{n_{k+1}-1}\,)\}$
   is at most countable.
   \begin{proof}Assume, by contradiction, that $Z$ is uncountable. For every $\alpha\in Z$, let $l_\alpha\in\omega$ be such that $\alpha\in F_{n_{i+1}-1}$ for each $i>l_\alpha$ and $F\in\mathcal{F}_{i+1}$ with $\alpha\in F$. We can find $W\in[Z]^{\omega_1}$ and $l\in\omega$ for which  $l_\alpha=l$ for all $\alpha\in W$.\\
   Since $\mathcal{F}$ is $2$-capturing, there are $\alpha,\beta\in W$, $l<i\in\omega$ and $F\in \mathcal{F}_{i+1}$ such that $F$ captures $\alpha,\beta.$ In particular, $\alpha\in F_0\backslash R(F)$ which is a contradiction.
   \end{proof}
   \end{lemma}
   Now, let us fix $M$ a countable elementary submodel of some large enough $H(\theta)$ with $\mathcal{F},\mathcal{B}_\mathcal{F}\in M.$ Define $\mathcal{A}=\mathcal{B}_\mathcal{F}\backslash M.$ The following lemma is a direct consequence of elementarity.
   \begin{lemma}\label{lemmaelementarysubmodel}Let $s\in\omega^{<\omega}.$ If $\mathcal{A}\cap[s]\not=\emptyset,$ then $\mathcal{A}\cap[s]$ is uncountable.
   \end{lemma}
   Remember that, by Definition \ref{definitionlex} and Remark \ref{remarklex}, we can think of $(\mathcal{A},<_{lex})$ as a subset of $\mathbb{R}.$
   \begin{lemma} $(\mathcal{A},<_{lex})$ is $\omega_1$-dense.
   \begin{proof}Let $f_\alpha,f_\beta\in \mathcal{A}$ with $f_\alpha<_{lex} f_\beta$, and let $l=\Delta(\alpha,\beta)$. Since the set $Z=\{\alpha\in\omega_1\,|\,\forall^{\infty}k\in\omega\,\exists F\in\mathcal{F}_{k+1}\,(\,\alpha\in F_{n_{k+1}-1}\,)\}$  is countable by Lemma \ref{countableZ} and is definable from $\mathcal{F}$ we have that $Z\subseteq M$. Hence, $\{f_\xi\,|\, \xi\in Z\}$ is disjoint from $\mathcal{A}$. In this way, there is $k>l$ and $F\in\mathcal{F}_{k+1}$ such that $\alpha\in F\backslash F_{n_{k+1}-1}.$ Let $\gamma\in F_{n_{k+1}-1}$ such that $|(\{\gamma\})_k|=|(\{\alpha\})_k|$. By Lemma \ref{lemmaelementarysubmodel}, it follows that $S=\mathcal{A}\cap[f_{\gamma}|_{k+2}]$ is uncountable. In order to finish, we will show that $S$ is contained in the open interval given by $f_\alpha$ and $f_\beta$. For this, let $f_\xi\in S$. By definition, $f_\xi|_{k+1}=f_\alpha|_{k+1}$ so $f_\xi<_{lex}f_\beta$. Also  $f_\alpha(k+1)<f_\gamma(k+1)$, which means that $f_\alpha<_{lex}f_\xi$. 
   \end{proof}
   \end{lemma}
   Let $\alpha\in\omega_1.$ Define $h_\alpha:\omega\longrightarrow \omega$ as:$$h_\alpha(i)=m_i-f_{\alpha}(i.)$$
  Also, let $-\mathcal{A}=\{h_\alpha(i)\,|\,\alpha\in \omega_1\backslash M\}$. It is easy to see that $f_\alpha<_{lex} f_\beta$ if and only if $h_\beta<_{lex}h_\alpha.$ Hence, we have the following corollary.
  \begin{corollary}$(-\mathcal{A},<_{lex})$ is $\omega_1$-dense.
  
 \end{corollary}
 \begin{proposition} There is no increasing function from $\mathcal{A}$ to $-\mathcal{A}.$
 \begin{proof}By contradiction, assume that there is $\Psi:\mathcal{A}\longrightarrow-\mathcal{A}$ increasing. Define $\psi:\omega_1\backslash M\longrightarrow \omega_1\backslash M$ in such way that $\Psi(f_\alpha)=h_{\psi(\alpha)}$.\\We claim that $\psi$ has at most one fixed point. Suppose this is not true and let $\alpha,\beta\in \omega_1$  fixed points  of $\psi$ such that $f_\alpha<_{lex}f_\beta$. since $\Psi$ is increasing we have that \begin{align*}h_\alpha= h_{\psi(\alpha)}=\Psi(f_\alpha)<_{lex} \Psi(f_\beta)=h_{\psi(\beta)}=h_\beta.
 \end{align*}
 But this is a contradiction since, in fact, $h_\beta<_{lex}h_\alpha.$\\\\
 Now, let $X=\{\alpha\in\omega_1\,|\,\alpha<\psi(\alpha)\}$ and $Y=\{\alpha\in\omega_1\,|\,\alpha>\psi(\alpha)\}$. By the previous claim, one of this sets is  uncountable.\\\\
 Assume $X$ is uncountable. For each $\alpha\in X$, let $b_\alpha=\{\alpha,\phi(\alpha)\}$. Since $\mathcal{F}$ is $2$-capturing, we can find $\alpha,\beta\in X$ and $F\in\mathcal{F}$ which captures $\{b_\alpha,b_\beta\}.$ Observe that $f_\alpha<_{lex} f_\beta$ and $f_{\psi(\alpha)}<_{lex}f_{\psi(\beta)}$, or equivalently, $\Psi(f_\alpha)=h_{\psi(\alpha)}>_{lex}h_{\psi(\beta)}=\Psi(f_\beta)$. This is a contradiction, so we are done.
 \end{proof}
 \end{proposition}
 \begin{corollary}[$CA_2$] There are two $\omega_1$-dense sets of reals which are not isomorphic.
 \end{corollary}
\section{Capturing schemes from $\Diamond$}\label{diamondsection}
In \cite{schemeseparablestructures} there is a proof that $\Diamond$-principle implies $FCA$. Unfortunately, the proof presented there is incomplete. Fortunately, the Theorem is true. We present a correct proof in this section.
As we have done in the previous sections, we also fix a type $\{(m_k,n_{k+1},r_{k+1})\}_{k\in\omega}$. \\\\
We start by proving various results about the interaction between construction schemes over different sets of ordinals.
\begin{lemma}\label{lemmafiniteinterval}Let $X$ be a set of ordinals, $\mathcal{F}$ be a construction scheme over $X$ and $k<l\in\omega$. For every $F\in\mathcal{F}_l$, $G\in\mathcal{F}_k$ with $G\subseteq F$ and $\alpha\in G$ there is $E\in\mathcal{F}_k$ such that:\begin{itemize}
    \item $E\subseteq F,$
    \item $E\cap( \alpha+1)=G\cap (\alpha+1) (in particular, \alpha\in E),$
    \item $E\backslash \alpha$ is an interval in $F.$
\end{itemize}
\begin{proof}The proof of this lemma is by induction over $l.$\\\\
{(Base step)} If $l=k+1$, let $E=F_0$ if $\alpha\in R(F)$ and $E=G$ otherwise. By using part $(b)$ of Lemma \ref{interseccioninicial}, we know there is $i<n_l$ for which $E=F_i$. Trivially $E\subseteq F$ and $E\cap(\alpha+1)=G\cap(\alpha+1)$ independently of the case. If $\alpha\in R(F)$ we use that $F_0$ is an interval in $F$ to conclude that so is $E\backslash \alpha$. On the other hand, if $\alpha\not\in R(F)$ we use that $F_i\backslash R(F)$ is an interval in $F$ to finish the proof of this case.\\\\
{(Induction step)} Suppose we have proved the lemma for $k<j<l$. Using Lemma \ref{interseccioninicial} we get $j<n_l$ for which $G\subseteq F_j$. By induction hypothesis we know there is $E'\in\mathcal{F}_k$ such that $E'\subseteq F_j$, 
$E'\cap(\alpha+1)=G'\cap(\alpha+1)$ and $E'\backslash \alpha$ is an interval in $F_j$. If $\alpha\in F_j\backslash R(F)$ let $E=E'$. Since $F_j\backslash R(F)$ is an interval in $F$
so is $E\backslash \alpha$, so we are done. If $\alpha\in R(F)$, let $h:F_j\longrightarrow 
F_0$ be the unique increasing bijection. By Corollary \ref{corollaryisomorphismfiniteschemes}, 
 $E=h[E']\in \mathcal{F}$. As $h$ is increasing and $h|_R(F)$ is the 
identity, $E\cap(\alpha+1)=E'\cap(\alpha+1)$ and $E\backslash \alpha$ is an 
interval in $F_0$. Since $F_0$ is an interval in $F$, we are done.
\end{proof}
\end{lemma}
\begin{definition}Let $\mathcal{F}$ be a construction scheme over some set of ordinals $Y$ and let $X\subseteq Y.$ We define $\mathcal{F}|_X$ as $\{F\in \mathcal{F}\,|\,F\subseteq X\}.$
\end{definition}
\begin{proposition}\label{lemmainterval}Fix $\mathcal{F}$ a construction scheme over some limit ordinal $\gamma$. Let $G\in \mathcal{F}$ and $\alpha\in G$. Then $\big(G\cap\alpha+1\big)\cup \big((\alpha+n)\backslash \alpha\big)\in \mathcal{F}$ where $n=|G\backslash \alpha|$. 
\begin{proof}Let $F\in \mathcal{F}$ be such that $ G\cup \big((\alpha+n)\backslash \alpha\big)\subseteq F$. Making use of Lemma \ref{lemmafiniteinterval}, we conclude there is $E\in \mathcal{F}$ satisfying the following properties:\begin{enumerate}
\item $\rho^E=\rho^G,$
\item $E\subseteq F,$
\item $E\cap(\alpha+1)=G\cap(\alpha+1),$
\item $E\backslash \alpha$ is an interval in $F$.
\end{enumerate}
Observe that both $E\backslash \alpha$ and $(\alpha+n)\backslash \alpha$ are two intervals in $F$ of cardinality $n$. Furthermore, $\alpha$ is the first element of both. Hence, $E\backslash \alpha=(\alpha+n)\backslash \alpha$. By (3), we conclude that $E=\big(G\cap\alpha+1\big)\cup \big((\alpha+n)\backslash \alpha).$
\end{proof}
\end{proposition}
\begin{definition}Let $p$ be a nonempty finite set of ordinals and  $\delta$ be a limit ordinal. Whenever $p\cap\delta\not=\emptyset$, we define $red_\delta(p)$ as $\big(p\cap\delta\big)\cup \big((\alpha+n)\backslash \alpha\big)$ where $\alpha=\max(p\cap \delta)$ and $n=|p\backslash \alpha|$. On the other hand, if $p\cap\delta=\emptyset$ we define $red_\delta(p)$ as $|p|.$
\end{definition}
\begin{corollary}\label{schemerestriction}Fix $\mathcal{F}$ a construction scheme over some limit ordinal $\gamma$. If $\delta<\gamma$ is a limit ordinal, then $\mathcal{F}|_\delta$ is a construction scheme over $\delta$.
\begin{proof} It should be clear that $\mathcal{F}|_{\delta}$ satisfies conditions (2), (3) and (4) of Definition \ref{constructionscheme}. In order to show that condition (1) is satisfied, take an arbitrary $A\in [\delta]^{<\omega}$. Since $\mathcal{F}$ is a construction scheme over $\gamma$, there is $F\in \mathcal{F}$ containing $A$. By Proposition $2$, we know $red_\delta(F)\in \mathcal{F}|_{\delta}$. To finish, just note that $A\subseteq red_\delta(F)$.
\end{proof}
\end{corollary}
The following is easy.
\begin{proposition}\label{unionschemes}Let $\mathcal{A}$ be a family of sets of ordinals linearly ordered by inclusion and for each $X\in \mathcal{A}$ let $\mathcal{F}_X$ be a construction scheme over $X$ in such way that $\mathcal{F}_X\subseteq \mathcal{F}_Y$ whenever $X\subseteq Y$. The $\mathcal{F}=\bigcup\limits_{X\in \mathcal{A}}\mathcal{F}_X$ is a construction scheme over $Z=\bigcup \mathcal{A}$.
\end{proposition}
\begin{corollary}\label{corollaryomega}There is a unique construction scheme over $\omega$.
\begin{proof}By Proposition \ref{unionschemes}, $\mathcal{F}=\bigcup_{k\in\omega}\mathcal{F}(m_k)$ is a construction scheme over $\omega$. In order to see that it is unique, take another construction scheme over $\omega$, namely $\mathcal{G}$. Let $k\in\omega $ and $F\in \mathcal{G}_k$ be such that $0\in F$. We make use of Proposition \ref{lemmainterval} to conclude $m_k=F\cap 1\cup m_k\backslash 0\in \mathcal{G} $. This implies $\mathcal{F}(m_k)\subseteq \mathcal{G}$. Thus, $\mathcal{F}\subseteq \mathcal{G}.$ By Lemma \ref{equalconstructions}, we are done.
\end{proof}
\end{corollary}
From now on, we will call $\mathcal{F}(\omega)$ the unique construction scheme over $\omega.$
\begin{definition}\label{definitionforcing} Let $\mathcal{F}$ be a construction scheme over some ordinal $\gamma$. We define $\mathbb{P}(\mathcal{F})$ to be the set of all $p\in [\gamma+\omega]^{<\omega}$ such that:\begin{enumerate}[label=$(\alph*)$]
\item $\exists k\in\omega\,(\,|p|=m_k\,),$
\item $\exists F\in \mathcal{F}\,(|F|=|p|\,\textit{ and }p\cap \gamma\sqsubseteq F\,),$
\item $p\backslash \gamma\sqsubseteq (\gamma+\omega)\backslash \gamma.$
\end{enumerate}
We also  let $p<q$ if and only if $q\in \mathcal{F}(p).$ Finally, we denote $\{p\in \mathbb{P}(\mathcal{F})\,|\,|p|=m_k\}$ as $\mathbb{P}_k(\mathcal{F})$ for each $k\in\omega.$
\end{definition}
\begin{remark}\label{remarkreductionforcing}Let $p\in [\gamma+\omega]^{<\omega}$ which satisfies conditions $(a)$ and $(b)$ of Definition \ref{definitionforcing}. Then $p\in \mathbb{P}(\mathcal{F})$ if and only if $red_\delta(p)\in \mathcal{F}$.
\end{remark}
The following Definition appeared for the first time in \cite{schemeseparablestructures}. As the third author showed, this hypothesis is sufficient to build construction schemes in a recursive manner.
\begin{definition}[$IH_1$]Let $\mathcal{F}$ be a construction scheme on some set of ordinals $X$. We say that $\mathcal{F}$ satisfies $IH_1$ if for every $A\in [X]^{<\omega}$ and $\alpha\in X$, there is $F\in \mathcal{F}$ such that:
\begin{enumerate}
    \item $A\subseteq F_0,$
    \item $R(F)=F\cap \alpha$.
\end{enumerate}
\end{definition}
As for now, we have barely used the condition $(c)$ of Definition \ref{definitiontype}. The next proposition illustrates its utility. 
\begin{proposition}\label{ih1omega}$\mathcal{F}(\omega)$ satisfies $IH_1.$
\begin{proof}Let $A\in [\omega]^{<\omega}$ and $\alpha\in \omega.$ By condition $(c)$ of Definition \ref{definitiontype}, there are infinitely many $k\in\omega$ for which $r_k=\alpha$. Fix such $k$ with the property that $m_{k-1}>\max(A)$. Then $F=m_k\in \mathcal{F}(\omega)$, $A\subseteq m_{k-1}=F_0$ and $F\cap \alpha=r_k=R(F)$. Hence, we are done.
\end{proof}
\end{proposition}
\begin{proposition}\label{ih1unions}Let $\mathcal{A}$ be a family of sets of ordinals linearly ordered by inclusion and for each $X\in \mathcal{A}$ let $\mathcal{F}_X$ be a construction scheme over $X$ in such way that $\mathcal{F}_X\subseteq \mathcal{F}_Y$ whenever $X\subseteq Y$. Then $\mathcal{F}=\bigcup\limits_{X\in \mathcal{A}}\mathcal{F}_X$ satisfies $IH_1$ if and only if the same is true for all $\mathcal{F}_X.$
\end{proposition}
The following theorem was first proved in \cite{schemeseparablestructures}.
\begin{definition}Let $\mathcal{F}$ be a construction scheme over some limit ordinal $\gamma$ and let $\mathcal{G}$ be a filter over $\mathbb{P}(\mathcal{F})$. We define $\mathcal{F}^\mathcal{G}$ as $\bigcup_{p\in \mathcal{G}}\mathcal{F}(p)$. Finally, $\mathcal{F}^{Gen}$ denotes the name for $\mathcal{F}^G$ where $G$ is a generic filter.
\end{definition}

\begin{theorem}Let $\mathcal{F}$ be a construction scheme over some limit ordinal $\gamma$ which satisfies $IH_1$. Then there is a countable family $\mathcal{D}$ of dense sets in $\mathbb{P}(\mathcal{F})$ such that whenever $\mathcal{G}$ is a filter intersecting each member of $\mathcal{D}$,  then $\mathcal{F}^\mathcal{G}$ is a construction scheme over $\gamma+\omega$ satisfying $IH_1$ and containing $\mathcal{F}$. 
\end{theorem}

Observe that Theorem \ref{theoremexistenceconstructionschemes} is a direct consequence of Corollary \ref{corollaryomega}, Proposition \ref{ih1omega} and  Proposition \ref{ih1unions}. Furthermore, Corollary \ref{schemerestriction} and Lemma \ref{equalconstructions} tell us that every construction scheme over $\omega_1$ is completely determined by an increasing sequence of construction schemes over countable limit ordinals.
\begin{Problem}Is there a construction scheme over $\omega_1$ which do not satisfy $IH_1$?
\end{Problem}
\subsection{The property $IH_2$}Otherwise stated, $\mathcal{F}$ will denote a construction scheme over some limit ordinal $\gamma$ which satisfies $IH_1.$

\begin{definition}Let $k\in \omega$ and $\mathbb{I}=\{\mathbb{I}(i)\}_{i<j}$ be sequence of nonempty finite sets of ordinals. We say that $\mathbb{I}$ is a block interval sequence if:
\begin{enumerate}[label=$(\alph*)$]
\item $\mathbb{I}(j)$ is an interval for each $i<j,$
    \item $\mathbb{I}(j)<\mathbb{I}(j')$ for every $j<j'<i.$
\end{enumerate} Given $k\in \omega\backslash 2$ and $t<m_k$ we let $Bl(t,k)$  be the set of all the block interval sequences contained in $\mathscr{P}(\,[t,m_k))$. Given $\beta\in\gamma$ we make an abuse of notation by naming $Bl(|(\beta)^-_k|,k)$ as $Bl(\beta,k)$.
\end{definition}
Note that if $\mathbb{I}=\{\mathbb{I}(i)\}_{i<j}$ is a block interval sequence, then each $\mathbb{I}(i)$ is nonempty (for $i<j$). However, $\mathbb{I}$ itself might be empty. This is not a weird pathological case, it will be important for us in the future.
\begin{definition}Let $k\in\omega\backslash 2$, and $t<m_k$. We say that $(\mathbb{I},z_\mathbb{I})\in Bl(t,k)\times m_k$ is $(t,k)$-good:
\begin{itemize}
\item $z_\mathbb{I}\in [t,m_k)$, 
\item Either $\mathbb{I}=\emptyset$ or $\max(\bigcup\mathbb{I})\leq z_\mathbb{I}.$
\end{itemize}
For a nonempty $S\subseteq BL(t,k)$, we say that a sequence $T=\{(\mathbb{I},z_\mathbb{I})\}_{\mathbb{I}\in S}$ is $(t,k)$-good if $(\mathbb{I},z_\mathbb{I})$ is $(t,k)$-good for every $\mathbb{I}\in S.$ The set of all $(t,k)$-good sequences is denoted as $Good(t,k)$. Given $\beta \in \gamma$, we make as slight abuse of notation and denote $Good(|(\beta)_k^-|,k)$ as $Good(\beta,k).$
\end{definition}
In particular, if $\mathbb{I}=\emptyset$ then  $(\emptyset,z_\mathbb{I})$ is $(t,k)$-good if and only if $t\leq z_\mathbb{I}<m_k$.\\\\
In this moment, it may be convenient to recall some notation. If $a$ is a set of ordinals (or in general, a linear order) of size $n$, we decided to denote as $a(i)$ the $i$-th element of $a$ with respect to the increasing order. Following this notation, $a[S]$ denotes the set $\{a(i)\,|\, i\in S\}$ for each $S\subseteq n$. For the rest of this section we will apply this notation, in particular, to sets of the form $(\xi)_k$ and elements of the construction scheme.
\begin{definition}\label{definitionprojection}Let $\xi \in \gamma$, $k<l\in \omega\backslash 2$, and $(\mathbb{I},z_\mathbb{I})$ be $(0,k)$-good. Whenever $z_\mathbb{I}\leq |(\xi)^-_k|$, we define the $(\xi,l)$-projection of $\mathbb{I}$  as $\pi^{k,l}_\xi(\mathbb{I})=\{\pi^{k,l}_\xi(\mathbb{I})(j)\}_{j<|\mathbb{I}|}$ where $$\pi^{k,l}_\xi(\mathbb{I})(j)=\{|(\alpha)^-_l|\,|\,\alpha\in (\xi)_k[\mathbb{I}(j)]\}.$$
\end{definition}
A crucial fact is that, in general, the set $\pi^{k,l}_\xi(\mathbb{I})(j)$ is just a finite set, not an interval. In order for $\pi^{k,l}_\xi(\mathbb{I})$ to be an element of $Bl(0,l)$ it is a necessary an sufficient condition that each $\pi^{k,l}_\xi(\mathbb{I})(j)$ is an interval in $m_l$. This is equivalent to saying that $(\xi)_k[\mathbb{I}(j)]$ is an interval in $(\xi)_l.$\\\\
Recall that a $\Diamond$-sequence is a sequence $\{D_\alpha\}_{\alpha\in LIM}$ of 
countable subsets of $\omega_1$ such that:\begin{enumerate}[label=$(\arabic*)$]
    \item $\forall \alpha\in LIM\,(D_\alpha\subseteq \alpha),$
    \item $\forall D\subseteq \omega_1\,(\,\{\alpha\in LIM\,|\,D\cap \alpha=D_\alpha\}\textit{ is stationary }).$
\end{enumerate}
$\Diamond$-principle (or simply $\Diamond$) is the statement that asserts the existence of a $\Diamond$-sequence.   
For the rest of this section, we fix a $\Diamond$-sequence $\{D_\alpha\}_{\alpha\in LIM}$. 
\begin{definition}\label{definitioncheckmark}Let $\beta,\xi \in \gamma$, $k\leq l\in \omega\backslash 2$ and $T=\{(\mathbb{I},z_\mathbb{I})\}_{\mathbb{I}\in S}\in Good(\beta,k)$. We say that $(\beta,\xi,k,l)$ approves $T$ $($and write it as $(\beta,\xi,k,l) \checkmark T)$ there is $\mathbb{I}\in S$ for which the following happens:
\begin{enumerate}[label=$(\arabic*)$]
\item[$(0)$] $|(\beta)_l|\leq |(\xi)_l|$,
\item $z_\mathbb{I}=|(\xi)^-_k|$,
\item $\pi^{k,l}_\xi(\mathbb{I})\in Bl(0,l)$,
\item $(\xi)_l(|(\beta)^-_l|)\in (\xi)_k.$
\end{enumerate}
\end{definition}
Frequently, we will apply the previous definition in cases where $\xi<\beta.$
\begin{remark}\label{remarkcheckmark} If $(\beta,\xi,k,l)\checkmark T$ and $\beta',\xi'\in \gamma$ are such that $|(\beta)_l|=|(\beta')_l|$ and $|(\xi)_l|=|(\xi')_l|$, then $(\beta',\xi',k,l)\checkmark T.$
\end{remark}
\begin{definition}\label{acceptance}Let $\delta<\beta<\gamma$, $k<l\in \omega\backslash 2$ and $T\in Good(\beta, k)$. Given $C\in [D_\delta]^{<\omega}\backslash\{\emptyset\}$ and $G\in \mathcal{F}_{l}$, we say that $(C, G)$ is accepted by $(k,l,\beta,\delta,T)$ if :\begin{itemize}
\item $G$ captures $C$,
\item $(\beta,C(0),k,l-1)\checkmark T,$
\item $(\beta)_{l-1}\cap \delta=R(G).$
\end{itemize}
Finally, we define $j(k,l,\beta,\delta,T)$ as $$\max\{j\leq n_l\,|\,\exists C\in [D_\delta]^{j}\,\exists G\in (\mathcal{F}|_\delta)_{l}\,\big((C, G)\textit{ is accepted by }(k,l,\beta,\delta,T)\big)\}\cup\{0\}$$
\end{definition}
\begin{definition}[$IH_2$]\label{definitionih2} We say that $\mathcal{F}$ satisfies $IH_2$ if for each limit $\delta< \gamma$ one of the following mutually excluding conditions holds:
\begin{enumerate}[label=$(\arabic*)$]
\item There are infinitely many $l\in\omega$ for which there $C\in [D_\delta]^{<\omega}$ which is fully captured by some element of $\mathcal{F}_l.$\\
\item For every $\beta \in \gamma\backslash \delta$, $k\in \omega\backslash 2$ and $T\in Good(\beta,k)$ there are infinitely many $l \in \omega\backslash k$ for which:
\begin{enumerate}[label=$(\alph*)$]
\item $|(\beta)_{l-1}\cap \delta|=r_{l},$
\item There is no $C\in [D_\delta]^{<\omega}$ which is fully captured by any $F\in \mathcal{F}_l.$
\end{enumerate}
Furthermore, either $j(k,l,\beta,\delta,T)=0$ and there is $F\in \mathcal{F}_l$ such that $\beta\in F_0\backslash R(F)$ or $j(k,l,\beta,\delta,T)>0$ and there are $C\in [D_\delta]^{j(k,l\beta,\delta,T)}$ and $F\in \mathcal{F}_{l}$ for which $(C,F)$ is accepted by $(k,l,\beta,\delta,T)$ and  $\beta\in F_{j(k,l,\beta,\delta,T)}$.
\end{enumerate}

\end{definition}

\begin{proposition}\label{ih2unions}Let $\mathcal{A}$ be a set of countable ordinals. For each $\gamma\in \mathcal{A}$ let $\mathcal{F}_\gamma$ be a construction scheme over $\gamma$ in such way that $\mathcal{F}_\delta\subseteq \mathcal{F}_{\gamma}$ whenever $\delta<\gamma$. Then $\mathcal{F}_\mathcal{A}=\bigcup\limits_{\gamma\in \mathcal{A}}\mathcal{F}_\gamma$ satisfies $IH_2$ if and only if the same is true for all $\mathcal{F}_\gamma.$
\end{proposition}
Since $\omega$ is the first limit ordinal, we trivially conclude the following.
\begin{remark}\label{ih2omega}$\mathcal{F}(\omega)$ satisfies $IH_2.$
\end{remark}
The following Lemma already appeared in \cite{lopezschemethesis} and \cite{schemeseparablestructures}. We give a proof for convenience of the reader.
\begin{lemma}\label{lemmacapturingsubset}Let $\mathcal{F}$ be a construction scheme over $\omega_1$ such that for every uncountable $S\subseteq \omega_1$, there are infinitely many $l\in\omega$ for which there is $C\in [S]^{<\omega}$ which is fully captured by some element of $\mathcal{F}_l$. Then $\mathcal{F}$ is fully capturing.
\begin{proof}Let $S\subseteq [\omega_1]^{<\omega}$ be uncountable. By refining $S$, we can suppose without loss of generality that there is $k\in \omega$ such that:
\begin{itemize}
    \item $\forall c\in S\,(\,\rho^c=k\,)$, 
    \item $\forall c,c'\in S\,(\,|(c)_k|=|(c')_k|\,),$
    \item $\forall c,c'\in S\, (\,\phi[c]=c'\textit{ where }\phi\textit{ is the increasing bijection from }(c)_k\textit{ to }(c')_k\,).$    
\end{itemize}
The set $S'=\{\max(c)\,|\,c\in S\}$ is uncountable. Thus, by the hypotheses we conclude there are infinitely many $l\in \omega$ for which there is $C\in [S']^{<\omega}$ which is fully captured by some element of $\mathcal{F}_l$. If such $l$ is greater than $k$, let $c_i\in S$ such that $\max(c_i)=C(i)$ for each $i<n_l$. Now, consider $\mathcal{C}=\{c_i\}_{i<n_l}$. We claim that $F$ fully captures $\mathcal{C}.$ For this, first observe that since $\rho^{c(i)}=k\leq l-1$ and $C(i)=max(c_i)\in F_i\backslash R(F)$ for each $i$, then $c_i\subseteq F_i$ and $c_i\backslash R(F)\not=\emptyset$. Now, fix $i<n_l$ and  let $\Psi:F_0\longrightarrow F_i$ be the increasing bijection. Since $F$ fully captures $C$, then $\Psi(C(0))=C(i)$.  In this way $\Psi|_{(C(0))_k}$ is the increasing bijection from $(C(0))_k$ to $((C(i))_k$. To finish, just notice that $(C(0))_k=(c_0)_k$ and $(C(i))_k=(c_i)_k.$ So by the third point of our \say{refinement} assumptions, we are done.
\end{proof}

\end{lemma}
\begin{theorem}\label{maintheorem}Let $\mathcal{F}$ be a construction scheme over $\omega_1$ which satisfies $IH_2$. Then $\mathcal{F}$ is fully capturing.
\begin{proof}We will use Lemma \ref{lemmacapturingsubset}. Let $S$ be an uncountable subset of $\omega_1$. Since $\{D_\delta\}_{\delta\in LIM}$ is a $\Diamond$-sequence, we can choose a limit ordinal $\delta$ such that: \begin{enumerate}[label=$(\arabic*)$]
    \item $\delta$ is a limit point of $S$,
    \item $D_\delta= S\cap \delta,$
    \item $(\delta,<,D_\delta,\mathcal{F}|_\delta)$ is an elementary submodel of $(\omega_1,<,S,\mathcal{F}).$
\end{enumerate}
It is enough to prove that for each $k\in \omega$ there are $C\in [D_\delta]^{<\omega},  l>k$ and $F\in \mathcal{F}_l$ which fully captures $C.$ Suppose towards a contradiction that this is not the case, and take $k\in\omega\backslash 2$ for which the given assertion is false. Now fix $\beta\in S\backslash \delta$ and let $T= \{(\emptyset, z_\emptyset)\}$ where $z_\emptyset=|(\beta)_k^-|$. We know $\mathcal{F}$ satisfies $IH_2$ and the condition (1) of Definition \ref{definitionih2} does not hold. In this way,  we can take $l\in \omega\backslash k$ which testifies condition (2) of Definition  \ref{definitionih2} applied to $\beta, k$ and $T$. \\
For $j=j(k,l,\beta,\delta, T)$ we have two cases. If $j=0$ there is $F\in \mathcal{F}_l$ such $\beta\in F_0\backslash R(F)$. By elementarity there are $G\in (\mathcal{F}|_\delta)_l$ and $\xi_0\in D_\delta$ for which the following happens:
\begin{itemize}
\item $\xi_0\in G_0\backslash R(G), $
\item $|(\xi_0)_{l-1}|=|(\beta)_{l-1}|,$
\item $(\beta)_{l-1}\cap \delta=R(G).$
\end{itemize}
The first property implies that $G$ captures $\{\xi_0\}$. The second one implies that $|(\xi_0)_k^-|=|(\beta)^-_k|=z_\emptyset$ and $(\xi_0)_{l-1}(|(\beta)_{l-1}^-|)=\xi_0\in (\xi_0)_k.$ Since $\pi^{k,l}_\xi(\emptyset)=\emptyset$, we conclude that $(\{\xi_0\},G)$ is accepted by $ (k,l,\beta,\delta,T).$ Thus, $j\geq 1$ which is a contradiction.\\\\
Now, we deal with the case where $j> 0$. In this case there are $C\in [D_\delta]^j$ and $F\in \mathcal{F}_l$ for which $(C,F)$ is accepted by $(k,l,\beta,\delta,T)$ and  $\beta\in F_j$. By elementarity, there are $G\in (\mathcal{F}|_\delta)_l$ and $\xi_j\in D_\delta$ which satisfy the following properties:
\begin{itemize}
\item $\xi_j\in G_j$,
\item $|(\xi_j)_{l-1}|=|(\beta)_{l-1}|,$
\item $G_j$ captures $C$,
\item $(\beta)_{l-1}\cap \delta=R(G),$
\item $C\subseteq \xi_j.$
\end{itemize}
Since $(\beta,C(0),k,l-1)\checkmark T$, we know that $|(C(0))^-_k|=z_\emptyset=|(\beta)^-_k|$ by condition (1) of Definition \ref{definitioncheckmark}. Furthermore, if we define $\alpha$ as $(C(0))_{l-1}(|(\beta)^-_{l-1}|)$, condition (3) of Definition \ref{definitioncheckmark} assure us that  $\alpha\in (C(0))_k$. By definition of $\alpha$ it is also true that $|(\alpha)_{l-1}|=|(\beta)_{l-1}|$. 
In particular, $|(\alpha)_k|=|(\beta)_k|$. In this way, $$C(0)=(C(0))_k(|(C(0))^-_k|)=(C(0))_k(|(\alpha)^-_k|)=\alpha.$$  Thus, $|(C(0))_{l-1}^-|=|(\beta)_{l-1}^-|=|(\xi_j)^-_{l-1}|$. From this, we  conclude $G$ captures $C\cup\{\xi_j\}$ and consequently $(C\cup\{\xi_j\},G)$ is accepted by $(k,l,\beta,\delta,T)$. By definition of $j(k,l,\beta,\delta,T)$, we have that $j(k,l,\beta,\delta,T)>j$ which is a contradiction. This completes the proof.
\end{proof}
\end{theorem}
\subsection{Forcing $IH_2$} Now we have a clear path to follow in order to prove the existence of a fully capturing construction scheme. As in the previous subsection, we fix a construction scheme $\mathcal{F}$ over some limit ordinal $\gamma$. This time, we assume $\mathcal{F}$  satisfies $IH_1$ and $IH_2$. 
 \begin{definition}Let $F\in \mathcal{F}$ and $\alpha\in F$. We define $Cut(F,\alpha)$ as  $(F\cap \alpha )\cup [\gamma,\gamma +|F\backslash \alpha|).$
\end{definition}
Note that $Cut(F,\alpha)$ is always a condition in $\mathbb{P}(\mathcal{F})$ as testified by $F$. In fact, every $p\in \mathbb{P}(\mathcal{F})$ can be written as $Cut(F,\alpha)$ for each $F\in \mathcal{F}$ testifying that $p$ is a condition. The following lemma is the main tool we have to extend conditions in $\mathbb{P}(\mathcal{F}).$
\begin{lemma}\label{lemmacut}Let $k\in\omega$, $p\in \mathbb{P}_k(\mathcal{F})$ and $\alpha\in \gamma$ be such that $(\alpha)^-_k=p\cap\gamma$. Then $Cut(F,\alpha)\leq p$ for each $F\in \mathcal{F}$ for which $\alpha\in F$ and $\rho^F\geq k.$
\begin{proof}Let $F\in \mathcal{F}$ be as in the hypotheses. Consider $G\in \mathcal{F}_k$ for which $\alpha\in G$ and $G\subseteq F.$ Observe that $G\cap\alpha=(\alpha)_k^-=p\cap \gamma.$ This means $Cut(G,\alpha)=p$. Since $Cut(G,\alpha)\in \mathcal{F}(Cut(F,\alpha))$, we are done. 
\end{proof}
\end{lemma}
\begin{lemma}\label{lemmabetagamma}Let $p\in \mathbb{P}(\mathcal{F})$ and $\alpha \in \gamma$. Then there is $q\leq p$ such that $\alpha\in q$.
\begin{proof} Let $\xi=\min(red_\gamma(p)\backslash p)$. Since $\mathcal{F}$ satisfies $IH_1$, we can take $F\in \mathcal{F}$ such that $\{\alpha\}\cup red_\gamma(p)\subseteq F_0$ and $\xi\cap F=R(F).$ Let $\xi'=\min(F_1\backslash R(F))$, $q=Cut(F,\xi')$ and $l=\rho^{F}$. Then $l>\rho^{red_\gamma(p)}$ because $red_\gamma(p)\subsetneq F$. Furthermore, $(\xi)^-_{l-1}=R(F)=(\xi')^-_{l-1}$ which means that $(\xi')^-_k=(\xi)^-_k=p\cap \gamma$. In this way, $q\leq p$ due to Lemma \ref{lemmacut}. To finish, just notice that $\alpha\in F_0\subseteq q.$
    
\end{proof}
\end{lemma}
\begin{lemma}\label{lemmaextensionr}Let $p\in \mathbb{P}(\mathcal{F})$, $\beta \in (\gamma+\omega)\backslash \gamma$ and $k\in \omega$. Then there are $k'\in \omega\backslash k$ and $q\in \mathbb{P}_{k'}(\mathcal{F})$ such that:\begin{itemize}
    \item $q\leq p$,
    \item $\beta \in q,$
    \item $|q\cap \gamma|= r_{k'+1},$
    \end{itemize}  
\begin{proof}
Let  $l\in\omega$ be such that $p\in \mathbb{P}_l(\mathcal{F})$. Also, let $s\in \omega$ be such that $\beta=\gamma+s$ and let $\alpha=\min (red_\gamma(p)\backslash p)$. Using that $\mathcal{F}$ satisfies $IH_1$, we can take an $F\in \mathcal{F}$ for which:\begin{enumerate}
\item $red_\gamma(p)\cup (\alpha+s+m_k+1)\backslash \alpha\subseteq F_0$,
    \item $R(F)=F\cap \alpha.$
\end{enumerate}
Let $q=Cut(F_0,\alpha)$ and $k'=\rho^{F_0}$. As $(\alpha+s+m_k+1)\backslash\alpha\subseteq F_0$, then $k'\geq k$. Also $k'\geq l$ because $red_\gamma(p)\subseteq F_0$. In this way, by Lemma \ref{lemmacut} we conclude that $q\leq p$. Since $q\cap \gamma=F\cap \alpha=R(F)$, then $|q\cap \gamma|= r_{k'+1}$. To finish, just notice $\beta=\gamma+s\in(\gamma+s+m_k+1)\backslash \gamma\subseteq q$.
\end{proof}
\end{lemma}
\begin{definition}\label{transgood}Let $\gamma'$ be a limit ordinal and $\mathcal{F}'$ be a construction scheme over $\gamma'$. Let $k<k'\in \omega\backslash 2$, $\beta\in \gamma'$, $S\subseteq Bl(\beta,k)$ and $T=\{(\mathbb{I},z_\mathbb{I})\}_{\mathbb{I}\in S}\in Good(\beta,k)$. Given $\alpha\in (\beta)^-_{k'}$, we define $\mathcal{Q}$ as the set of all pairs $(G,\mathbb{I})\in \mathcal{F}_k(m_{k'})\times S$ such that:\begin{enumerate}[label=$(\alph*)$]
    \item $|(\beta)^-_{k'}|\in G$,
    \item $G[\mathbb{I}(i)]$ is an interval in $m_{k'}$ for each $i<|\mathbb{I}|.$ 
    \end{enumerate}
    Furthermore, Given $(G,\mathbb{I})\in \mathcal{Q}$ we define the following objects:
\begin{itemize}
\item $\mathbb{J}_{(G,\mathbb{I})}$ is the sequence defined as $\mathbb{J}_{(G,\mathbb{I})}(0)=(\,|(\alpha)^-_{k'}|,|(\beta)^-_{k'}|\,]$ and $\mathbb{J}_{(G,\mathbb{I})}(i)=G[\mathbb{I}(i-1)]$ for  each $i<|\mathbb{I}|.$
\item $z_{\mathbb{J}_{(G,\mathbb{I})}}=G(z_\mathbb{I})$.
\end{itemize}
Finally, we let $S'=\{\mathbb{J}_{(G,\mathbb{I})}\,|\,(G,\mathbb{I})\in \mathcal{Q}\}$ and $Trans(k,k',\alpha,\beta,T)=\{(\mathbb{J},z_{\mathbb{J}})\}_{\mathbb{J}\in S'}$.
\end{definition}
Given $(G,\mathbb{I})\in \mathcal{Q}$, the property $(b)$ of the previous definition implies that $\mathbb{J}_{(G,\mathbb{I})}(i)$ is an interval in $m_{k'}$ 
for each $i\leq \mathbb{I}$. Furthermore, since $|(\beta)^-_{k'}|\in G$ then $G(|(\beta)^-_k|)=|(\beta)^-_{k'}|.$ Let $i<\mathbb{I}$. As a direct consequence of 
the previous fact and using that $\mathbb{I}(i+1)\subseteq (\,|(\beta)^-_k|,z_\mathbb{I}\,]$, we conclude that $$\{|(\alpha)^-_{k'}|\}<\mathbb{J}_{(G,\mathbb{I})}(0)<G(\mathbb{I}(i+1))\leq \{G(z_\mathbb{I})\}=\{z_{\mathbb{J}_{(G,\mathbb{I})}}\}.$$
Where the second inequality holds because $\max(\mathbb{J}_{(G,\mathbb{I})}(0))=G(|(\beta)^-_k|)$. From this, it should be clear that $Trans(k,k',\alpha,\beta,T)$ belongs to $Good(\alpha,k')$.
\begin{theorem}\label{maintheoremlemma}Let $\gamma'$, $\mathcal{F}'$, $k,k',\alpha,\beta$, $\mathcal{Q}$, $S$, and $T$ be as in Definition \ref{transgood}. Suppose  $\xi\in \gamma'$ and $l\in \omega\backslash(k'+1)$ are  such that $|(\beta)_l|\leq |(\xi)_l|$, $(\xi)_l(|(\beta)^-_l|)\in (\xi)_{k'}$ and $(\beta)_{k'}[\,(|(\alpha)^-_{k'}|,|(\beta)^-_{k'}|]\,]$ is an interval in $(\beta)_l$. Then the following statements are equivalent:
\begin{enumerate}[label=$(\alph*)$]
\item $(\alpha,\xi,k',l)\checkmark Trans(k,k',\alpha,\beta,T)$,
\item $(\beta,\xi,k,l)\checkmark T$.
\end{enumerate}
\begin{proof}Let $F\in \mathcal{F}'_{k'}$ be such that $\xi\in F$ and let $\phi:F\longrightarrow m_{k'}$ be the increasing bijection (which is given by $\phi(\alpha)=|(\alpha)_{k'}^-|$ and its the inverse of the $F(\,\_\,)$ function).\\\\
First we prove that (a) implies (b). Let $(G,\mathbb{I})\in \mathcal{Q}$ for which $\mathbb{J}_{(G,\mathbb{I})}$ satisfies conditions (0), (1), (2) and (3) of 
Definition \ref{definitioncheckmark} applied to $\alpha,\xi,k', l$ and $Trans(k,k',\alpha,\beta,T).$ That is, $|(\alpha)_l|\leq |(\xi)_l|$ and:\begin{enumerate}[label=(\Roman*)]
\item $|(\xi)^-_{k'}|=z_{\mathbb{J}_{(G,\mathbb{I})}}=G(z_\mathbb{I})$,
\item$\pi^{k',l}_\xi(\mathbb{J}_{(G,\mathbb{I})})\in Bl(0,l)$. So in particular, $(\xi)_{k'}[\mathbb{J}_{(G,\mathbb{I})}(i+1)]=(\xi)_{k'}[G[\mathbb{I}(i)]]$ for each $i<|\mathbb{I}|$,
\item $(\xi)_l(|(\alpha)^-_l|)\in (\xi)_k$.
\end{enumerate}
We claim  that $\mathbb{I}$ testifies $(\beta,\xi,k,l)\checkmark T.$ The condition (0) of Definition \ref{definitioncheckmark} is satisfied by hypothesis, so we will only prove that the remaining conditions hold: 
\begin{enumerate}
\item By (I), it  follows that $\phi(\xi)=G(z_\mathbb{I})$. Hence, $\phi[(\xi)^-_k]=(G(z_\mathbb{I}))^-_k=G(z_\mathbb{I})\cap G$. Since the function $\phi$ is injective, then $|(\xi)^-_k|=|G(z_{\mathbb{I}})\cap G|=z_\mathbb{I}.$ 

\item  We want to show that $\pi^{k,l}_\xi\in Bl(0,l)$. As we mentioned just after the Definition \ref{definitionprojection}, it is enough prove that $(\xi)_k[\mathbb{I}(i)]$ is an interval in $(\xi)_l$ for each $i<\mathbb{I}$. Observe that since $G(z_\mathbb{I})=|(\xi)^-_{k'}|$ then $\xi\in F[G]$ and $(\xi)_k(j)=F[G](j)=F(G(j))=(\xi)_{k'}(G(j))$ for each $j<|(\xi)^-_k|$. Consequently \begin{align*}(\xi)_{k'}[ G(z_\mathbb{I})\cap G]&=(\xi)_{k'}[\{G(j)\,|\,j<z_{\mathbb{I}}\,\}]\\
&=\{(\xi)_{k'}(G(j))\,|\,j<|(\xi)^-_k|\,\}\\
&=\{(\xi)_{k}(j)\,|\,j<|(\xi)^-_k|\,\}\\
&=(\xi)^-_{k}.   
\end{align*} As $G[\mathbb{I}(i)]\subseteq G\cap G(z_\mathbb{I})$, then  $(\xi)_{k'}[G[\mathbb{I}(i)]]=(\xi)_{k}[\mathbb{I}(i)]$. This means that this last set is an interval in $(\xi)_l$. Since $i$ was arbitrary, we conclude  $\pi^{k,l}_\xi(\mathbb{I})\in Bl(0,l).$

\item Finally, we prove that $(\xi)_l(|(\beta)^-_l|)\in (\xi)_k$. Remember that by hypothesis we have 
that $(\xi)_l(|(\beta)^-_l|)\in(\xi)_{k'}$. In this way, $(\xi)_l(|(\beta)^-_l|)=(\xi)_{k'}(|(\beta)^-_{k'}|)$. But both $\phi(\beta)=|(\beta)^-_{k'}|$ and $\phi(\xi)=|(\xi)^-_{k'}|$ are members of $G$. $\phi(\xi)\in G$ by $(I)$ above and $\phi(\beta)\in G$ since $(G,\mathbb{I})\in \mathcal{Q}$. This means  $(\phi(\xi))_{k'}(|(\beta)^-_{k'}|)=|(\beta)^-_{k'}|\in G\cap (G(\phi(\xi))+1)= (\phi(\xi))^-_k$. Thus, $(\xi)_{k'}(|(\beta)^-_{k'}|)\in (\xi)_k$.\\
\end{enumerate}
Now we prove (b) implies (a). To achieve this, fix $\mathbb{I}\in S$ which satisfies conditions (1), (2) and (3) of Definition \ref{definitioncheckmark} applied to $\beta,\xi,k, l$. That is:\begin{enumerate}[label=(\Roman*)]
    \item $|(\xi)^-_k|=z_\mathbb{I},$
    \item $\pi^{k,l}_\xi(\mathbb{I})\in Bl(0,l),$
    \item $(\xi)_l(|(\beta)^-_l|)\in (\xi)_k$. So in particular, $(\xi)_l(|(\beta)^-_l|)\in (\xi)_{k'}.$
    \end{enumerate} Since  $\alpha\in (\beta)_{k'}\subseteq (\beta)_l$ then $|(\alpha)_l|\leq |(\beta)_l|$. But by hypothesis we have that $|(\beta)_l|\leq |(\xi)_l|$. In this way, $|(\alpha)_l|\leq |(\xi)_l|$.  Hence, it suffices to find a $G\in \mathcal{F}_k(m_{k'})$ for which $(G,\mathbb{I})\in \mathcal{Q}$ and $\mathbb{J}_{(G,\mathbb{I})}$ satisfies the conditions (1), (2) and (3) of 
    Definition \ref{definitioncheckmark} when applied to $\alpha,\xi,k'$ and $l$. 
For this, just take a fixed $G\in \mathcal{F}_k(m_{k'})$ for which $\phi(\xi)\in G$. We claim this $G$ works. First we prove $(G,\mathbb{I})$ is an element of $\mathcal{Q}$. Observe that since $(\xi)_l(|(\beta)^-_l|)\in (\xi)_k$ then $(\xi)_l(|(\beta)^-_l|)=(\xi)_{k'}(|(\beta)^-_{k'}|)$, which means that $|(\beta)^-_{k'}|\in (|(\xi)^-_{k'}|)_{k} $. But $|(\xi)^-_{k'}|=\phi(\xi)\in G$. In this way, $(|(\xi)^-_{k'}|)_{k}=G\cap(\phi(\xi)+1)$ 
    and consequently $|(\beta)^-_{k'}|\in G$. To prove that $G[\mathbb{I}(i)]$ is an interval in $m_{k'}$ for each $i<|\mathbb{I}|$, we use (II) and apply a similar analysis. The details 
    of this are left for the reader. Now, we show the satisfaction of the remaining requirments:\begin{enumerate}
 \item Note that $z_{(G,\mathbb{I})}=G(z_\mathbb{I})=G(|(\xi)^-_{k}|)$. Since $\phi(\xi)\in G$ then $G\cap \phi(\xi)=(\phi(\xi))^-_k$. 
 On the other side, $G\cap \phi(\xi)=(\phi(\xi))_k^-$ because $\phi(\xi)\in G$ and $G\in \mathcal{F}_k(m_{k'})$. This means that $|G\cap \phi(\xi)|=|(\phi(\xi))^-_k|=|(\xi)^-_k|$. That is, there are exactly $|(\xi)^-_k|$ elements of $G$ below $\phi(\xi)$. 
 Thus, $G(|(\xi)^-_{k}|)=\phi(\xi)=|(\xi)^-_{k'}|$. In this way, $z_{(G,\mathbb{I})}=|(\xi)^-_{k'}|.$

\item Now we prove that $\pi^{k',l}_\xi(\mathbb{J}_{(G,\mathbb{I})})\in Bl(0,l)$. By the hypotheses, $(\beta)_{k'}[\,(|(\alpha)^-_{k'}|,|(\beta)^-_{k'}|]\,]$ is an interval in $(\beta)_l$. As $(\xi)_l(|(\beta)^-_l|)\in (\xi)_{k'}$, it is also true that $(\xi)_{k'}
[\mathbb{J}_{(G,\mathbb{I})}(0)]=(\xi)_{k'}[\,(|(\alpha)^-_{k'}|,|(\beta)^-_{k'}|]\,]$ is also an interval in $(\xi)_l$. Furthermore, by condition (II) we have that $(\xi)_{k'}[\mathbb{J}_{(G,\mathbb{I})}(i+1)]=(\xi)_{k'}[G[\mathbb{I}(i)]]=(\xi)_{k}[\mathbb{I}(i)]$ is also an interval in $(\xi)_l$ for each $i<|\mathbb{I}|$. Hence, 
$\pi^{k',l}_\xi(\mathbb{J}_{(G,\mathbb{I})})\in Bl(0,l)$ as we wanted. 

\item We finish by proving the last condition of Definition \ref{definitioncheckmark}. For this, just observe that since $\alpha\in (\beta)_{k'}$, then $|(\alpha)^-_l|\in (|(\beta)^-_l|)_{k'}.$ Furthermore, since $(\xi)_l(|(\beta)^-_l|)\in (\xi)_{k'}$ then $|(\beta)^-_l|\in (|(\xi)^-_l|)_{k'}$. This means that $|(\alpha)^-_l|\in (|(\xi)^-_l|)_{k'}$ which is the same as saying that $(\xi)_l(|(\alpha)^-_l|)\in (\xi)_{k'}$. 
\end{enumerate}
 \end{proof}
\end{theorem}
\begin{corollary}\label{maintheoremlemmacorollary}Let $\gamma'$ be a limit ordinal and let $\mathcal{F}'$ be a construction scheme over $\gamma'$. Also, let $k<k'\in\omega\backslash 2$, $\beta\in \gamma$, $\subseteq Bl(\beta,k)$, $T=\{(\mathbb{I},z_\mathbb{I})\}_{\mathbb{I}\in S}\in Good(\beta,k)$ and $\alpha\in (\beta)^-_{k'}$. Suppose $l\in \omega\backslash (k'+1)$ and $\delta\leq \alpha $ are such that:
\begin{enumerate}
\item $\delta$ is limit and $|(\alpha)_{l-1}\cap \delta|=r_l$,
\item $(\beta)_{k'}[\,(|(\alpha)^-_{k'}|,|(\beta)^-_{k'}|]\,]$ is an interval in $(\beta)_l$. 
 \end{enumerate}
 If $C\in[D_\delta]^{<\omega}\backslash \{\emptyset\}$ and $G\in \mathcal{F}_l$ are such that $(C,G)$ is accepted by $(k,l,\beta,\delta,T)$ then $(C,G)$ is also accepted by $(k',l,\beta,\delta,T')$ where $T'=Trans(k,k',\alpha,\beta,T)$. In particular, $j(k,l,\beta,\delta, T)\leq j(k',l,\alpha,\delta,T')$.
 \begin{proof}Since $l-1\geq k'$, $\alpha\in (\beta)_{k'}$ and $\delta\leq \alpha$, then $(\alpha)_{l-1}\cap \delta=(\beta)_{l-1}\cap \delta= R(G)$. By the hypothesis, we also have that $G$ captures $C$. Thus, the only thing left to prove is that $(\alpha,C(0),k,l-1)\checkmark T'$. For this, it is enough to show that the hypotheses of Theorem \ref{maintheoremlemma} are fulfilled for $\xi=C(0)$. Indeed, $(\beta)_{k'}[\,(|(\alpha)^-_{k'}|,|(\beta)^-_{k'}|]\,]$ is an interval in $(\beta)_l$ by the assumptions of this corollary. Furthermore, since $(\beta,C(0),k,l)\checkmark T$, then $|(\beta)_l|\leq |(C(0))_l|$ and $(C(0))_l(|(\beta)^-_l|)\in (C(0))_{k}$. As $(C(0))_k\subseteq (C(0))_{k'}$, this finishes the proof. 
     
 \end{proof}
\end{corollary}
\begin{theorem}$\mathbb{P}(\mathcal{F})\Vdash\lq\lq\, \mathcal{F}^{Gen}\textit{ satisfies }IH_2\,".$

\begin{proof}
Let $p\in \mathbb{P}(\mathcal{F})$ and $\delta<\gamma+\omega$ be a limit ordinal. Since $p\Vdash\lq\lq \,\mathcal{F}_\delta\subseteq (\mathcal{F}^\mathcal{G})_\delta\,"$, it is straightforward that $\mathcal{F}$ satisfies condition (1) of Definition \ref{definitionih2} if and only if $p$ forces the same for $\mathcal{F}^G$. Thus, we will suppose that $\mathcal{F}$ satisfies condition (2) for $\delta$ 
and find a $q\leq p$ which forces the same for $\mathcal{F}^G.$\\
Let $\beta\in \gamma+\omega$ and $k\in \omega\backslash 2$. If $\beta<\gamma$, there is nothing to do since $p\Vdash\lq\lq\,\mathcal{F}\subseteq \mathcal{F}^\mathcal{G}\,"$ and $\mathcal{F}$ already satisfies  condition $(2)$. Thus, we may suppose that $\beta\geq \gamma.$ By doing this, we can summon Lemma \ref{lemmaextensionr} to get a condition $p'\leq p$ for which $\beta \in p'$,  $p'\in \mathbb{P}_{k'}$ for some  $k'>k$ and $p\cap \gamma= r_{k'+1}$. Without loss of generality $p'=p$. In this way, $p\Vdash\lq\lq\,(\beta)_{k'}= p\cap(\beta+1)\,"$. In particular, this means $(\beta)_k$ and $Good(\beta,k)$ is fully determined by $p$. By this, we  make a small abuse of notation 
and continue the proof taking an arbitrary $\emptyset\not=S\in Bl(\beta,k)$ and $T=\{(\mathbb{I},z_\mathbb{I})\}_{\mathbb{I}\in S}\in Good(\beta,k)$. The rest of the proof is divided into two cases.\\\\
If $\delta=\gamma$, let $l=k'+1.$ The information needed to calculate  $j(k,l,\beta,\delta,T)$ lies in $\mathcal{F}|_{\delta}$ and $(\beta)_{l-1}$. Since $p$ already knows these two sets, we conclude that there is $j\in \omega$ for 
which $p\Vdash\lq\lq\,j=j(k,l,\beta,\delta,T)"$. If $j=0$, just take an arbitrary $F\in \mathcal{F}_l$ for which $R(F)=p\cap \gamma.$ Such $F$ exists because $|p\cap \gamma|=r_l.$ Now, define  $\alpha$ as $\min F_0\backslash R(F)$. Then $(\alpha)^-_{l-1}=F_0\cap \alpha=p\cap \gamma$. In this way, Lemma \ref{lemmacut} assures that $q=Cut(F,\alpha)\leq p.$ It is 
straightforward that $q\Vdash\lq\lq\beta\in q_0\backslash R(q)"$. This is just because $p=q_0$. Now we deal with the case where  $j>0$. Here, we take $F\in \mathcal{F}_l$ and $C\in [D_\delta]^j$ for which $p\Vdash\lq\lq\, (C,F)\textit{ is accepted by }(k,l,\beta,\delta,T)\,"$. By the third point in Definition \ref{acceptance}, $p\Vdash\lq\lq\, (\beta)_{l-1}\cap \gamma= R(F)$. Since $p$ also forces that $(\beta)_{l-1}=(\beta+1)\cap p$, then $p\cap\gamma=R(F)$. Thus, if we let $\alpha= 
\min F_j\backslash R(F)$ then  $(\alpha)^-_{l-1}=p\cap \gamma.$ Again by Lemma \ref{lemmacut} get that $q=Cut(F,\alpha)\leq p.$  From 
this, it is straightforward that $$q\Vdash\lq\lq\,(C,q)\textit{ is accepted by }(k,l,\beta,\delta,T)\textit{ and }\beta\in q_j\backslash R(q)".$$\\
Now we consider the case where $\delta<\gamma$. By applying Lemma \ref{lemmabetagamma} and using once more Lemma \ref{lemmaextensionr}, we get a condition whose intersection with $\gamma\backslash \delta$ is nonempty. Without loss of generality we can suppose $p\cap \gamma\backslash \delta\not=\emptyset.$  In this way, let $\alpha=\max p\cap \gamma$.  Finally, let $S'$, $\mathcal{Q}$ and $T'=Trans(\alpha,\beta,k,k',T)$ as in Definition \ref{transgood}.
\\
Since $\mathcal{F}$ satisfies $IH_2$, let us take an $l\in \omega\backslash {k'}$ which testifies the condition (2) of Definition \ref{definitionih2} for $\alpha, k'$ and $T'$. That is:
\begin{enumerate}[label=$(\alph*)$]
\item $|(\alpha)_{l-1}\cap \delta|=r_{l},$
\item There is no $C\in [D_\delta]^{<\omega}$ which is fully captured by any $F\in \mathcal{F}_l.$ Furthermore, either $j(k',l,\alpha,\delta,T')=0$ and there is $F\in \mathcal{F}_l$ such that $\alpha\in F_0\backslash R(F)$ or $j(k',l,\alpha,\delta,T')>0$ and there are $C\in [D_\delta]^{j(k',l,\alpha,\delta,T')}$ and $F\in \mathcal{F}_{l}$ for which $(C,F)$ is accepted by $(k',l,\alpha,\delta,T')$ and  $\alpha\in F_{j(k,l,\beta,\delta,T)}$.
\end{enumerate}
Let $j=j(k',l,\alpha,\delta,T')$.  As in the previous case, we should consider the subcases where $j=0$ or $j>0.$\\\\
If $j=0$, let $F\in \mathcal{F}_l$ with $\alpha\in F_0\backslash R(F)$. Note that $F\cap \delta=R(F)$. Now, take $G\in \mathcal{F}_{k'}$ for which $\alpha\in G\subseteq F$. Since $|(\alpha)^-_{k'}|<|p|=m_{k'}$ then $G\backslash (\alpha+1)\not=\emptyset.$ In this way, let $\nu=\min G\backslash(\alpha+1)$. Observe 
that $(\nu)^-_{k'}=(\alpha)_{k'}=p\cap \gamma$. Thus, if we let $q=Cut(F,\nu)$ then $q\leq p$. Indeed, $q\Vdash \lq\lq\, \beta\in q_0\backslash R(q)"$ by the part  (c) of Lemma \ref{lemmaxi} and because $q\Vdash\lq\lq\, \alpha\in q_0\backslash R(q)\,".$ We also have that $$q\Vdash\lq\lq\,(\beta)_{l-1}\cap \delta=(\alpha)_{l-1}\cap \delta= R(q)\,".$$ In this way, it has sense to consider $l$ as a candidate to fulfil condition (2) of Definition \ref{definitionih2} for $\beta,k,l$ and $T$. We claim this $q$ forces what we 
want. We need to show that $j(k,l,\beta,\delta,T)$ is forced to be $0$. We will achieve this by making use of Theorem \ref{maintheoremlemma}.  For this, note that $p\Vdash\lq\lq \,(\beta)_{k'}[\,(|(\alpha)^-_{k'}|,|(\beta)^-_{k'}|]\,]=\beta\backslash \gamma\,"$. Since $q$ forces the same and $\beta\backslash \gamma$ is already an interval in $\omega_1$, then $$q\Vdash\lq\lq\,(\beta)_{k'}[\,(|(\alpha)^-_{k'}|,|(\beta)^-_{k'}|]\,]\textit{ is an interval in }(\beta)_l\,".$$ In this way, we can use Corollary \ref{maintheoremlemmacorollary} to conclude that $q\Vdash\lq\lq\,j(k,l,\beta,\delta,T)=j=0\,"$.
As $q\Vdash\lq\lq\, q\in \mathcal{F}^{Gen}_l\textit{ and }\beta\in q_0\backslash R(F)\,"$, we are done.\\\\
Now we deal with the case where $j>0$. Here, we can take $C\in [D_\delta]^{j}$ and $F\in \mathcal{F}_l$ for which $(C,F)$ is accepted by $(k',l,\alpha,\delta,T')$ and 
$\alpha\in F_j.$ Let $\eta\in F_j$ for which $|(\eta)_{l-1}|=|(C(0))_{l-1}|$. Since  $(\alpha,C(0),k',l-1)\checkmark T'$, we get that $(\alpha,\eta,k',l-1)\checkmark T'$ as a consequence of Remark \ref{remarkcheckmark}. Thus,  there is $(G,\mathbb{I})\in \mathcal{Q}$ for which the following conditions hold:
\begin{enumerate}
\item[(0)]$|(\alpha)_{l-1}|\leq |(\eta)_{l-1}|,$
\item $G(z_\mathbb{I})=|(\eta)^-_{k'}|$,
\item $\pi^{k',l-1}_{\eta}(\mathbb{J}_{(G,\mathbb{I})})\in Bl(0,l-1),$
\item $(\eta)_{l-1}(|(\alpha)^-_{l-1}|)\in (\eta)_{k'}.$
\end{enumerate}
Since $\alpha,\eta\in F_j$ and $|(\alpha)_{l-1}|\leq |(\eta)_{l-1}|$, we can use (3) to conclude that $\alpha=(\eta)_{k'}(|(\alpha)^-_{k'}|)$. Let $\nu=(\eta)_{k'}(|(\alpha)^-_{k'}|+1)$. By definition,  $(\nu)^-_{k'}=(\alpha)_{k'}=p\cap \gamma$. Thus, we  can use  Lemma \ref{lemmacut} to conclude that $q=Cut(F,\nu)\leq p$. As in the case where $j=0$, we claim that this $q$ forces what we want.\\\\
 By condition (2) stated above we get that $(\eta)_{k'}[\mathbb{J}_{(G,\mathbb{I})}(0)]=(\eta)_{k'}[\,(|(\alpha)^-_{k'}|,|(\beta)^-_{k'}|]\,]$ is an interval in $(\eta)_{l-1}$. Since $(\eta)_{k'}[\,(|(\alpha)^-_{k'}|,|(\beta)^-_{k'}|]\,]\subseteq \eta\backslash\nu \subseteq F_j\backslash R(F)$, we conclude that this set is also an interval in $F$ whose first element is $\nu$. By definition of the Cut function, this interval is moved to the interval $(\beta+1)\backslash \gamma$.  In this way, we conclude the following:\begin{enumerate}[label=$(\alph*)$]
\item $q\Vdash \lq\lq\,|(\beta)_l|=|((\eta)_{k'}(b))_l|\," ,$
\item $q\Vdash \lq\lq\,(\beta)_l=(\nu)^-_l\cup (\beta+1\backslash \gamma)\,".$
\end{enumerate}
In particular, this implies  $q\Vdash\lq\lq\, |(\beta)_{l-1}\cap \delta|=r_l\,"$. Thus, it has sense to consider $l$ as a candidate to fulfil condition (2) of Definition \ref{definitionih2} for $\beta$, $k$ and $T.$\\
In order to finish, just note that $q$ forces $k,k',\alpha, T$ and $C(0)$ to fulfil the hypotheses of Theorem \ref{maintheoremlemma}. The following statements  are direct consequences of the application of that theorem and the Corollary \ref{maintheoremlemmacorollary}.
\begin{itemize}
\item $q\Vdash\lq\lq\,j(k,l,\beta,\delta,T)=j\,",$
\item $q\Vdash\lq\lq \, (C,q)\textit{ is accepted by }(k,l,\beta,\delta,T)\textit{ and }\beta\in q_j\," $
\end{itemize}
Thus, the proof is over.
\end{proof}
\end{theorem}
It is not hard to see that the property forced in the previous theorem can be coded by countable many dense subsets of $\mathbb{P}(\mathcal{F})$. In this way, we have the following corollaries.
\begin{theorem}\label{theoremih2}Let $\mathcal{F}$ be a construction scheme over some limit ordinal $\gamma$ which satisfies $IH_1$ and $IH_2$. Then there is a countable family $\mathcal{D}$ of dense sets in $\mathbb{P}(\mathcal{F})$ such that whenever $\mathcal{G}$ is a filter intersecting each member of $\mathcal{D}$,  then $\mathcal{F}^\mathcal{G}$ is a construction scheme over $\gamma+\omega$ satisfying $IH_1$, $IH_2$ and containing $\mathcal{F}$. 
\end{theorem}
In view of Propositions \ref{ih1unions} and \ref{ih2unions}, as well as Theorems \ref{maintheorem} and \ref{theoremih2} and the fact that $\mathcal{F}(\omega)$ satisfies both $IH_1$ and $IH_2$, we get to the Theorem we were looking for.
\begin{theorem}\label{diamondimpliesfca}$\Diamond$ implies $FCA$.
\end{theorem}
Suppose we have a partition $\mathcal{P}$ which is compatible with the type of $\mathcal{F}$. We can modify the definitions of $IH_1$ and $IH_2$ in the following way.

\begin{definition}[$IH_1(\mathcal{P})$]Let $\mathcal{F}$ be a construction scheme on some set of ordinals $X$. We say that $\mathcal{F}$ satisfies $IH_1(\mathcal{P})$ if for every $P\in \mathcal{P}$, $A\in [X]^{<\omega}$ and $\alpha\in X$, there is $F\in \mathcal{F}$ such that:
\begin{enumerate}
    \item $A\subseteq F_0,$
    \item $R(F)=F\cap \alpha$,
    \item $\rho^F\in \mathcal{P}.$
\end{enumerate}
\end{definition}
\begin{definition}[$IH_2(\mathcal{P})$] We say that $\mathcal{F}$ satisfies $IH_2(\mathcal{P})$ if for each limit $\delta<\gamma$ and $P\in \mathcal{P}$ one of the following mutually excluding conditions holds:
\begin{enumerate}[label=$(\arabic*)$]
\item There are infinitely many $l\in P$ for which there $C\in [D_\delta]^{<\omega}$ which is fully captured by some element of $\mathcal{F}_l.$\\
\item For every $\beta \in \gamma\backslash \delta$, $k\in \omega\backslash 2$ and $T\in Good(\beta,k)$ there are infinitely many $l \in P\backslash k$ for which:
\begin{enumerate}[label=$(\alph*)$]
\item $|(\beta)_{l-1}\cap \delta|=r_{l},$
\item There is no $C\in [D_\delta]^{<\omega}$ which is fully captured by any $F\in \mathcal{F}_l.$
\end{enumerate}
Furthermore, either $j(k,l,\beta,\delta,T)=0$ and there is $F\in \mathcal{F}_l$ such that $\beta\in F_0\backslash R(F)$ or $j(k,l,\beta,\delta,T)>0$ and there are $C\in [D_\delta]^{j(k,l\beta,\delta,T)}$ and $F\in \mathcal{F}_{l}$ for which $(C,F)$ is accepted by $(k,l,\beta,\delta,T)$ and  $\beta\in F_{j(k,l,\beta,\delta,T)}$.
\end{enumerate}
\end{definition}
All the results stated in this section are true when $IH_1$ is changed by $IH_1(\mathcal{P})$ and $IH_2$ is changed by $IH_2(\mathcal{P})$. Moreover, the proof are completely analogous. In this way, we conclude the following.
\begin{theorem}\label{diamondimpliesfcapart}$\Diamond$ implies $FCA(part).$
\end{theorem}
\section{Final remarks and open problems}
We include a lot of problems that the authors do not know how to solve.\\\\
We do not know how to distinguish any capturing axiom with its partition counterpart.  Moreover, the proof of Theorem \ref{diamondimpliesfcapart} is the same as Theorem \ref{diamondimpliesfca} with only very minor modifications. Thus, we ask:
\begin{Problem}\label{capartproblem}Does $CA_n$ (respectively $FCA$) imply $CA_{n}(part)$ (respectively $FCA(Part)$) for some $n\in\omega$?
\end{Problem}
There is a natural way to define a filter from a $2$-capturing construction scheme. Let $\mathcal{F}$ be a $2$-capturing construction scheme. For each $S\in [\omega_1]^{\omega_1}$ let $A_S$ be the set of all $k\in\omega$ for which there are $F\in \mathcal{F}_k$ and $C\in [S]^2$ such that $F$ captures $C$. Observe that $A_S$ is always infinite since $\mathcal{F}$ is $2$-capturing.
\begin{proposition}$\{A_S\,|\,S\in\omega_1\}$ is downwards directed with respect to $\subseteq^*$.
\begin{proof}
Let $S_0$ and $S_1$ be uncountable subsets of $\omega_1$. Recursively, we build two increasing sequences $\{\alpha_\xi\}_{\xi\in\omega_1}\subseteq S_0$ and $\{\beta_\xi\}_{\xi\in\omega_1}\subseteq S_1$ such that $\alpha_\xi<\beta_\xi<\alpha_{\xi+1}$ for every $\xi\in\omega_1$. By a refining 
argument, we can suppose that there is $k\in\omega$ for which $\rho(\alpha_\xi,\beta_\xi)=k$ for all $\xi\in \omega_1$. Furthermore, we can suppose that $|(\alpha_\xi)_k|=|(\alpha_\delta)_k$ for ever $\xi,\delta\in \omega_1.$ Let $S=\{\beta_\xi\}_{\xi\in\omega_1}$. We claim $A_S\backslash (k+1)\subseteq A_{S_0}\cap A_{S_1}$.   Let $l\in A_S\backslash (k+1)$, $\delta<\xi \in \omega_1$ and $F\in \mathcal{F}_l$ which captures $\{\beta_\delta,\beta_\xi\}$. Clearly $l\in A_{S_0}$ since $S\subseteq S_1$. In this 
way, it suffices to prove that $\{\alpha_\delta,\alpha_\xi\}$ is captured by $F.$ For this, take $\phi:F_0\longrightarrow F_1$ the increasing bijection. Then $\phi(\beta_\delta)=\beta_\xi$. Furthermore, by regularity of $\rho$ we have that $\phi[(\beta_\delta)_k]=(\beta_\xi)_k$. Since $|(\alpha_\xi)_k|=|(\alpha_\delta)_k|$, $\alpha_\delta\in (\beta_\delta)_k$ and $\alpha_xi\in (\beta_\xi)_k$, it follows that $\phi(\alpha_\delta)=\alpha_\xi.$ This concludes the proof.
\end{proof}
\end{proposition}
Let $\mathcal{U}_\mathcal{F}$ be the filter induced by $\{A_S\,|\,S\in[\omega_1]^{<\omega}\}$. Note that $\mathcal{U}_\mathcal{F}$ is an ultrafilter if and only if there is no partition $\mathcal{P}$ of $\omega$ in at least two infinite sets for which $\mathcal{F}$ is $\mathcal{P}$-2 capturing. We do not know if this filter is maximal.
\begin{Problem}Is there a $2$-capturing construction scheme, say $\mathcal{F}$, for which $\mathcal{U}_\mathcal{F}$ is an ultrafilter?
\end{Problem}
If the answer to the previous question is negative, all known results implied from $CA_2(Part)$ would also follow from $CA_2$ even if the answer to Problem \ref{capartproblem} is negative. On the other side, it would be also very interesting if $\mathcal{U}_\mathcal{F}$ could be an ultrafilter. It would be an ultrafilter on $\omega$ defined from a structure on $\omega_1$. A similar phenomenon was studied by the third author in \cite{walksselectiveultrafilters}.\\\\
We do not know if there are construction scheme in most of the usual models studied in the set theory of the reals (see \cite{Barty}). In particular we ask:
\begin{Problem} Are there $n$-capturing construction schemes in the Sacks model for some $n\in \omega\backslash 2$? What about full capturing construction schemes?
\end{Problem} 
The reader can learn more about the Sacks model in \cite{totalfailuremartin}, \cite{iteratedperfectsetforcing}, \cite{LifeintheSacksModel} and \cite{CPAbook}.
\begin{Problem}Are there $n$-capturing construction schemes in the Random model for some $n\in \omega\backslash 2$? What about full capturing construction schemes?
\end{Problem}

In \cite{forcingandconstructionschemes}, D. Kalajdzievski and F. L\'opez proved that $FCA(part)$ holds after adding $\omega_1$ Cohen reals. Still, we do not know the situation when only one Cohen real is added. 
\begin{Problem}Does adding one Cohen real also adds a fully capturing construction scheme?
\end{Problem}
We also do not know the relationship between cardinal invariants of the continuum and capturing construction schemes (see \cite{HandbookBlass} and \cite{ForcingIdealized}). We ask:
\begin{Problem}Is there any relationship between $FCA$ and cardinal invariants of the continuum?
\end{Problem}
It is not hard to see that $\sigma$-centered forcings preserve fully capturing construction schemes. Applying this result we may deduce that $\mathfrak{b}>\omega_1$ is compatible with $FCA(part)$.\\\\
In \cite{ParametrizedDiamonds}, M. D\v{z}amonja, M. Hru\v{s}\'ak and J. Moore introduced parametrized diamonds, which are weakenings of $\Diamond$ that may hold in models where $CH$ fails. Both capturing construction schemes and parametrized diamonds are useful tool for building interesting objects of size $\omega_1$. We ask:
\begin{Problem}Is there any relation between $FCA$ (respectively $CA_n$) and any parametrized $\Diamond$ principle?
\end{Problem}
It would be specially interesting to know if there is any realtionship between $FCA$ and $\Diamond(non(\mathcal{M})).$\\\\
In a similar way we ask:
\begin{Problem} What is the relationship between $FCA$(respectively $CA_n$) and $\clubsuit$?   
\end{Problem}
To learn more about $\clubsuit$, the reader may consult \cite{CardinalInvariantsoftheContinuumandCombinatoricsonUncountableCardinals}, \cite{SimilarbutnottheSame}, \cite{SticksandClubs}, \cite{Diamantesubd}, \cite{genericsplayinggames} and \cite{ProperandImproper}. Recently, A. Rinot and R. Shalev introduced the principle $\clubsuit_{AD}$. This principle follows both from $\clubsuit$ and from the existence of a Suslin Tree. Thus, it also follows from $FCA$. In this way, we may ask:
\begin{Problem}Does $CA_n(part)$ imply $\clubsuit_{AD}$ for some $n\in\omega\backslash 2$?
    
\end{Problem}
In \cite{fewsubalgebras}, M. Rubin constructed an uncountable Boolean algebra with only $\omega_1$ subalgebras. We wonder if the same can be achieved from a construction scheme.
\begin{Problem}Does $FCA(part)$ imply the existence of an uncountable Boolean algebra with only $\omega_1$ subalgebras?
\end{Problem}
The Boolean algebra $B$ constructed in \cite{fewsubalgebras} satisfies even more properties. In particular, every uncountable set of $B$ contains three distinct elements $a,b,c$ such that $a\wedge b=c.$ We remark that the algebra of clopen sets of a space constructed by the third author in \cite{schemeseparablestructures} already satisfies this property.

\section{Acknowledgment}
Part of the research was done while the authors participated in \textit{Thematic Program on Set Theoretic Methods in Algebra, Dynamics and Geometry} of the Fields Institute. The authors would also like to thank C\'esar Corral Rojas, Michael Hru\v{s}\'ak and Alejandro R\'ios Herrej\'on for some valuable comments regarding the paper. 
\bibliographystyle{plain}
\bibliography{bibliografia}

Jorge Antonio Cruz Chapital\\
Centro de Ciencias Matem\'aticas, UNAM.\\
chapi@matmor.unam.mx\\

Osvaldo Guzmán\\
    Centro de Ciencias Matemáticas, UNAM.\\
    oguzman@matmor.unam.mx\\

     Stevo Todor\v{c}evi\'c\\
    Department of Mathematics, University of Toronto, Canada\\
    stevo@math.toronto.edu\\
    Institut de Mathematiques de Jussieu, CNRS, Paris, France\\
    stevo.todorcevic@imj-prg.fr\\
    Matematički Institut, SANU, Belgrade, Serbia\\
    stevo.todorcevic@sanu.ac.rs

\end{document}